\documentclass[a4paper,reqno]{amsart}
\usepackage{amssymb,amsmath}
\usepackage[latin1]{inputenc}

    \def\RCSid$#1 #2,v #3 #4 #5${\def\version{{\underline{Version #3}} of
        file \texttt{#2}. Date: #4}}
    \RCSid$Id: JDE.tex,v 2.11 2007/07/31 16:07:10 pomet Exp $

    \newcommand{\CCC}[1]{\ensuremath{\mathbf{C}^{#1}}}
    \newcommand\RR{{I\!\!R}}
    \newcommand\NN{I\!\!N}
    \newcommand\carre{\hbox{\,{\vrule width1.0ex height1.0ex}\,}}
    \newcommand\D{\mathrm{d}}
    \newcommand\Dt{\D t}
    \renewcommand\emptyset{\varnothing} 

    \newtheorem{thm}{Theorem}[section]
    \newtheorem{lem}[thm]{Lemma}
    \newtheorem{cor}[thm]{Corollary}
    \newtheorem{prop}[thm]{Proposition}
    
    \newtheorem{openq}[thm]{Open Question}
  \theoremstyle{definition}
    \newtheorem{defn}[thm]{Definition}
  \theoremstyle{remark}
    \newtheorem{rem}[thm]{Remark}

\title{On local linearization of control systems}

\author{Laurent Baratchart}
\address{INRIA, B.P. 93, 06902 Sophia Antipolis cedex, France}
\email{Laurent.Baratchart@sophia.inria.fr} 

\author{Jean-Baptiste Pomet}
\address{INRIA, B.P. 93, 06902 Sophia Antipolis cedex, France}
\email{Jean-Baptiste.Pomet@sophia.inria.fr} 

\dedicatory{\textup{INRIA, B.P. 93, 06902 Sophia Antipolis cedex, France\\
\texttt{Laurent.Baratchart@sophia.inria.fr}, 
\texttt{Jean-Baptiste.Pomet@sophia.inria.fr} }}

\keywords{Control systems, linearization, topological equivalence, Grobman-Hartman Theorem. \\
\hspace*{1.2em}January 20,~2009. To appear in \emph{Journal of Dynamical and Control Systems}} 

\begin{document}

\begin{abstract} 
We consider the problem of topological linearization of smooth
($\CCC\infty$ or $\CCC\omega$) control systems, i.e. of their local
equivalence to a linear controllable system \emph{via} point-wise
transformations on the state and the control (static feedback
transformations) that are topological but not necessarily differentiable. 
We prove that local topological linearization implies local smooth 
linearization, at generic points. 
At arbitrary points, it implies local conjugation to a linear system
{\it via} a homeomorphism that induces a smooth diffeomorphism on the
state variables, and, except at ``strongly'' singular points, this
homeomorphism can be chosen to be a smooth mapping (the inverse map
needs not be smooth). 
Deciding whether the same is true at ``strongly'' singular points is
tantamount to solve an intriguing open question in differential topology.
\end{abstract} 

\maketitle

\section{Introduction}
Throughout the paper, \emph{smooth} means of class $\CCC\infty$.

In the early works \cite{Jaku-Res80,Hunt-Su-Mey83,VdSc84scl},  nice 
necessary and sufficient conditions were obtained for a smooth 
control system $\dot{x}=f(x,u)$, with state $x\in\RR^n$ and 
control $u\in\RR^m$, to be locally smoothly linearizable, 
{\it i.e.} locally equivalent to a controllable
linear system by means of a diffeomorphic change of variables on the
state and the control.
The afore-mentioned conditions require certain distributions 
of vector fields to be
integrable, hence locally smoothly linearizable control systems are highly 
non generic among smooth control systems. Similar results
hold for real analytic 
control systems with respect to real analytic linearizability.

Consider now the \emph{topological} linearizability
of a smooth control system, namely the property that it is
locally equivalent to a controllable
linear system {\it via} a homeomorphism
on the state and the control which may \emph{not}, this time, 
be differentiable. Obviously,
smooth linearizability implies topological linearizability; 
the extend to which the converse holds will be the main 
concern of the present paper. 
We address the real analytic case in the same stroke.

In brief, our goal
is to
{\it
  describe the class of smooth control systems that are 
locally topologically linearizable, yet not smoothly locally linearizable.}
This class in nonempty~: the smooth (even real-analytic) scalar system 
\begin{equation}
  \label{eq:u3}
  \dot{x}\;=\;u^3\ \, \ \ \ u\in\RR,\;x\in\RR,
\end{equation}
gets linearized locally around $(0,0)$ by the homeomorphism 
$(x,u)\mapsto(x,u^3)$, whereas the conditions for
smooth linearizability fail at this point.
However, we observe on this example that the conjugating homeomorphism has
much more regularity than prescribed {\it a priori}:
\\1. it is a smooth (even real-analytic) 
local diffeomorphism around all points $(x,u)$ such that
$u\neq0$,
\\2. it is triangular and induces a smooth (even real-analytic)
diffeomorphism on 
the state variable ({\it i.e.} the identity map
$x\mapsto x$),
\\3. it is a smooth (even real-analytic)
map that fails to be a diffeomorphism
only because its inverse is not smooth.

Theorem~\ref{lem-C1} of the present paper states that this example essentially
depicts the general situation.
More precisely, if a smooth control system is locally topologically
linearizable at some point $(\bar{x},\bar{u})$ in the state-control space, then
\\$1'$. in a neighborhood of $(\bar{x},\bar{u})$, the system is
locally smoothly linearizable 
around each point outside a closed subset of empty interior 
(an analytic
variety of positive co-dimension in the analytic case),
\\$2'$. around $(\bar{x},\bar{u})$, there is a triangular
linearizing homeomorphism 
that induces a smooth diffeomorphism on the state variable,
\\$3'$. the above-mentioned homeomorphism is smooth (although
its inverse may not), at least 
if $\partial f/\partial u$ has constant rank around $(\bar{x},\bar{u})$
or if $\sup_{x,u}{\rm Rank}\partial f/\partial u(x,u)=m$ on every 
neighborhood of $(\bar{x},\bar{u})$. 

Similar results hold for real-analytic linearization of a real-analytic system.

A homeomorphism 
satisfying $2'$ will be called 
\emph{quasi-smooth} (see Definitions~\ref{def-quasi}, \ref{def-lindiff}),
hence our main result is that local topological linearizability
implies local quasi-smooth linearizability.
A point $(\bar{x},\bar{u})$ where the first rank 
condition in $3'$ is satisfied is 
called \emph{regular}, and at such points local smooth linearizability
is equivalent to local topological linearizability ({\it cf.} Theorem 
\ref{th-2}).
A point $(\bar{x},\bar{u})$ where none of the rank 
conditions in $3'$ are satisfied is 
called \emph{strongly singular}.
Whether the conclusion of $3'$ continues to hold
at strongly singular points raises 
an intriguing question in differential topology,
namely can one redefine the last components of a local homeomorphism
whose first few components are smooth so as to obtain a new 
homeomorphism which is smooth? The answer seems not to be known,
see the discussion in section \ref{subsec-diff}.

\subsubsection*{Motivations} They include the following.

\textit{1.}\quad
For systems without controls, \textit{i.e.} ordinary differential
  equations, local linearization around an equilibrium has generated a sizable
  literature, see Section~\ref{sec-ode} for a small sample. It tells us
  that, even for a real analytic o.d.e., linearizability
  much depends on the admissible class of transformations (formal, real
  analytic, $C^k$ or topological). For instance,
  although analytic linearization requires subtle conditions relying 
  upon a refined analysis of resonances and small divisors, the
  Grobman-Hartman theorem says nevertheless that topological linearization 
  is always
  possible at a hyperbolic equilibrium. As one might suspect
  (this is indeed shown in section \ref{sec-generic}), 
  no naive analog to the  Grobman-Hartman theorem 
  can hold for control systems because they feature a family of 
  vector fields rather than a single one.  However, it might still be expected
  that  relaxing the smoothness of the allowable transformations 
  increases the class of  linearizable control systems. It is in fact 
  hardly so:  we knew already  from
  \cite{Jaku-Res80,Hunt-Su-Mey83,VdSc84scl} that $C^1$ linearizability of a
  smooth control system implies smooth linearizability, and we prove here
  that for $C^0$ linearizability this class does not get much bigger.
  In particular, there are no subtle questions about
  resonances and one may say that the most prominent feature of a control 
  system is to be, or not to be
  linearizable, regardless of smoothness. 

\smallskip

\textit{2.}\quad
Linearizable control systems are
  systems with linear dynamics, whose nonlinear character lies in their
  input-to-state and state-to-output maps only. Such models are 
  advocated in \cite{Judi-Hja-Ben-D-L-S-Z95,Sjob-Zha-Lju-B-D-G-H-J95}
  for identification (in the discrete-time case), as their reduced 
  complexity makes them more amenable to standard techniques. It is therefore 
  natural to investigate this class, and topological equivalence is about
  the weakest possible from the point of view of identification.

\smallskip

\textit{3.}\quad 
From a control engineering point of view, it is common practice
to design locally stabilizing feedback laws for a given system based on its
linear  approximation when the latter is controllable... and to a
certain extent one believes that the latter and the former
\emph{locally} ``look alike''.
It is therefore legitimate to ask about the
relationship between them. Since no
discriminating topological invariants are known, topological conjugacy
might appear as a good candidate. The present paper shows that the
relationship is almost never that strong: topological conjugacy to the
linear approximation is almost as rare as
differential conjugacy.

Incidentally, a system whose linear approximation is not controllable
may still happen to be locally topologically linearizable, \textit{i.e.} equivalent to a 
linear controllable system (which is \emph{not} its linear approximation).
This phenomenon is clarified in section \ref{lve}.

\subsubsection*{Techniques} 
The conditions for smooth linearizability derived in
\cite{Jaku-Res80,Hunt-Su-Mey83,VdSc84scl} come up naturally
in some sense.
Indeed, to any control system, one may associate a sequence
of distributions defined via a construction using Lie brackets of
vector fields attached to the system; it turns out that the instance
of this sequence of
distributions for linear systems yields  ``constant'' --hence
integrable-- distributions that span the entire state space in 
a finite number of
steps if the system is controllable.
Since Lie brackets and integrability of distributions are preserved under
local \emph{diffeomorphisms}, this translates at once into necessary conditions 
for smooth linearizability, shown in
\cite{Jaku-Res80,Hunt-Su-Mey83,VdSc84scl} to be sufficient. 
In contrast,
homeomorphisms do not allow to pull back Lie brackets or tangent
vector fields; hence the same conditions need not be necessary for
topological linearization, and the proofs in the present paper are more intricate.
Specifically,
we have to rely upon the notion of orbits of families of smooth
vector fields rather than integral manifolds. The proof of Theorem \ref{lem-C1} uses
classical results concerning  such orbits, 
first established in \cite{Suss73}, that we recall and slightly expand in
Appendix \ref{app-orbits}.
Incidentally, the lack of a theory dealing with orbits of  $\CCC{k}$
vector fields ($k\in\NN$) is the main reason why the results of the present paper
restrict to $\CCC\infty$ or $\CCC{\omega}$ ({\it i.e.} real analytic)
control systems.

Hopefully our method can be useful to study local 
topological equivalence to 
other classes of systems than linear ones; this is not investigated
here.

\subsubsection*{Organization of the paper}
Section \ref{sec-ode} recalls classical facts on 
local linearization of ordinary differential equations.
Section \ref{sec-topequiv} introduces conjugation for \emph{control systems}
(under a homeomorphism, a diffeomorphism, etc.) and establishes
basic properties of conjugating maps. 
Section \ref{sec-linear} reviews
(topological, smooth, linear) conjugacy between \emph{linear} control systems
after \cite{Brun70,Will80}.
Section \ref{sec-topodiff} states the main result of the paper (Theorem
\ref{lem-C1}), namely that local topological linearizability implies
local quasi-smooth linearizability for smooth control systems 
(smooth meaning either $\CCC\infty$ or $\CCC\omega$), and discusses the gap 
between smooth and quasi-smooth linearizability,
including geometric characterizations thereof. 
Section \ref{sec-proofs} contains the proofs of these results;
the proof of Theorem \ref{lem-C1}, given in 
subsection \ref{sec-proofth}, relies upon
section \ref{sec-topequiv}, results from \cite{Suss73} stated in 
Appendix \ref{app-orbits}, and
technical lemmas from Appendix \ref{blode1}. 

\section{Local linearization for ordinary differential equations}
\label{sec-ode}

Consider the differential equation 
\begin{equation}
\label{eq1}
\dot x(t)=f(x(t)),
\end{equation}
where $f\in \CCC{k}(U,\RR^n)$ with $U$ an open subset of ${\RR}^n$ and 
$k\in\NN\cup\{\infty,\omega\}$, $k\geq1$.

It is well known 
(the ``flow box theorem'', see \textit{e.g.} \cite{Arno80}) that,
around each $x_0\in U$ such that
$f(x_0)\neq0$, there is a change of coordinates of class $\CCC{k}$ 
that conjugates \eqref{eq1} to the equation
$\dot{x}_1=1$, $\dot{x}_2=0$, \ldots,
$\dot{x}_n=0$. Hence all differentiable vector fields are
equivalent to each other, at points where they do not vanish, \textit{via} a
diffeomorphism having the same degree of smoothness 
(including real analyticity).

At a point $x_0\in U$ such that $f(x_0)=0$, \textit{i.e.} at an
equilibrium of the dynamical system (\ref{eq1}), 
its linear approximation is the system 
\begin{equation}
\label{eq2}
\dot x(t)=Ax(t)-Ax_0
\end{equation}
where $A=Df(x_0)$ is the derivative of $f$ at $x_0$.
The equilibrium $x_0$ is said to be \emph{hyperbolic} if
the matrix $A$ has no purely imaginary eigenvalue.

The problem of locally linearizing \eqref{eq1} is that of finding a
local homeomorphism $h:V\rightarrow W$  around $x_0$ 
mapping the trajectories of \eqref{eq1} in $V$ onto
trajectories of \eqref{eq2} in $W$ in a time-preserving manner. In other words,
if $\phi_t$ denotes  the flow of \eqref{eq1}, we should have for each $x\in V$
that 
\[h\circ \phi_t(x)=e^{At}\bigl(h(x)-h(x_0)\bigr)+h(x_0)\] 
provided that $\phi_\rho(x)\in V$ for $0\leq\rho\leq t$.
When this is the case we say that $h$ conjugates \eqref{eq1} and \eqref{eq2}, 
and we speak of topological, $\CCC k$, smooth, or analytic 
linearization depending on the
regularity of $h$ and $h^{-1}$.

Local linearization at an equilibrium is a very old issue. At the beginning of 
the twentieth century, H. Poincar\'e already identified the obstructions to the
existence of a
\emph{formal} change of variables $h$ that removes all 
the nonlinear terms when $f$ is analytic. 
These are the so-called resonances, see {\it e.g.}
\cite{Hart82,Arno80}. In fact, resonant
monomials of order $\ell$ are  obstructions to linearizing the
Taylor expansion of $f$ at order $\ell$ and consequently also obstructions
to $\CCC\ell$ linearization. However, although 
there exists a formal power series expansion for $h$ when there are no 
resonant terms, the existence of a 
\emph{convergent} power series for $h$ (analytic 
linearization) is a delicate issue. When the eigenvalues of the Jacobian 
belong to the so-called Poincar\'e domain, the absence of resonances indeed
implies analytic linearizability (the Poincar\'e theorem). 
If it is not the case, a famous theorem by
Siegel gives additional Diophantine conditions on these eigenvalues 
to the same conclusion.
These conditions are generically
satisfied in the measure-theoretic sense \cite{Arno80}.
If no
eigenvalue of the Jacobian is purely imaginary, it turns out \cite{Rous75}
that the absence of resonances is also sufficient
for smooth ($h$, $h^{-1}$ of class $\CCC\infty$)
but in general not real analytic linearization. This is still valid when  
$f$ is merely of class $\CCC\infty$.

In contrast, if one allows conjugation {\it via}
a topological but not necessarily differentiable homeomorphism,
the Grobman-Hartman theorem asserts that every ordinary
differential equation with no purely imaginary eigenvalue of the Jacobian
(hyperbolicity) can be locally linearized around an equilibrium, that is, resonances are no
longer an obstruction. A proof of this classical result can be found in 
\cite{Hart82}:
\begin{thm}[Grobman-Hartman]
\label{th1}
Under the assumption that $x_0$ is a hyperbolic equilibrium point,
system \eqref{eq1} is topologically conjugate to system \eqref{eq2} at $x_0$.
\end{thm}
In fact, it is proved in
\cite{VanS90} that the conjugating homeomorphism $h$ 
(together with its inverse $h^{-1}$) can be chosen Hölder-continuous, and
even differentiable at $x_0$ (but not in a neighborhood).
This brings additional rigidity to the mapping $h$.

The above theorem entails that the only invariant under local topological 
conjugacy,
around a hyperbolic equilibrium, is the number of eigenvalues 
with positive real
part in the Jacobian matrix, counting multiplicity. Indeed, as is well-known 
({\it cf.} \cite{Arno74}), the linear system ${\dot x}=Ax$
where $A$ has no pure imaginary eigenvalue is topologically conjugate to
${\dot x}=DX$, where $D$ is diagonal with diagonal 
entries $\pm 1$, the number of $+1$ being the number of
eigenvalues of $A$ with positive real part.

\clearpage

\section{Preliminaries on topological equivalence for control systems}
\label{sec-topequiv}

\subsection{Control systems and their solutions}
\label{sec-sys-def}

Consider two {\em control systems} where $n,m,n',m'$ are natural integers~:
\begin{eqnarray}
\label{sysnl}
        \dot{x} &=& f(x,u)\ ,\ \ x\in\RR^{n}\ ,\ \ u\in\RR^{m}\ ,
\\
\label{sysnl2}
        \dot{z} &=&  g(z,v)\ ,\ \ z\in\RR^{n'}\ ,\ \ v\in\RR^{m'}\ ,
\end{eqnarray}
or expanded in coordinates~:
\begin{displaymath}
  \begin{array}{rcllrcl}
\dot{x}_1 &\!\!=\!\!& f_1(x_1,\ldots,x_n,u_1,\ldots,u_m)
&\ \ \ \ &
        \dot{z}_1 &\!\!=\!\!&  g_1(z_1,\ldots,z_{n'},v_1,\ldots,v_{m'})
\\
&\vdots&&&&\vdots
\\
\dot{x}_{n} &\!\!=\!\!& f_{n}(x_1,\ldots,x_n,u_1,\ldots,u_m)
&\ \ \ \ &
        \dot{z}_{n'} &\!\!=\!\!&  g_{n'}(z_1,\ldots,z_{n'},v_1,\ldots,v_{m'})\;
  \end{array}
\end{displaymath}
where $x$ or $z$ is called the \emph{state} and 
$u$ or $v$ the \emph{control}. 

Although our main
results are stated (in section \ref{sec-topodiff}) for infinitely
differentiable ---or real analytic--- control systems, their
proofs deal with non-smooth objects because the transformations 
we consider are only assumed to
be continuous. 
This leads us to keep smoothness assumptions to a minimum in 
the present section. Accordingly, 
the maps $f_i:\RR^n\times\RR^m\to\RR$ and
$g_i:\RR^{n'}\times\RR^{m'}\to\RR$ are assumed to be at least continuous; {\em
  any additional regularity assumption will be stated explicitly}. 
We do not restrict their domains of definition; this is no real loss of 
generality
because they could anyway be extended using partitions of unity (real
analyticity plays no role in the present section), and whenever
a result is stated, the domain where it holds true is precisely stated and the
value of $f$ and $g$ outside this domain does not matter.

If $m$ is zero or $f$ does not depend on $u$, equation \eqref{sysnl} reduces
to the ordinary differential equation \eqref{eq1}. 
Of course ``genuine'' control systems are those
whose right hand side does depend on the control. 

\begin{defn}
\label{def-sol}
  By \emph{a solution of \eqref{sysnl} that remains in an open set
    $\Omega\subset\RR^{n+m}$}, we mean a mapping $\gamma$ defined on a real
   interval $I$, say
\begin{equation}
\label{gamma}
\begin{array}{cccl}
        \gamma\,: & I & \rightarrow & \Omega  \\
         & t & \mapsto & \gamma(t)\ =\ 
         (\,\gamma_{\mathrm{I}}(t)\,,\,\gamma_{\mathrm{I\!I}}(t)\,)\ 
\end{array}
\end{equation}
with $\gamma_{\mathrm{I}}(t)\in\RR^n$ and $\gamma_{\mathrm{I\!I}}(t)\in\RR^m$,
such that :
\begin{itemize}
\item $\gamma$ is measurable, locally bounded, and $\gamma_{\mathrm{I}}$ is
  absolutely continuous,
\item whenever $[T_{1},T_{2}]\subset I$, we have~:
\begin{equation}
\label{eqdi}
        \gamma_{\mathrm{I}}(T_{2})\,-\,\gamma_{\mathrm{I}}(T_{1})\ \;=\;\ 
        \int_{T_{1}}^{T_{2}}
f(\,\gamma_{\mathrm{I}}(t)\,,\,\gamma_{\mathrm{I\!I}}(t)\,) \,\Dt \ .
\end{equation} 
\end{itemize}

Solutions of \eqref{sysnl2} that remain in
$\Omega'\subset\RR^{n'+m'}$ are likewise defined to be mappings
\begin{equation}
\label{gammaprime}
\begin{array}{cccl}
        \gamma'\,: & I & \rightarrow & \Omega'  \\
         & t & \mapsto & \gamma'(t)\ =\ 
         (\,\gamma'_{\mathrm{I}}(t)\,,\,\gamma'_{\mathrm{I\!I}}(t)\,)
\end{array}
\end{equation}
having the corresponding properties with respect to $g$.
\end{defn}

If $(\bar{x},\bar{u})$ is a point  in $\Omega$, 
$\mathcal{U}$ a neighborhood of $\bar{u}$ such that
$\{\bar{x}\}\times\mathcal{U}\subset\Omega$, $J$ a real interval, and
$\gamma_{\mathrm{I\!I}}:J\to\mathcal{U}$ a measurable and locally bounded map,
then, by \cite[Ch. 2, Theorem 1.1]{Codd-Lev55} and the continuity of $f$,
there exists, on a possibly smaller interval $I\subset J$, a solution $\gamma$
of \eqref{sysnl} that remains in $\Omega$
subject to the initial condition 
$\gamma_{\mathrm{I}}(0)=\bar{x}$.
This solution may not be unique without further assumptions on $f$, for 
instance that it is continuously differentiable, or merely locally Lipschitz 
in the first argument.

\begin{rem}
\label{ptw}
Observe that Definition \ref{def-sol} assigns a definite value to
$\gamma_{\mathrm{I\!I}}(t)$ for \emph{each} $t\in I$. Of course,
since $\gamma_{\mathrm{I}}$ remains a solution to \eqref{eqdi} when
the control $\gamma_{\mathrm{I\!I}}$ gets redefined over a set of measure 0,
one could identify two control functions whose values agree a.e. on $I$,
as is customary in integration theory.
However, these values are in any case subject to the constraint that
$\gamma(t)\in\Omega$ for \emph{every} $t\in I$, and altogether we find it
more convenient to adopt Definition \ref{def-sol}.
\end{rem}

\subsection{Feedbacks}
\label{sec-feed-def}

In the terminology of control, a solution in the sense of 
Definition \ref{def-sol} would be termed 
\emph{open loop} to emphasize that the value of the control at time $t$ is 
a function of time only, namely that $\gamma_{\mathrm{I\!I}}(t)$
bears no relation to the state $x$ whatsoever.
A central concept in control theory, though, is that of \emph{closed loop} or
\emph{feedback} 
control, where the value of the control at time $t$ is computed 
from the corresponding value of the state, namely is of the 
form $\alpha(x(t))$.
To make a formal definition of a feedback defined on an arbitrary open set, we 
need one more piece of notation~: 
if $\Omega\subset\RR^{n}\times\RR^{m}$ is open, we let
$\pi_{n}: \Omega\rightarrow\Omega_{\RR^{n}}$ the natural projection that
selects the first $n$ components, where
$\Omega_{\RR^{n}}=\pi_n(\Omega)\subset\RR^{n}$.

\begin{defn}
\label{def-feed}
Given an open set $\Omega\subset\RR^{n+m}$, a \emph{feedback} 
on $\Omega$ is a 
continuous mapping $\alpha:\,\Omega_{\RR^{n}}\rightarrow\RR^m$ 
such that
$(x,\alpha(x))\in\Omega$ for all $x\in\Omega_{\RR^{n}}$.
A \CCC\infty (resp. $\CCC \omega$) feedback on $\Omega$ is one
of class \CCC\infty (resp. $\CCC \omega$).
\end{defn}

A feedback is nothing but a mapping $\alpha$ such that
$x\mapsto(x,\alpha(x))$ is a continuous section of the natural fibration 
$\pi_{n}: \Omega\rightarrow\Omega_{\RR^{n}}$.
Of course, there are sets $\Omega$ whose topology prevents
the existence of any feedback.
However, if there is one there are plenty, among which
\CCC\infty feedbacks are uniformly dense. 
This is the content of the
next proposition, that will be used in the proof of Theorem \ref{lem-C1}.
To fix notations, let us agree throughout that
the symbol $\|~\|$ designates the Euclidean
norm on $\RR^\ell$ irrespectively of the positive integer $\ell$, while
$B(x,r)$ stands for the open ball centered at $x$ of radius $r$
and $\overline{B}(x,r)$ for the corresponding closed ball.
\begin{prop}
\label{prop-feed-lisse}
Let $\Omega$ be open in $\RR^{n+m}$, and 
$\alpha:\Omega_{\RR^n}\to\RR^m$ be a
feedback on $\Omega$. To each $\varepsilon>0$, there is a \CCC\infty feedback
$\beta:\Omega_{\RR^n}\to\RR^m$ such that 
$\|\alpha(x)-\beta(x)\|<\varepsilon$ for $x\in\Omega_{\RR^n}$.
\end{prop}
\begin{proof} 
Let 
$
\emptyset=\mathcal{K}_0\subset\mathcal{K}_1\cdots
\subset\mathcal{K}_k\subset\mathcal{K}_{k+1}\cdots
$
be an increasing sequence of compact subsets of $\Omega_{\RR^n}$, each of which 
contains the previous one in its interior, and whose union is all of
$\Omega_{\RR^n}$.
For each $x\in\Omega_{\RR^n}$, define an integer
\begin{equation}
  \label{kx}
  k(x)\ \;\stackrel\Delta=\;\ 
\min\{k\in\NN;\;x\in\mathcal{K}_{k}\}\ .
\end{equation}
To each $k$, by the continuity of $\alpha$ and the compactness of 
$\mathcal{K}_k$, there is $\mu_{k}>0$ such that
\begin{equation}
  \label{muK}
  x\in\mathcal{K}_k\;\Rightarrow\left\{
    \begin{array}{l}
\bullet\;
B(x,\mu_{k})\times\mathrm{Conv}\left\{\,\alpha\bigl(B(x,\mu_{k})\bigr)\,\right\}
\ \subset\ \Omega\,,
\\
\bullet\;
\forall u_1,u_2\in
\mathrm{Conv}\left\{\,\alpha\bigl(B(x,\mu_{k})\bigr)\,\right\}
,\;\|u_1-u_2\|<\varepsilon\,,
    \end{array}
\right.
\end{equation}
where the symbol $\mathrm{Conv}$ designates the convex hull.
In addition, we may assume that the sequence $(\mu_k)$ is non increasing.

Denote by
$\stackrel{\circ}{\mathcal{K}}_k$ the interior of $\mathcal{K}_k$, 
set $\mathcal{D}_k=\mathcal{K}_k\setminus \stackrel{\circ}{\mathcal{K}}_{k-1}$
for $k\geq1$, and cover the compact set $\mathcal{D}_k$
with a finite collection $\mathcal{B}_k$ of open balls having the following
properties~:
\begin{itemize}
\item each of these balls is centered at a point of  $\mathcal{D}_k$
and is contained in the open set 
$\stackrel{\circ}{\mathcal{K}}_{k+1}\setminus \mathcal{K}_{k-2}$
(with the convention that $\mathcal{K}_{-1}=\emptyset$),
\item each of these balls has radius at most 
$\displaystyle\frac{\mu_{k+1}}{2}$.
\end{itemize}
The union $\mathcal{B}=\bigcup_{k\geq1}\mathcal{B}_k$ is a countable locally finite 
collection of open balls that covers $\Omega_{\RR^n}$, and it has the
property that \emph{every ball in $\mathcal{B}$ is included in $B(x,\mu_{k(x)})$ as soon
as it contains $x$}. Let $B_j$, for $j\in\NN$, enumerate $\mathcal{B}$,
and $h_j$ be a smooth partition of unity where $h_j$ has support 
$\mathrm{supp}h_j\subset B_j$.
If we pick $x_j\in B_j$ for each $j$, the map $\beta:\Omega_{\RR^n}\to\RR^m$
defined by 
\begin{equation}
\label{equafeeds}
\beta(x)=\sum_{j\in\NN}h_j(x)\alpha(x_j)
\end{equation}
is certainly smooth. In addition, since by construction 
$x_j$ belongs to $B(x,\mu_{k(x)})$ 
whenever $h_j(x)\neq0$, we get that $\beta(x)$ lies in the convex hull
of $\alpha\bigl(B(x,r)\bigr)$ for some $r<\mu_{k(x)}$, and therefore, from 
\eqref{muK} and \eqref{kx}, that
$(x,\beta(x))\in\Omega$ and $\|\alpha(x)-\beta(x)\|<\varepsilon$.
Hence $\beta$ is a smooth feedback on $\Omega$ such that
$\|\alpha(x)-\beta(x)\|<\varepsilon$ for all $x\in\Omega_{\RR^n}$.
\end{proof}

\subsection{Conjugacy}
\label{sec-conj-def}

We turn to the notion of conjugacy for control systems, 
which is the central topic of the paper. 
\begin{defn}
\label{def-conj}
Let
\begin{equation}
\label{homeo}
\begin{array}{cccl}
    \chi\,:&\Omega&\rightarrow&\Omega'  \\
    & (x,u) & \mapsto & \chi(x,u)
    \ =\ (\,\chi_{\mathrm{I}}(x,u)\,,\,\chi_{\mathrm{I\!I}}(x,u)\,)
\end{array}
\end{equation}
be a bijective mapping between two open subsets of $\RR^{n+m}$ and
$\RR^{n'+m'}$ respectively.
We say that $\chi$ \emph{conjugates} systems
\eqref{sysnl} and \eqref{sysnl2} if,
for any real interval $I$, a map $\gamma\,:\,I\rightarrow \Omega$
is a solution of \eqref{sysnl} that remains in $\Omega$ if, and only if,
$\chi\circ\gamma$ is a solution of \eqref{sysnl2} that remains in $\Omega'$.
\end{defn}
Although this definition makes sense without any regularity assumption, we
only consider the case when $\chi$ and $\chi^{-1}$ are at least continuous.
Then Brouwer's invariance of the domain (see \emph{e.g.}
\cite{Munk84}) implies that $n'+m'=n+m$ if \eqref{sysnl} and \eqref{sysnl2}
are conjugate via such a $\chi$. 
Proposition \ref{prop-0} below asserts that 
more in fact is true.

\begin{prop}
\label{prop-0}
If the map $\chi$ in \eqref{homeo} is a homeomorphism 
that conjugates \eqref{sysnl} to
\eqref{sysnl2},
then $n=n'$, $m=m'$, and $\chi_{\mathrm{I}}$ depends only on $x$:
\begin{equation}
\label{homeotri}
\chi(x,u)\ \;=\;\ 
(\,\chi_{\mathrm{I}}(x)\,,\,\chi_{\mathrm{I\!I}}(x,u)\,)\ .
\end{equation}
Moreover, $\chi_{\mathrm{I}}:\Omega_{\RR^{n}}\to\Omega^\prime_{\RR^{n}}$ is
a homeomorphism. Here, one should recall the notation $\Omega_{\RR^n}$
that was introduced before Definition \ref{def-feed}.
\end{prop}
\begin{proof}
Let $\bar{x}$, $\bar{u}$, $\bar{u}'$ be such
that $(\bar{x},\bar{u})$ and $(\bar{x},\bar{u}')$ 
belong to $\Omega$. Let
further $x(t)$ be a solution\footnote{This solution is not necessarily unique
  since here $f$ and $g$ are merely assumed to be continuous.}
to \eqref{sysnl} with $x(0)=\bar{x}$ and
\begin{displaymath}
        \begin{array}{l}
                \ u(t)=\bar{u}~~{\rm if}~~t\leq 0,\\
                \ u(t)=\bar{u}'~~{\rm if}~~t>0~.
        \end{array}
\end{displaymath}
By conjugacy, $z(t)=\chi_{\mathrm{I}}(x(t),u(t))$ is
a solution to
\eqref{sysnl2} with $v$ given by $v(t)=\chi_{\mathrm{I\!I}}(x(t),u(t))$,
for $t\in(-\epsilon,\epsilon)$ and some $\epsilon>0$.
In particular
$\chi_{\mathrm{I}}(x(t),u(t))$ is continuous in $t$ so its values at
$0^+$ and $0^-$ are equal. Hence
$\chi_{\mathrm{I}}(\bar{x},\bar{u})=\chi_{\mathrm{I}}(\bar{x},\bar{u}')$
so that
$\chi_{\mathrm{I}}:\Omega_{{\RR^n}}\rightarrow \Omega^\prime_{\RR^{n'}}$ is
well defined and continuous. Similarly, 
$\left(\chi^{-1}\right)_{\mathrm{I}}$ induces a continuous inverse
$\Omega^\prime_{\RR^{n'}}\rightarrow\Omega_{\RR^{n}}$. By invariance of the 
domain $n=n'$.
\end{proof}
In view of this proposition, we will only consider conjugacy between
systems having the same number of states and inputs.
Hence the distinction between $(n,m)$ and $(n',m')$ from now on disappears.

\begin{rem}
\label{rmk-last}In the literature, there seems to be no general agreement on 
what should be
called a solution of a control system, nor on the concept of equivalence. 
We discuss and compare some notions in use in section~\ref{sec-xeq1}.
\end{rem}

\begin{rem}
\label{rmk-conjdiff}
Taking into account the triangular structure of $\chi$ in Proposition
\ref{prop-0}, one may describe conjugacy as resulting from a change of
coordinates in the state-space (upon setting $z=\chi_{\mathrm{I}}(x)$) 
and then  
feeding the system with a function both of the state and of a new control
variable $v$ (upon setting 
$u=(\chi^{-1})_{\mathrm{I\!I}}(z,v)$), in such a way
that the correspondence $(x,u)\mapsto(z,v)$ is invertible. 
In the language of control, this is known as a \emph{static feedback
  transformation}, and two systems conjugate in the sense of 
Definition~\ref{def-eqtop} would be termed 
\emph{equivalent under static feedback}.

This notion has received considerable attention 
(see for instance  \cite{Jaku90}),
albeit only in the \emph{differentiable} case (i.e. when 
$\chi$ is a diffeomorphism).
Differentiability has the following advantage~:
when $\chi_{\mathrm{I}}$ and $(\chi_{\mathrm{I}})^{-1}$ are
differentiable, $\chi$ conjugates
systems \eqref{sysnl} and \eqref{sysnl2} on some domain if, and only if
\begin{equation}
  \label{formuleConj}
  g(\chi_{\mathrm{I}}(x),\chi_{\mathrm{I\!I}}(x,u))\ \;=\;\ 
\frac{\partial\chi_{\mathrm{I}}}{\partial x}(x)\,f(x,u)
\end{equation}
holds true on this domain. Hence one may replace
Definition~\ref{def-eqtop}, which is based on \emph{solutions} to 
\eqref{sysnl} and \eqref{sysnl2}, by the equality above expressing
the way 
in which $\chi$ transforms the \emph{equations}. Note that 
the differentiability of $\chi_{\mathrm{I\!I}}$ is not required.
\end{rem}

Various degrees of regularity for $\chi$
give rise to corresponding notions of conjugacy in Definition~\ref{def-eqtop} below.

\begin{defn}
\label{def-quasi}
For $k\in\NN\cup\{\infty,\omega\}$, $k\geq1$, a map $\chi$ as in (\ref{homeotri}) is
called a \emph{quasi-\CCC{k} diffeomorphism} if and only of it is $\CCC0$
homeomorphism and $\chi_{\mathrm{I}}$ is a \CCC{k} diffeomorphism
$\Omega_{\RR^{n}}\to\Omega^\prime_{\RR^{n}}$, \textit{i.e.}
$\chi_{\mathrm{I}}$ and ${\chi_{\mathrm{I}}}^{-1}$ are of class $\CCC{k}$.
\end{defn}

\begin{defn}
\label{def-eqtop} 
Let $k\in\NN\cup\{\infty,\omega\}$, $k\geq1$.

Systems \eqref{sysnl} and \eqref{sysnl2} are \emph{topologically}
(resp. \emph{\CCC{k}}, resp. \emph{quasi-\CCC{k}}) \emph{conjugate
  over the pair $\Omega$,$\Omega'$} if there exists a homeomorphism
(resp. $\CCC{k}$ diffeomorphism, resp. quasi-$\CCC{k}$ diffeomorphism)
$\chi:\Omega\to\Omega'$  that conjugates the two systems.

System \eqref{sysnl} is \emph{locally} topologically ($\CCC{k}$, quasi-$\CCC{k}$)
conjugate to system \eqref{sysnl2} at
$(\bar{x},\bar{u})\in\RR^{n+m}$
if\footnote{It would be more natural to
  say that system \eqref{sysnl} {\bf\boldmath at $(\bar{x},\bar{u})$}$\in\RR^{n+m}$
  is locally conjugate to system \eqref{sysnl2} {\bf\boldmath at $(\bar{x}',\bar{u}')$}$\in\RR^{n+m}$
  if the two systems are conjugate over a pair
  $\Omega$, $\Omega'$, where $\Omega$ is a neighborhood of $(\bar{x},\bar{u})$ 
  and $\Omega'$ is a neighborhood of $(\bar{x}',\bar{u}')$.
  However, prescribing
  $(\bar{x}',\bar{u}')$ would increase notational burden and add no 
  relevant information.
}  the two systems are topologically 
($\CCC{k}$, quasi-$\CCC{k}$)
conjugate over a pair $\Omega$, $\Omega'$,
where $\Omega$ is a neighborhood of $(\bar{x},\bar{u})$.
\end{defn}

\begin{rem}
\label{rmk-f}
All definitions are invariant under linear time re-parameterization, 
namely~:\hspace{1em}
\emph{if $\chi:\Omega\rightarrow\Omega'$ conjugates systems \eqref{sysnl} and \eqref{sysnl2}, then for any $\lambda\in\RR$ (if
    $\lambda<0$, this reverses time) the map
$\chi$ also conjugates the systems
  \begin{displaymath}
    \dot{x}=\lambda f(x,u)\ \ \ \ \mbox{and}\ \ \ \ \dot{z}=\lambda g(z,v)\ .
  \end{displaymath}
}Indeed, this is trivial for $\lambda=0$, otherwise,
if $t\mapsto(x(t),u(t))$ is a solution of $\dot{x}=\lambda f(x,u)$
on a time-interval
$[t_{1},t_{2}]$, and $\tilde{x}(t)$ and $\tilde{u}(t)$ denote respectively
$x(t/\lambda)$ and $u(t/\lambda)$, then $t\mapsto(\tilde{x}(t),\tilde{u}(t))$
is a solution of \eqref{sysnl} on $[\lambda t_1,\lambda t_2]$, hence $\chi$ sends
$(\tilde{x}(t),\tilde{u}(t))$ to $(\tilde{z}(t),\tilde{v}(t))$ satisfying
$\dot{\tilde{z}}(t)=g(\tilde{z}(t),\tilde{v}(t))$. Consequently, $\chi$ maps
$(x(t),u(t))$ to $(z(t),v(t))=((\tilde{z}(\lambda t),\tilde{v}(\lambda t)))$,
which is a solution of $\dot{z}=\lambda g(z,v)$.
\end{rem}

In case there is no control ({\it i. e.} $m=m'=0$) so that neither $u$ nor
$\chi_{\mathrm{I\!I}}$ appear in \eqref{homeo},
Definition~\ref{def-eqtop} coincides with the usual notion of 
local conjugacy for ordinary differential equations.

\subsection{Properties of conjugating maps}
Below we derive some technical facts about conjugacy and feedback that are
fundamental to the proof of Theorem \ref{lem-C1}, although they are not needed
to understand the result itself.

In the proof of Proposition \ref{prop-0}, we only used conjugacy 
on a very small class of
solutions, namely those corresponding to piecewise constant controls
with a single discontinuity.
This raises the question whether smaller classes of solutions than
prescribed in Definition~\ref{def-sol} are still sufficiently rich
to check for conjugacy. Under mild conditions on $f$ and $g$, 
as we will see in the forthcoming proposition, 
conjugacy essentially holds if
it is granted for a class of inputs that locally uniformly approximates
piecewise continuous functions, and this fact will be of technical
use in the proof of Lemma \ref{lem-recur2}.
To fix terminology, we agree that a function $I\to\RR^m$,
where $I$ is a real interval, is called
{\em piecewise continuous} if it is continuous except possibly
at \emph{finitely many} interior points of $I$ where it 
has limits from both sides and is either right or 
left continuous. If in addition the function is constant (resp. affine, resp. \CCC\infty)
on every open interval not containing a discontinuity point, we say that
it is \emph{piecewise constant} (resp. \emph{piecewise affine}, resp. \emph{piecewise \CCC\infty}).
\begin{prop}[Conjugacy from restricted classes of inputs]
\label{prop-uspecial}
Assume that $f$ and $g$ are continuous $\RR^n\times\RR^m\to\RR^n$
and locally Lipschitz-continuous with respect to their first
argument\footnote{This means that each $(\bar{x},\bar{u})\in\Omega$ has a
  neighborhood $\mathcal{N}$ such that $\|f(x',u)-f(x,u)\|\leq c\,\|x'-x\|$ for
  some constant $c$ whenever $(x,u)$ and $(x',u)$ lie in $\mathcal{N}$.
}.
Let $\chi:\Omega\rightarrow\Omega'$ be a homeomorphism
between two open subsets of $\RR^{n+m}$, and denote by
$\Omega_{\mathrm{I\!I}}$ and $\Omega'_{\mathrm{I\!I}}$ 
respectively the open subsets of
$\RR^m$ obtained by projecting $\Omega$ and $\Omega'$ onto the second 
factor.
Let further $\mathcal{C}$ and $\mathcal{C}'$ be collections of locally bounded 
measurable functions $\RR\to\RR^m$ whose restrictions $\mathcal{C}{|_J}$ 
and $\mathcal{C}'{|_J}$ to any compact interval $J$ contain in their 
respective closures, for the topology of uniform convergence,
the set of all piecewise continuous functions 
$J\to\Omega_{\mathrm{I\!I}}$
and $J\to\Omega'_{\mathrm{I\!I}}$ respectively.
If $\chi$ maps every solution \eqref{gamma} of \eqref{sysnl}
such that $\gamma_\mathrm{I\!I}(t)\in\mathcal{C}{|_I}$ 
to a solution of \eqref{sysnl2} while, conversely,
$\chi^{-1}$ maps every solution
\eqref{gammaprime} of \eqref{sysnl2} such that
$\gamma'_\mathrm{I\!I}(t)\in\mathcal{C}'{|_I}$ to a
solution of \eqref{sysnl},
then the restriction of $\chi$ to any relatively compact open subset
$O\subset\Omega$ conjugates systems \eqref{sysnl} and \eqref{sysnl2}
over the pair $O$, $\chi(O)$.
\end{prop}
\begin{proof}
Let us first show that
\begin{equation}
  \label{pieceCont}
\left.
  \begin{array}{l}
\textrm{for any solution}\ \gamma:I\to\Omega\ \textrm{of \eqref{sysnl}}
\textit{ such that }\gamma_{\mathrm{I\!I}}\textit{ is}\\
\textit{piecewise continuous},
\chi\circ\gamma\textrm{ is a solution to \eqref{sysnl2}.}
  \end{array}
\right\}
\end{equation}
Since the property of being a
solution is local with respect to time, we may suppose that $I$ is a compact
interval. Then, there is an open set $\mathcal{O}$ and a compact set 
$\mathcal{K}$ such that $\gamma(I)\subset \mathcal{O}\subset
\mathcal{K}\subset\Omega$. 
By the hypothesis on $\mathcal{C}$, there exists a sequence 
of functions $\gamma_{\mathrm{I\!I},k}:I\to\RR^m$ converging
uniformly to $\gamma_{\mathrm{I\!I}}$ such that 
$\gamma_{\mathrm{I\!I},k}\in\mathcal{C}{|_I}$.
Define for each $k\in\NN$ a time-varying vector field $X^k$ by
$X^k(t,x)=f(x,\gamma_{\mathrm{I\!I},k}(t))$.
By the continuity of $f$, this sequence converges 
uniformly on compact subsets of $I\times\RR^n$ to 
$X(t,x)=f(x,\gamma_{\mathrm{I\!I}}(t))$; moreover,
since $\gamma_{\mathrm{I\!I}}$ is bounded (being piecewise continuous)
$\gamma_{\mathrm{I\!I},k}$ is also bounded,
thus the local Lipschitz character of $f(x,u)$ with respect to $x$ implies
by compactness that $X(t,x)$ and $X^k(t,x)$ are themselves locally 
Lipschitz with respect to $x$ on $I\times\mathcal{O}_{\RR^n}$.
Pick $t_0\in I$ and apply Lemma~\ref{lem-ode1} with
$I=[t_1,t_2]$, $x_0=\gamma_{\mathrm{I}}(t_0)$, and 
$\mathcal{U}=\mathcal{O}_{\RR^n}$.
This yields, say for $k>K$, that the solution $\gamma_{\mathrm{I},k}$
to the Cauchy problem 
\[
\dot{\gamma}_{\mathrm{I},k}(t)=X^k(t,\gamma_{\mathrm{I},k}(t))\,,\ \ \ \ \ \ 
~~~~~~~~~~~~\gamma_{\mathrm{I},k}(t_0)=\gamma_{\mathrm{I}}(t_0)\,,\]
maps $I$ into $\mathcal{O}_{\RR^n}$ 
and that the sequence
$(\gamma_{\mathrm{I},k})_{k>K}$ converges uniformly on $I$ to 
$\gamma_{\mathrm{I}}$.
Hence, if we let
\begin{displaymath}
  \gamma_k(t)=(\gamma_{\mathrm{I},k}(t),\gamma_{\mathrm{I\!I},k}(t)),
\end{displaymath}
the sequence $(\gamma_k)_{k>K}$ converges to $\gamma$, uniformly on $I$.
In particular $\gamma_k(I)\subset\mathcal{K}\subset\Omega$ for
$k$ large enough.

Now, since $\gamma_k:I\to\Omega$ is a solution to 
\eqref{sysnl} with $\gamma_{\mathrm{I\!I},k}\in\mathcal{C}{|_I}$,
it follows from the hypothesis that $\chi\circ\gamma_k$ is a solution to
\eqref{sysnl2} that remains in $\Omega'$, \textit{i.e.} 
with the notations of \eqref{homeo} we have, for $k$
large enough,
\begin{equation}
\label{integprimep}
  \chi_{\mathrm{I}}\circ\gamma_k(t)\;-\;\chi_{\mathrm{I}}\circ\gamma_k(t_0)
\ \;=\;\ \int_{t_0}^{t}
g(\chi\circ\gamma_k(s))\,\D s,~~~~t\in I.
\end{equation}
By the continuity of $\chi$, the convergence of $\gamma_k(t)$ to $\gamma(t)$,
and the fact that $g$ remains bounded on the compact set $\chi(\mathcal{K})$,
we can apply the dominated convergence theorem to the right hand-side of
\eqref{integprimep} to obtain in the limit, as $k\rightarrow\infty$, that
\[
\chi_{\mathrm{I}}\circ\gamma(t)\;-\;\chi_{\mathrm{I}}\circ\gamma(t_0)
\ \;=\;\ \int_{t_0}^{t}
g(\chi\circ\gamma(s))\,\D s,~~~~t\in I.
\]
Thus $\chi\circ\gamma:I\to\RR^{n+m}$ is a solution to
\eqref{sysnl2} that remains in $\Omega'$, thereby
proving \eqref{pieceCont}.

The next step is to observe from \eqref{pieceCont} that,
since piecewise constant controls are in particular piecewise
continuous, the proof of Proposition \ref{prop-0} applies to show that 
$\chi:\Omega\to\Omega'$
has a triangular structure of the form \eqref{homeotri}.

With \eqref{pieceCont} and \eqref{homeotri} at our
disposal, let us now prove the proposition in its generality.
Choose an arbitrary open subset $\mathcal{O}$ with compact closure 
$\overline{\mathcal{O}}$ in $\Omega$,
and fix two compact subsets $\mathcal{K}$ and $\mathcal{K}_1$ of $\Omega$ 
such that
\[\mathcal{O}\ \;\subset\;\ \overline{\mathcal{O}}\ \;\subset\;\ \stackrel{\circ}{\mathcal{K}}
\ \;\subset\;\ \mathcal{K}\ \;\subset\;\ \stackrel{\circ}{\mathcal{K}}_1
\ \;\subset\;\ \mathcal{K}_1\ \;\subset\;\ \Omega.
\]
where $\stackrel{\circ}{\mathcal{K}}$ stands for the {\em interior} of
${\mathcal{K}}$. 

Let $\gamma:I\to\mathcal{O}$ be a solution of \eqref{sysnl}.
We need to prove that $\chi\circ\gamma$ is a solution to
\eqref{sysnl2} and again, since the property of being a
solution is local with respect to time, we may suppose that $I$ is compact.
Notations being as in \eqref{gamma}, it follows by definition
of a solution that $\gamma_{\mathrm{I\!I}}$ is a bounded 
measurable function $I\to\RR^m$.
We shall proceed as before in that we again approximate $\gamma$ 
by a sequence $\gamma_k$ of trajectories 
of \eqref{sysnl} that are mapped by $\chi$ to trajectories of \eqref{sysnl2}. 
This time, however, the approximation process is slightly more delicate,
because it is no longer granted by the hypothesis on $\mathcal{C}$
but it will rather depend on general point-wise approximation properties
to measurable functions by continuous ones.

By the compactness of $\mathcal{K}$, there is $\varepsilon_\mathcal{K}>0$ such that
\begin{equation}
\label{compOmega}
(x,u)\in\mathcal{K}\ \Rightarrow\ B\bigl((x,u)\,,\,\varepsilon_\mathcal{K}\bigr)
\subset \stackrel{\circ}{\mathcal{K}}_1.
\end{equation}
Let $u_{\gamma_{\mathrm{I}}}:I\to\RR^m$ be an auxiliary function with the
following properties~:
\begin{itemize}
\item[(i)] $u_{\gamma_{\mathrm{I}}}$ is piecewise constant on $I$, 
\item[(ii)] $(\xi(t),u_{\gamma_{\mathrm{I}}}(t))\in
\stackrel{\circ}{\mathcal{K}}_1$ 
for all $t\in I$ and every map $\xi:I\to\RR^n$ that satisfies 
\begin{equation}
\label{epsilonu}
\sup_{t\in I}\|\xi(t)-\gamma_{\mathrm{I}}(t)\|<\varepsilon_{\mathcal{K}}/2.
\end{equation}
\end{itemize}
Such a function $u_{\gamma_{\mathrm{I}}}$ certainly exists. 
Indeed, by definition of a solution, 
$\gamma_{{\mathrm I}}$ is absolutely continuous thus {\it a fortiori}
continuous $I\to\RR^n$, and therefore we know for each $t\in I$ that the set
\[
\gamma_{\mathrm{I}}^{-1}\bigl(B(\gamma_{\mathrm{I}}(t),\varepsilon_{\mathcal{K}}/2)\bigr)
\]
is an open neighborhood of $t$ in $I$, hence a disjoint union of 
open intervals in $I$ one of which contains $t$; call this particular
interval $U_t$. By the compactness of $I$, 
we may cover the latter with finitely many intervals $U_{t_j}$
for $1\leq j\leq\nu$.
Let now $j(t)$ denote, for each $t\in I$, the smallest index 
$j\in\{1,\ldots,\nu\}$ such that $t\in U_{t_j}$.
Then, the map
\[u_{\gamma_{\mathrm{I}}}(t)=\gamma_{\mathrm{I\!I}}(t_{j(t)})\] 
clearly satisfies (i), and 
since $(\gamma_{\mathrm{I}}(t_{j(t)}),
\gamma_{\mathrm{I\!I}}(t_{j(t)}))\in\mathcal{O}\subset
\mathcal{K}$, it follows from \eqref{compOmega}
and the fact that 
$\|\gamma_{\mathrm{I}}(t)-\gamma_{\mathrm{I}}(t_{j(t)})\|
<\varepsilon_{\mathcal{K}}/2$ by 
definition of $j(t)$ that $u_{\gamma_{\mathrm{I}}}$ also satisfies (ii).

Next, recall that $\gamma_{\mathrm{I\!I}}$ is a bounded
measurable function $I\to\RR^m$
so, by Lusin's theorem \cite[Theorem 2.23]{Rudi75} applied component-wise,
there is, for every integer $k\geq1$, a continuous function
$h_k:I\to\RR^m$ that coincides with 
$\gamma_{\mathrm{I\!I}}$ outside some set $\mathcal{T}_k\subset I$ of Lebesgue
measure strictly less than $1/k^2$, and in addition such that
\begin{equation}
\label{bornh}
\sup_{t\in I}\|h_k(t)\|\leq \sqrt{m}
\sup_{t\in I}\|\gamma_{\mathrm{I\!I}}(t)\|.
\end{equation}
Put
\(E_k=\{t\in I;~(\gamma_{\mathrm{I}}(t),h_k(t))\notin
\stackrel{\circ}{\mathcal{K}}\}\).
Since $h_k$ is continuous $E_k$ is compact, and since
$\gamma(I)\subset \mathcal{O}\subset \stackrel{\circ}{\mathcal{K}}$
it is clear that $E_k\subset\mathcal{T}_k$ hence $E_k$ has Lebesgue
measure strictly less than $1/k^2$. Consequently, by the outer regularity
of Lebesgue measure, $E_k$ can be covered by finitely many open real
intervals $I_{k,1},\ldots,I_{k,N_k}$ whose lengths add up to no 
more than $1/k^2$.

We now define the sequence of functions 
$\gamma_{\mathrm{I\!I},k}$ on $I$ by setting, for $k\geq1$,
\begin{equation}
\label{defgammaIIk}
    \begin{array}{rll}
    \gamma_{\mathrm{I\!I},k}(t)&=&  h_k(t)\ \mbox{if}\ 
t\in I\setminus \bigcup_{j=1}^{N_k} I_{k,j},\\
    \gamma_{\mathrm{I\!I},k}(t)&=&  u_{\gamma_{\mathrm{I}}}(t)\ \mbox{if}\ 
t\in \bigcup_{j=1}^{N_k} I_{k,j}.\\
    \end{array}
\end{equation}
By construction $\gamma_{\mathrm{I\!I},k}$ is piecewise continuous,
and uniformly bounded independently of $k$ in view of \eqref{bornh}
and the fact that $u_{\gamma_{\mathrm{I}}}$, being piecewise constant,
is bounded. Moreover,
as $\sum_{k\geq1} 1/k^2<\infty$, the 
measure of the set $\cup_{j=1}^{N_k} I_{k,j}$ is the general term, 
indexed by $k$, of
a convergent series, hence almost every $t\in I$ belongs at most 
to finitely many of 
these sets so that $\gamma_{\mathrm{I\!I},k}$ converges
point-wise a.e. to $\gamma_{\mathrm{I\!I}}$ on $I$ as $k\to\infty$.

Redefine now $X^k(t,x)=f(x,\gamma_{\mathrm{I\!I},k}(t))$,
$X(t,x)=f(x,\gamma_{\mathrm{I\!I}}(t))$, and observe from what we just said
and the continuity of $f$ that $X^k(t,x)$ converges to $X(t,x)$ when
$k\to\infty$, locally uniformly with respect to $x\in\mathcal{O}_{\RR^n}$, 
as soon as $t\notin E$ where $E\subset I$
is a set of zero measure which is independent of $k$.
Moreover, again from the
boundedness of $\gamma_{\mathrm{I\!I},k}$, $\gamma_{\mathrm{I\!I}}$ and the
local Lipschitz character of $f$, we have that $X^k(t,x)$, $X(t,x)$ 
are locally Lipschitz with respect to $x$. 
Pick $t_0\in I$ and
apply Lemma~\ref{lem-ode1} with $\mathcal{U}=\mathcal{O}_{\RR^n}$,
$I=[t_1,t_2]$, and 
$x_0=\gamma_{\mathrm{I}}(t_0)$. We get, say for $k>K$, 
that the solution $\gamma_{\mathrm{I},k}$
to the Cauchy problem 
\[\dot{\gamma}_{\mathrm{I},k}(t)=
X^k(t,\gamma_{\mathrm{I},k}(t))\,,~~~~~~\ \ \ \ 
\gamma_{\mathrm{I},k}(t_0)=\gamma_{\mathrm{I}}(t_0)\,,\]
is defined over $I$, maps the latter into $\mathcal{O}_{\RR^n}$,
and that the sequence
$(\gamma_{\mathrm{I},k})_{k>K}$ converges uniformly on $[t_1,t_2]$ 
to $\gamma_{\mathrm{I}}$.

{\em We claim that
$\gamma_k(t)=(\gamma_{\mathrm{I},k}(t),
\gamma_{\mathrm{I\!I},k}(t))$ lies 
in $\stackrel{\circ}{\mathcal{K}}_1$ for all 
$t\in I$ when $k$ is so large that}
\begin{equation}
\label{verifpc}
\sup_{t\in I}\|\gamma_{\mathrm {I},k}(t)-
\gamma_{\mathrm{I}}(t)\|<\varepsilon_\mathcal{K}/2.
\end{equation}
Indeed,  if $t\in \cup_j I_{k,j}$, this 
follows automatically from definition \eqref{defgammaIIk}
by property (ii) of $u_{\gamma_{\mathrm{I}}}$;
if $t\notin \cup_j I_{k,j}$,
then $(\gamma_{\mathrm{I}}(t),h_k(t))
\in\stackrel{\circ}{\mathcal{K}}$
by the very definition of $\cup_j I_{k,j}$,
and since $\gamma_k(t)=(\gamma_{\mathrm{I},k}(t),h_k(t))$ 
in this case, we deduce from \eqref{compOmega} and
\eqref{verifpc} that 
$\gamma_k(t)\in\stackrel{\circ}{\mathcal{K}}_1$.
{\em This proves the claim}.

Altogether, we have shown that $\gamma_k:I\to\stackrel{\circ}{\mathcal{K}}_1$ 
is a solution of \eqref{sysnl}
as soon as $k$ is large enough, with $\gamma_{\mathrm{I\!I},k}$ a
piecewise continuous function on $I$ by construction.
By \eqref{pieceCont}, we now deduce that,
for $k$ large enough,
$\gamma'_k=\chi\circ\gamma_k$
is a solution of \eqref{sysnl2} that stays in $\Omega'$.
Let us block-decompose $\gamma'_k$ into
\begin{displaymath}
  \gamma'_{\mathrm{I},k}(t)\ =\ 
\chi_{\mathrm{I}}(\,\gamma_{\mathrm{I},k}(t)\,)
\ ,\ \ \ 
\gamma'_{\mathrm{I\!I},k}(t)\ =\ 
\chi_{\mathrm{I\!I}}(\,\gamma_{\mathrm{I},k}(t)\,,
\,\gamma_{\mathrm{I\!I},k}(t)\,)\ ,
\end{displaymath}
where we have taken into account the triangular structure of $\chi$.
That $\gamma'_k:I\to\Omega'$
is a solution of \eqref{sysnl2} means exactly that
\begin{equation}
\label{integprime}
  \gamma'_{\mathrm{I},k}(t)\;-\;\gamma'_{\mathrm{I},k}(t_0)
\ \;=\;\ \int_{t_0}^{t}
g(\,\gamma'_{\mathrm{I},k}(s)\,,\,\gamma'_{\mathrm{I\!I},k}(s)\,)\D s,~~~~t\in 
I.
\end{equation}
Due to the continuity of $\chi$, the functions $\gamma'_{\mathrm{I},k}$ and
$\gamma'_{\mathrm{I\!I},k}$ respectively converge uniformly and
point-wise almost everywhere
to $\gamma'_{\mathrm{I}}=\chi_{\mathrm{I}}\circ\gamma_{\mathrm{I}}$ and
$\gamma'_{\mathrm{I\!I}}=\chi_{\mathrm{I\!I}}\circ\gamma$ on $I$.
Since $g$ is bounded on the compact set $\chi(\mathcal{K}_1)$ that contains
$\gamma_k(I)$ for $k$ large enough, we get on the one hand, 
by dominated convergence, that the right-hand side of
\eqref{integprime} converges,
as $k\rightarrow\infty$, to
$\int_{t_0}^{t}
g(\,\gamma'_{\mathrm{I}}(s)\,,\,\gamma'_{\mathrm{I\!I}}(s)\,)\D s$,
and on the other hand that the left-hand side converges to
$\gamma'_{\mathrm{I}}(t)-\gamma'_{\mathrm{I}}(t_0)$.
Therefore 
$(\gamma'_{\mathrm{I}},\gamma'_{\mathrm{I\!I}})=
\chi_\circ\gamma:I\to\Omega'$ is a solution of \eqref{sysnl2}. 

This way we have shown that $\chi$ maps
any solution of \eqref{sysnl} that stays in a relatively compact 
open subset $\mathcal{O}$ of
$\Omega$ to a solution of \eqref{sysnl2} that stays in $\Omega'$.
This achieves the proof, for the converse is obtained symmetrically upon 
swapping $f$ and $g$, $\mathcal{C}$ and $\mathcal{C}'$,
and replacing $\chi$ by $\chi^{-1}$.
\end{proof}

The triangular structure of conjugating homeomorphisms
asserted by Proposition \ref{prop-0} is to the effect
that any such homeomorphism
$\chi:\Omega\rightarrow\Omega'$ 
is a fiber preserving map from the bundle
$\Omega\rightarrow\Omega_{\RR^{n}}$ to the bundle
$\Omega'\rightarrow\Omega'_{\RR^{n}}$.
Since feedbacks are naturally
associated to sections of these bundles
by Definition~\ref{def-feed}, $\chi$
gives rise to a natural transformation from feedbacks
on $\Omega$ to feedbacks on $\Omega'$.
This transformation will prove important enough to deserve a notation~:
to any feedback $\alpha$ on $\Omega$, we associate a feedback 
$\chi\carre\alpha$ on $\Omega'$ by the formula
\begin{equation}
\label{defcarre}
  \chi\carre\alpha(z)\ \;\stackrel\Delta=\;\ 
  \chi_{\mathrm{I\!I}}(\chi_{\mathrm{I}}^{-1}(z),
         \alpha(\chi_{\mathrm{I}}^{-1}(z)))\,.
\end{equation}
We leave it to the reader to check that the properties of an action are
satisfied, and in particular that
\begin{equation}
\label{act-inv}
  \chi^{-1}\carre\left(\chi\carre\alpha\right)\ \;=\;\ \alpha\ .
\end{equation}

Naturally associated to a control system \eqref{sysnl} and a feedback $\alpha$
is the following continuous vector field $f_{\alpha}$ on $\Omega_{\RR^{n}}$~:
\begin{equation}
  \label{falpha}
  f_{\alpha}(x) \ \;=\;\ f(x,\alpha(x))\ .
\end{equation}
If the homeomorphism $\chi$ in \eqref{homeotri} 
conjugates system \eqref{sysnl} to system \eqref{sysnl2}, then it is clear that
$\chi_{\mathrm{I}}$ maps the solutions of the ordinary differential equation
$\dot{x}=f_{\alpha}(x)$ to the solutions of the ordinary differential equation
$\dot{z}=g_{\chi\carre\alpha}(z)$. 
Indeed if $x(t)$ is a solution of the former, then
$(x(t),\alpha(x(t)))$ is a solution of the control system \eqref{sysnl} in the 
sense of Definition~\ref{def-sol} so the conjugacy assumption implies
that $(\chi_{\mathrm{I}}(x(t)),\chi_{\mathrm{I\!I}}(x(t),\alpha(x(t))))$ is a
solution of \eqref{sysnl2}, and setting 
$z(t)=\chi_{\mathrm{I}}(x(t))$ one
clearly has
$\chi_{\mathrm{I\!I}}(x(t),\alpha(x(t))))
=\chi\carre\alpha(z(t))$; hence
$z(t)$ is a solution to $\dot{z}=g_{\chi\carre\alpha}(z)$ because
$(z(t),\chi\carre\alpha(z(t)))$ is a solution of \eqref{sysnl2}.

Now, if $\alpha_1$ and $\alpha_2$ are two feedbacks on $\Omega$, 
and the two vector fields $f_{\alpha_1}$ and $f_{\alpha_2}$ are 
defined on $\Omega_{\RR^n}$ by \eqref{falpha}, we denote their 
difference by $\delta f_{\alpha_1,\alpha_2}$~:
\begin{equation}
  \label{deltafalpha}
  \delta f_{\alpha_1,\alpha_2}\ \;=\;\ 
f_{\alpha_1}\;-\;f_{\alpha_2}\ .
\end{equation}
Such vector fields are similar to the difference vector fields used in 
\cite{Jurd-Sus72},
except that we consider arbitrary feedbacks instead of constant ones.
To us, these vector fields will play an essential role.
The next proposition states that a homeomorphism that conjugates two
control systems also conjugates the integral curves of such difference vector
fields. 

\begin{prop}[preservation of difference vector fields]
\label{prop-01}
Suppose that $f$ and $g$ in \eqref{sysnl} and
\eqref{sysnl2} are continuous and locally Lipschitz continuous
with respect to their first argument.
Assume they
are locally topologically conjugate at $(0,0)$ over the pair 
$\Omega,\Omega'$.
Then, notations for $\chi_{\mathrm{I}}$ and $\chi_{\mathrm{I\!I}}$ being as in
Proposition~\ref{prop-0}, we have
for every pair of feedbacks $\alpha_1,\alpha_2$ on $\Omega$ 
that $\chi_{\mathrm{I}}$ conjugates any solution of
  \begin{equation}
\label{edof}
\dot{x}\ \;=\;\ \delta f_{\alpha_1,\alpha_2}(x)
  \end{equation}
that remains in $\Omega_{\RR^n}$ to a solution of
  \begin{equation}
\label{edog}
\dot{z}\ \;=\;\ \delta g_{\chi\carre\alpha_1,\chi\carre\alpha_2}(z)
  \end{equation}
that remains in $\Omega'_{\RR^n}$.
\end{prop}
It is perhaps worth emphasizing that the solutions of \eqref{edof} and
\eqref{edog} need not be unique since $\alpha$ is merely assumed to be
continuous.

\begin{proof}
Let
$\eta:[t_1,t_2]\rightarrow\Omega_{\RR^n}$ be an integral curve of $\delta
f_{\alpha_1,\alpha_2}$, and set 
\begin{equation}
\label{open}
u_1(t)\ \;=\;\ \alpha_{1}(\eta(t))\  \ ,\ \ \ \ 
u_2(t)\ \;=\;\ \alpha_{2}(\eta(t))\ .
\end{equation}
Let further $\widehat{f}:\RR^{n+m}\rightarrow\RR^n$ be bounded, continuous and
Lipschitz continuous with respect to its first argument,
and coincide with $f$ on some compact neighborhood of
\begin{displaymath}
  \eta([t_1,t_2])\times\left(
\alpha_{1}(\eta([t_1,t_2]))\,\bigcup\,\alpha_{2}(\eta([t_1,t_2]))
\right)\ .
\end{displaymath}
Such a  $\widehat{f}$ is easily obtained upon multiplying $f$ by a function of
class \CCC\infty with compact support.
For $\ell\in\NN$, let $\eta^\ell$ be the solution to the Cauchy problem
\begin{equation}
    \label{eq:a12}
        \eta^\ell(t)\ \;=\;\ \eta(t_1)\;+\;
        \int_{t_{1}}^{t}G_\ell(\tau,\eta^\ell(\tau)) \D\tau\ ,
\end{equation}
with
  \begin{equation}
\label{eq:a15}
    \begin{array}{rl}
    G_\ell(t,x)\;=&  2\,\widehat{f}(x,u_1(t))\\
&\hspace{.8cm}
\mbox{if}\ 
t\in [t_1+\frac j \ell (t_2-t_1),t_1+(\frac j \ell + \frac{1}{2\ell})(t_2-t_1) ),
\\[1ex]
    G_\ell(t,x)\;=&-2\,\widehat{f}(x,u_2(t))\\
&\hspace{.8cm}
\mbox{if}\ 
t\in [t_1+(\frac j \ell + \frac{1}{2\ell})(t_2-t_1),t_1+\frac {j+1} \ell (t_2-t_1) ),
\\[1ex]
    G_\ell(t_2,x)\;=& -2\,\widehat{f}(x,u_2(t_2)),~~~~\ \ \ \ 
0\leq j\leq\ell-1.
    \end{array}
  \end{equation}
The definition of $\eta^\ell$ is valid because,
since $G_\ell(t,x)$ is bounded and locally Lipschitz with respect 
to the variable $x$, the solution to \eqref{eq:a12} uniquely exists.

From Lemma~\ref{lem-aver} applied to the case where
$X^{1,\ell}(t,x)=\widehat{f}(x,u_1(t))$ and
$X^{2,\ell}(t,x)=\widehat{f}(x,u_2(t))$
are in fact independent of $\ell$, any accumulation point of the sequence
$(\eta^\ell)$, say $\eta^\infty$, is a solution to
\begin{displaymath}
\dot{\eta}^\infty(t) \ \;=\;\ 
\widehat{f}(\eta^\infty(t),u_1(t))\;-\;\widehat{f}(\eta^\infty(t),u_2(t))
\ \ ,\ \ \ \ \eta^\infty(t_1)\;=\;\eta(t_1)\ .
\end{displaymath}
Since $\widehat{f}$ is locally Lipschitz continuous with respect 
to its first
argument, the solution to this Cauchy problem is unique and, since $f$ and
$\widehat{f}$ coincide at all points $(\eta(t),u_1(t))$ and
$(\eta(t),u_2(t))$, this entails $\eta^\infty=\eta$.
Thus $(\eta^\ell)$ converges uniformly to $\eta$ on $[t_1,t_2]$ and, for
$\ell$ large enough, $\eta^\ell$ remains a solution of \eqref{eq:a12} if
$\widehat{f}$ is replaced by $f$ in \eqref{eq:a15}. Moreover,
$\eta^\ell([t_1,t_2])\subset\Omega_{\RR^n}$ for $\ell$ large since
the same is true of $\eta$.
Since $\chi$ conjugates the two systems, hence also by
Remark~\ref{rmk-f} the systems where $f$ and $g$ are 
multiplied by $2$ or $-2$,
the map $\chi_{\mathrm{I}}\circ\eta^\ell:[t_1,t_2]
\rightarrow\Omega'_{\RR^n}$ 
is, for $\ell$ large enough, a solution to
\begin{equation}
\label{eq:a13b}
        \chi_{\mathrm{I}}\circ\eta^\ell(t)
\ \;=\;\ \chi_{\mathrm{I}}\circ\eta(t_1)
\;+\; \int_{t_{1}}^{t}
\widetilde{G}_\ell(\tau,\chi_{\mathrm{I}}\circ\eta^\ell(\tau)) \D\tau
\end{equation}
with
  \begin{equation}
    \label{eq:a13}
    \begin{array}{rl}
    \widetilde{G}_\ell(t,z)\;=& 2\,g(z,\chi_{\mathrm{I\!I}}(\chi_{\mathrm{I}}^{-1}(z),u_1(t)))\\
&\hspace{.8cm}
\mbox{if}\ 
t\in [t_1+\frac j \ell (t_2-t_1),t_1+(\frac j \ell + \frac{1}{2\ell})(t_2-t_1) ),
\\[1ex]
    \widetilde{G}_\ell(t,z)\;=& -2\,g(z,\chi_{\mathrm{I\!I}}(\chi_{\mathrm{I}}^{-1}(z),u_2(t)))\\
&\hspace{.8cm}
\mbox{if}\ 
t\in [t_1+(\frac j \ell + \frac{1}{2\ell})(t_2-t_1),t_1+\frac {j+1} \ell (t_2-t_1) ),
\\[1ex]
\widetilde{G}_\ell(t_2,z)\;=&
-2\,g(z,\chi_{\mathrm{I\!I}}(\chi_{\mathrm{I}}^{-1}(z),u_2(t_2))).
    \end{array}
  \end{equation}
Since $(\chi_{\mathrm{I}}\circ\eta^\ell)$
converges uniformly to $\chi_{\mathrm{I}}\circ\eta$ 
by the continuity of
$\chi$, replacing $g$ by a bounded and continuous
$\widehat{g}:\RR^{n+m}\rightarrow\RR^{n}$ that coincides with
$g$ on a compact neighborhood of
\begin{displaymath}
  \chi_{\mathrm{I}}\circ\eta([t_1,t_2])\times
\left(
\chi_{\mathrm{I\!I}}(\eta([t_1,t_2]),\alpha_{1}(\eta([t_1,t_2])))
\,\bigcup\,
\chi_{\mathrm{I\!I}}(\eta([t_1,t_2]),\alpha_{2}(\eta([t_1,t_2])))
\right)
\end{displaymath}
does not affect the validity of \eqref{eq:a13b}-\eqref{eq:a13} for $\ell$
large enough.
Lemma \ref{lem-aver} now implies that all accumulation points of the sequence
$(\chi_{\mathrm{I}}\circ\eta^\ell)$ in the uniform topology on $[t_1,t_2]$ are
solutions of
$$\dot{z}=g(z,\chi_{\mathrm{I\!I}}(\chi_{\mathrm{I}}^{-1}(z),u_1(t)))
\,-\,g(z,\chi_{\mathrm{I\!I}}(\chi_{\mathrm{I}}^{-1}(z),u_2(t))).$$
Because $\chi_{\mathrm{I}}\circ\eta$ is such an accumulation point, it is
by \eqref{open} a solution to
\begin{displaymath}
  \dot{z}\ \;=\;\ 
g(z,\chi_{\mathrm{I\!I}}(\chi_{\mathrm{I}}^{-1}(z),
\alpha_{1}(\chi_{\mathrm{I}}^{-1}(z)))
\,-\,g(z,\chi_{\mathrm{I\!I}}(\chi_{\mathrm{I}}^{-1}(z),
\alpha_{2}(\chi_{\mathrm{I}}^{-1}(z)))\ ,
\end{displaymath}
which is nothing but \eqref{edog}.\end{proof}

\subsection{Alternative notions of conjugacy and equivalence}
\label{sec-xeq1}

\subsubsection{Transformations in functional spaces}
\label{sec-Colo}
Following \cite{Colo-Kli93asp}, one may view the control system \eqref{sysnl} as
a flow on the product space $\RR^{n}\times\mathcal{U}$, 
where $\mathcal{U}$ is a functional space of 
admissible controls whose dynamics is induced by the time-shift.
Transformations on $\RR^{n}\times\mathcal{U}$ then naturally arise; they
involve the \emph{future and the past of the control}, 
unlike the mere homeomorphisms on finite dimensional spaces that we consider
here.
The corresponding notion of equivalence is obviously rather weak.
In \cite{Bara-Chy-Pom06}, a ``Grobman-Hartman theorem'' theorem is proved in
this setting, i.e. generic control systems \eqref{sysnl} are
locally conjugate to a linear system via this kind of transformation. 
With the much stronger notion of equivalence that we use here, 
we shall see (section~\ref{sec-generic}) that ``almost'' no
system is conjugate to a linear system.

Let us also mention \cite{Celi95}, where control systems are 
maps $(x(0),u(.))\mapsto x(.)$ that satisfy certain axioms, without reference
to differential equations, and where the notion of topological equivalence
involves transformations on the product $\RR^n\times\mathcal{U}$. 

\subsubsection{$x$-conjugacy}
\label{sec-xconj}
Let us call \emph{$x$-solution} of system \eqref{sysnl} any map
$t\mapsto\gamma_{\mathrm{I}}(t)$ \emph{such that there exists} a map
$\gamma_{\mathrm{I\!I}}$ for which
$\gamma=(\gamma_{\mathrm{I}},\gamma_{\mathrm{I\!I}})$ is a solution in the
sense of Definition \ref{def-sol};
the set of $x$-solutions is the projection on
the $x$ factor of the set of solutions.
Let then $x$-conjugacy be defined in the same way as Definition
\ref{def-eqtop} defines conjugacy, except that we replace solutions by
$x$-solutions and the homeomorphism $\chi$ that acts on state and control with
a homeomorphism $x\mapsto z=h(x)$ on the state only. 

In the literature, both notions are used (without
the prefix ``$x$-''). For instance \cite{Will80}, devoted to the topological
classification of {\em linear} control systems (see section \ref{sec-willems})
relies on $x$-conjugacy.
We favor Definitions \ref{def-conj} and \ref{def-eqtop} of conjugacy and
solutions because results have to be stated locally with
respect both to $x$ and $u$ for nonlinear control systems.

Conjugacy implies $x$-conjugacy: use Proposition \ref{prop-0},
take $h=\chi_{\mathrm{I}}$
and ignore $\chi_{\mathrm{I\!I}}$.
The converse is not true in general, as the reader may check easily. 

\section{The case of linear control systems}
\label{sec-linear}

\subsection{Kronecker indices}
\label{sec-Kroneker}

A \emph{linear} control systems is a special instance of \eqref{sysnl}, 
of the form
\begin{equation}
  \label{syslin}
  \dot x \ \;=\;\ Ax\;+\;Bu
\end{equation}
where $A$ and $B$ are constant $n\times n$ and
$n\times m$ matrices respectively.
When dealing with linear systems, it is natural to consider an equivalence
relation similar to that of Definition~\ref{def-eqtop}, but where $\chi$ is
restricted to be a linear isomorphism~:
\begin{defn}
\label{def-eqlin}
Two linear systems
\begin{displaymath}
  \dot x \ =\ Ax\;+\;Bu
\ \ \ \ \ \mbox{and}\ \ \ \ \ 
  \dot z \ =\ \widetilde{A}z\;+\;\widetilde{B}v
\end{displaymath}
are \emph{linearly conjugate} if and only if any of the following two
equivalent properties is satisfied~:
\begin{enumerate}
\item
There is a nonempty open set $\Omega\subset\RR^{n+m}$, and a linear
isomorphism $\chi$ of $\RR^{n+m}$ whose restriction
$\Omega\rightarrow\chi(\Omega)$ conjugates the two systems in the
sense of Definition~\ref{def-eqtop}.
\item
There exist matrices $P\in\RR^{n\times n}$, $Q\in\RR^{m\times m}$ and
$K\in\RR^{n\times m}$, with $P$ and $Q$ invertible, such that
\begin{equation}
  \label{lfeedABtilde}
  \begin{array}{rcl}
\widetilde{A} & = & P(A-BK)P^{-1}\ ,
\\
\widetilde{B} & = & PBQ^{-1}\ .
  \end{array}
\end{equation}
\end{enumerate}
\end{defn}
Since, by Proposition~\ref{prop-0}, a linear conjugating homeomorphism is 
necessarily of the form $(x,u)\mapsto(Px,Kx+Qu)$, the
equivalence between properties (1) and (2) follows at once
from differentiating the solutions. Provided it exists,
$\Omega$ plays absolutely no role in this context
since \eqref{lfeedABtilde} implies 
that the two systems are in fact linearly
conjugate on all of $\RR^{n+m}$.

Linear conjugacy actually defines an equivalence relation on linear control
systems or equivalently on pairs $(A,B)$, for which \eqref{lfeedABtilde} can be
read as ``$(A,B)$ is equivalent to $(\widetilde{A},\widetilde{B})$''. 
The classification of linear systems under this equivalence relation 
is well-known
\cite{Brun70}, and goes as follows.
Each equivalence class contains a pair $(A_c,B_c)$ of the form
(block matrices)~:
\begin{equation}
  \label{fBru}
  A_c\ \;=\;\ 
  \left(
    \begin{array}{cccc}
A_0^c&0&\cdots&0\\
0&A_1^c&\ddots &\vdots\\
\vdots&\ddots&\ddots&0 \\
0&\cdots&0&A_m^c
    \end{array}  \right)
\ \ ,\ \ \ 
  B_c\ \;=\;\ 
  \left(
    \begin{array}{ccc}
0&\cdots&0\\
b^c_1&\ddots &\vdots\\
0&\ddots&0 \\
\vdots&0&b^c_m
    \end{array}  \right)
\end{equation}
where
\begin{equation}
  \label{fBru2}
  A_i^c\;=\;\ 
  {\left(
    \begin{array}{cccccc}
0&1&0&\cdots&\ \ \ &0\\
\vdots &\ddots&\ddots&\ddots&&\vdots \\
\vdots& & & && 0 \\
 & & &&\!\!\!\!\!\!\!\!\ddots&1 \\
0&\cdots&&&\!\!\!\!\!\!\!\!\cdots&0
    \end{array}  \right),}_{\!\!\!(\kappa_i\times\kappa_i)}
  b^c_i\;=\;\ 
  {\left(
    \begin{array}{c}
0\\\vdots\\ \vdots\\0\\1
    \end{array}  \right),}_{\!\!\!(\kappa_i\times 1)}
\hspace{-2ex}1\leq i\leq m.
\end{equation}
The integers $(\kappa_1,\ldots,\kappa_m)$ are called
the controllability indices of
the control system, also known as the Kronecker indices of the matrix 
pencil $(A,B)$,
while $A_0^c$ is a square matrix of dimension $n-(\kappa_1+\cdots+\kappa_m)$
that may be assumed in Jordan canonical form.
Note that $\kappa_1+\cdots+\kappa_m\leq n$, and if
$\kappa_1+\cdots+\kappa_m=n$ there is no $A_0^c$~; also, it may well
happen that $\kappa_i=0$, in which case $A_i^c$ and $b_i^c$ are empty 
and do not occur in \eqref{fBru} to the effect that
there are less than $m$ blocks beyond $A_0^c$.
Normalizing so that
\begin{displaymath}
  \kappa_1\geq \cdots\geq \kappa_m\geq 0,
\end{displaymath}
and ordering the Jordan blocks arbitrarily,
there is one and only one such normal form per equivalence class.
A complete set of invariants is then the list of Kronecker indices and the
spectral invariants of the matrix $A_0^c$.

With the natural partition $z=(Z_0,Z_1,\ldots,Z_m)$ corresponding to the block
decomposition \eqref{fBru}, the control
system associated to the pair $(A_c,B_c)$ reads
\begin{displaymath}
  \dot{Z}_0\ \;=\;\ A_0Z_0\ ,\ \ 
  \dot{Z}_1\ \;=\;\ A_1Z_1\;+\;u_1b^c_1\ ,\ \ \cdots\ \ ,\ \ 
  \dot{Z}_m\ \;=\;\ A_mZ_m\;+\;u_mb^c_m\ ,
\end{displaymath}
where $Z_0$ is missing if $\kappa_1+\cdots+\kappa_m=n$ and
$Z_i$ is missing if $\kappa_i=0$. Because it is not influenced at all by the
controls, $Z_0$ is sometimes called the non-controllable part of the state.
In this paper, we are only interested in controllable linear systems, namely~:
\begin{defn}
\label{def-eqlincont}
A linear control system \eqref{syslin}
is said to be \emph{controllable} if, and only if, the following two
equivalent properties are satisfied~:
\begin{enumerate}
\item
There is no bloc $A_0^c$ in the associated
normal form \eqref{fBru}.
\item Kalman's criterion for controllability~:
\[\mathrm{Rank}(B,AB,\ldots,A^{n-1}B)=n.\]
\end{enumerate}
\end{defn}
To see the equivalence of the two properties, observe that
the $n-\kappa_1-\cdots-\kappa_m$ first rows of the matrix $P$ that puts
$(A,B)$ into canonical form (i.e. $z=Px$) form a
basis of the smallest dual subspace that annihilates the columns of
$B$ and at the same time is invariant under right 
multiplication by $A$, 
i.e. they are a basis of the left kernel of
$(B,AB,\ldots,A^{n-1}B)$.
For controllable linear systems, the only invariant under linear conjugacy
is thus the ordered list of Kronecker indices.
These can be computed from
$(B,AB, \ldots, A^{n-1}B)$ as follows~: if we put
\begin{equation}
  \label{indices}
  \begin{array}{rl}
r_j\ =\ \mathrm{Rank}(B,AB,\ldots,A^{j-1}B)\;,\ j\geq1\ ,\ \ &r_0=0,\ r_{-1}=-m\ ,
\\
s_j\ =\ r_j-r_{j-1}\;,\ j\geq1\ ,\ \ &s_0=m\ ,
  \end{array}
\end{equation}
then $s_j$ does not increase  with $j$ and a moment's
thinking will convince the reader that
the number of Kronecker indices that are equal to $i$ is $s_i-s_{i+1}$, or
equivalently that $s_{k}$ is the number of $\kappa_j$'s that are no smaller
than $k$.

To us, it will be more convenient to use as normal form the following
permutation of the previous one.
Let $\rho$ be the smallest integer such that $s_\rho=0$, so that
\begin{displaymath}
        0\ =\ s_{\rho}\ <\ s_{\rho-1}\ \leq\ s_{\rho-2}\ \leq\ \cdots\ \leq\ 
        s_{1}\ \leq\ s_{0}\ =\ m\ ,
\end{displaymath}
with $\displaystyle\sum_{j\geq 1}s_{j}=n$.
From these we define, for $0\leq i\leq\rho$~:
\begin{equation}
\label{sigmai}
        \sigma_{i}\ \;=\;\ \sum_{j\geq i}s_{j}\ \;=\;\ n-r_{i-1}\,,
\end{equation}
so that in particular $\sigma_{\rho}=0$, $\sigma_{\rho-1}=s_{\rho-1}>0$,
$\sigma_{1}=n$ and $\sigma_{0}=n+m$. 
Note that, from \eqref{indices}, $\sigma_{i}=n-r_{i-1}$ for $i\geq 1$.
We shall write our controllable canonical form as 
$\dot{z}=A_c z + B_c v$ with
\begin{equation}
\label{syslcan}
A_c \!=\!
\left(  \!\!\!\!
\begin{array}{c}
        {\setlength{\unitlength}{0.7pt}
        \begin{picture}(290,290)
        \thinlines
            \put(  0,270){\line(1,0){93}}
            \put( 20,290){\line(0,-1){58}}

            \put( 17,235){\line(1,0){126}}
            \put( 55,290){\line(0,-1){93}}

            \put( 52,200){\line(1,0){106}}
            \put( 90,273){\line(0,-1){81}}

            \put(140,238){\line(0,-1){46}}

            \put(157,130){\line(1,0){133}}
            \put(160,133){\line(0,-1){66}}

            \put(157, 70){\line(1,0){133}}
            \put(220,133){\line(0,-1){133}}

            \put( 10,280){\makebox(0,0){0}}
            \put( 35,255){\makebox(0,0){0}}
            \put( 70,220){\makebox(0,0){0}}
            \put(190,100){\makebox(0,0){0}}
            \put(255, 35){\makebox(0,0){0}}

            \put( 20,270){\makebox(35,20){$J^{s_{\rho-1}}_{s_{\rho-2}}$}}
            \put( 55,235){\makebox(35,35){$J^{s_{\rho-2}}_{s_{\rho-3}}$}}
            \put( 90,200){\makebox(50,35){$J^{s_{\rho-3}}_{s_{\rho-4}}$}}
            \put(220, 70){\makebox(70,60){$J^{s_{2}}_{s_{1}}$}}

           \multiput(100,190)(10,-10){6}{\circle*{2}}
           \multiput(150,190)(10,-10){6}{\circle*{2}}
        \end{picture}}
\end{array}
\!\!\!\!\right)
,\,
B_c \!=\!
\left(\!\!\!\!
\begin{array}{c}
        {\setlength{\unitlength}{0.7pt}
        \begin{picture}(72,290)
        \thinlines
            \put(  0,70){\line(1,0){72}}
            \put(  0,70){\makebox(72,220){0}}
            \put(  0,0){\makebox(72,70){$J^{s_{1}}_{s_{0}}$}}
        \end{picture}}
\end{array}
\!\!\!\!\right)\!\!
\end{equation}
where for any integers $r$ and $s$ with $s\leq r$, $J^{s}_{r}$ is the
$s\times r$ matrix
\begin{equation}
\label{syslcanJ}
        J^{s}_{r}\ \;=\;\ \left(
\begin{array}{c}
        {\setlength{\unitlength}{0.8pt}
        \begin{picture}(90,60)
        \thinlines
            \put(60,0){\line(0,1){60}}
            \put(0,0){\makebox(60,60){$I_{s}$}}
            \put(60,0){\makebox(30,60){0}}
        \end{picture}}
\end{array}
        \right)
\end{equation}
where $I_{s}$ is the $s\times s$ identity matrix.

\subsection{Topological classification of linear control systems}
\label{sec-willems}

In \cite{Will80}, which is devoted to the topological
classification of \emph{linear} control systems and uses the notion of
$x$-conjugacy rather than conjugacy ({\it cf}  Section~\ref{sec-xeq1}), the following result is proved:

\begin{thm}[Willems~\cite{Will80}]
\label{th-willems}
If two linear control systems $\dot{x}=Ax+Bu$ and
$\dot{z}=\widetilde{A}z+\widetilde{B}v$ are topologically $x$-conjugate,
then they have the same list of Kronecker indices, and the non-controllable
blocks $A_0^c$ and $\widetilde{A}_0^c$ in their respective canonical forms
\eqref{fBru} are such that the two linear differential equations
$\dot{X}_0=A_0^c X_0$
and $\dot{Z}_0=\widetilde{A}_0^c Z_0$ are topologically equivalent.
\end{thm}

As pointed out in Section~\ref{sec-xeq1}, topological conjugacy 
implies topological $x$-conjugacy but not conversely.
However, for linear control systems having the same number $m$ of inputs,
Theorem \ref{th-willems} implies that these notions are equivalent.
Indeed, if two systems are respectively brought into their canonical form
\eqref{fBru} by a linear change of variable on $\RR^{n+m}$, 
and if in addition they are $x$-conjugate, then
their non-controllable parts are topologically equivalent while
the remaining blocks are identical by equality of the Kronecker indices.
Hence, both in the above theorem and in the corollary below, one may use
indifferently ``$x$-conjugate'' or ``conjugate''

\begin{cor}
\label{cor-willems}
If two linear systems $\dot{x}=Ax+Bu$ and
$\dot{z}=\widetilde{A}z+\widetilde{B}v$ are topologically conjugate and one
of them is controllable, then the other one is controllable too and they are
linearly conjugate.
\end{cor}
\begin{proof}
Controllability is preserved, since Kronecker indices are by the theorem.
Linear conjugacy follows, as we saw that the list of Kronecker
indices is a complete invariant for controllable systems under 
linear conjugacy.
\end{proof}

In some sense, the results of section \ref{sec-topodiff} can be
viewed as a generalization of Corollary~\ref{cor-willems} to a local setting
where only one of the two systems is linear.

\section{Local linearization for control systems}
\label{sec-topodiff}

In this section, we consistently assume that the map $f$ defining system
\eqref{sysnl} is either smooth or real-analytic.

\medskip

\begin{defn}
\label{def-lindiff} Let $k\in\{\infty,\omega\}$.
The system \eqref{sysnl} is said to be 
\emph{locally topologically}
(resp. \emph{\CCC{k}}, resp. \emph{quasi-\CCC{k}}) 
\emph{linearizable at $(\bar{x},\bar{u})\in\RR^{n+m}$} if it
is locally topologically
(resp. \emph{\CCC{k}}, resp. \emph{quasi-\CCC{k}})
conjugate, in the sense of Definition~\ref{def-eqtop},
to a linear controllable system $\dot{z}=Az+Bv$ 
(\textit{cf.} Definition~\ref{def-eqlincont}). 
\end{defn}
This definition of smooth linearizability coincides with
linearizability by smooth static feedback as described in the
textbooks \cite{Isid89,Nijm-VdS90}. In subsection~\ref{sec-JakRes}, 
we recall classical 
necessary and sufficient geometric conditions for a system to be 
smoothly (resp. analytically) linearizable, and we complement them 
with a characterization of quasi-smooth (resp. quasi-analytic) 
linearizability. 

\subsection{Main result}
\label{subsec-diff}

If a smooth control system is locally topologically linearizable, then the 
conjugating homeomorphism has a lot more regularity than required
\textit{a priori}. This is in contrast with the Grobman-Hartman theorem for 
ODE's and constitutes the central result of the paper:
\begin{thm}
        \label{lem-C1}
Let $k\in\{\infty,\omega\}$ and 
assume that $f$ is of class \CCC{k} on an open set
$\Omega\subset\RR^{n+m}$.
Then system \eqref{sysnl} is
locally topologically linearizable at $(\bar{x},\bar{u})\in\Omega$ if, and
only if, it is locally quasi-\CCC{k} linearizable at $(\bar{x},\bar{u})$.
\end{thm}
\begin{proof}
  See section~\ref{sec-proofth}.
\end{proof}

Observe from \eqref{formuleConj}, that a quasi-\CCC{k}  
diffeomorphism $\chi$ is a linearizing homeomorphism if and only if it 
satisfies
\begin{equation}
  \label{eq:trsfpart}
  \frac{\partial\chi_{\mathrm{I}}}{\partial x}(x)\,f(x,u)\ \;=\;\ 
\;A\,\chi_{\mathrm{I}}(x)\ +\ B\,\chi_{\mathrm{I\!I}}(x,u)\ .
\end{equation}
Hence quasi-smooth linearizability is much easier to handle than 
topological linearizability, that relies on conjugating 
\emph{solutions} rather than equations.

\smallskip

System \eqref{eq:u3} of the introduction is topologically, quasi-\CCC\omega and
quasi-\CCC\infty linearizable at $(0,0)$ but fails to be even \CCC 1
linearizable; hence quasi-\CCC{k} cannot be
replaced with \CCC{k} in Theorem~\ref{lem-C1}. 
To study the gap between \CCC{k} and quasi-\CCC{k} linearizability,
note that \eqref{eq:trsfpart} imposes
additional regularity on a linearizing quasi-\CCC{k}  
diffeomorphism~:
\begin{prop}
\label{lissechi2}
Let $k\in\{\infty,\omega\}$ and $f$ in \eqref{sysnl} be \CCC{k}.
If $\chi:\Omega\to\Omega'$ is a quasi-\CCC{k} diffeomorphism that conjugates 
\eqref{sysnl} to the linear system $\dot{z}=Az+Bv$, then~:
\begin{enumerate}
\item \label{B-lisse} 
the map $B\chi_{\mathrm{I\!I}}:\Omega\to\RR^{m}$ is of class \CCC{k},
\item for any $(x,u)\in\Omega$ in the neighborhood of which the rank of
  $\partial f/\partial u$ is constant, one has
$\displaystyle
\mathrm{Rank}\;\frac{\partial f}{\partial u}(x,u)
=\mathrm{Rank}\,B
$.
  \label{eq:b=m}
\item for any open subset $O$ of $\Omega$, one has
$\displaystyle
\sup_{(x',u')\in O}\,\mathrm{Rank}\;\frac{\partial f}{\partial u}(x',u')
=\mathrm{Rank}\,B.
$
  \label{eq:b=mtilde}
\end{enumerate}
\end{prop}
\begin{proof} 
Point (\ref{B-lisse}) is direct consequence of \eqref{eq:trsfpart} and the
smoothness of $\chi_{\mathrm I}$ and $f$.
To establish (\ref{eq:b=m}) and (\ref{eq:b=mtilde}), differentiate
\eqref{eq:trsfpart} with respect to $u$ to obtain
\begin{equation}
\label{eq-rangB}
        \frac{\partial\chi_{\mathrm{I}}}{\partial x}(x)
        \frac{\partial{f}}{\partial{u}}(x,u)\ \;=\;\ 
        \frac{\partial (B\chi_{\mathrm{I\!I}})}{\partial u}(x,u)\ .
\end{equation}
Let $\mathcal{V}\subset\Omega$ be open and such that
$\mathrm{Rank}\,\partial f/\partial u(x,u)=\rho$ some integer $\rho$ and
all $(x,u)\in\mathcal{V}$. Define $\phi:\mathcal{V}\to\RR^{n+m}$ by 
$\phi(x,u)=(\chi_{\mathrm{I}}(x),B\chi_{\mathrm{I\!I}}(x,u))$.
On the one hand, since $\chi_{\mathrm{I}}$ is a diffeomorphism,
\eqref{eq-rangB} implies that the rank of the Jacobian of $\phi$ is $n+\rho$,
hence, by the constant rank theorem, $\phi(\mathcal{V})$ is a $(n+\rho)$
dimensional immersed sub-manifold of $\RR^{n+m}$; on the other hand, 
since $\chi$ is open, $\phi(\mathcal{V})$ is an open subset of the
$(n+\mathrm{Rank}\,B)$-dimensional linear range of $I_n\times B$; hence
$\rho=\mathrm{Rank}\,B$. This proves point \ref{eq:b=m}, and at the same time
point~\ref{eq:b=mtilde} because any $O\subset\Omega$ contains an open 
subset on
which the rank of $\partial f/\partial u$ is constant while
\eqref{eq-rangB} clearly implies that, for all $(x,u)\in\Omega$, the rank
of $\partial f/\partial u(x,u)$ is no larger than $\mathrm{Rank}\,B$.
\end{proof}

Based on Proposition \ref{lissechi2},  let us divide the
points of $\Omega$ into three classes.
\begin{itemize}
\item A point $(\bar{x},\bar{u})\in\Omega$ is called \emph{regular} if 
it has a neighborhood on which $\partial f/\partial u$ has constant rank. It
is easy to see that regular points form an open dense subset of $\Omega$.
\item If  $(\bar{x},\bar{u})$ is not regular, it is termed 
\emph{weakly singular} if each neighborhood $O\subset\Omega$ of this point
satisfies
\begin{equation}
  \label{eq:ws}
  \sup_{(x,u)\in O}
\mathrm{Rank}\;\frac{\partial f}{\partial u}(x,u)
=m.
\end{equation}
\item A point $(\bar{x},\bar{u})\in\Omega$ which is neither regular nor 
weakly singular is said to be \emph{strongly singular}.
This means it has a neighborhood $O\subset\Omega$, such that 
\begin{equation}
  \label{eq:ss}
  \sup_{(x,u)\in O}
\mathrm{Rank}\;\frac{\partial f}{\partial u}(x,u)
=m'<m\,,\ \ \ 
\mathrm{Rank}\;\frac{\partial f}{\partial u}(\bar x,\bar u)<m'\,.
\end{equation}
\end{itemize}
The distinction between topological and smooth linearizability 
may now be approached \textit{via} the following theorem that complements 
Theorem~\ref{lem-C1}.
\begin{thm}
\label{th-2}
Let $k\in\{\infty,\omega\}$ and 
$f$ be of class \CCC{k} on an open set $\Omega\subset\RR^{n+m}$.
\begin{itemize}
\item System \eqref{sysnl} is locally \CCC{k} linearizable at 
$(\bar{x},\bar{u})\in\Omega$ if, and
only if it is locally topologically linearizable 
at $(\bar{x},\bar{u})$ \emph{and} the latter is a regular point. 
\item If \eqref{sysnl} is
locally topologically linearizable 
at $(\bar{x},\bar{u})$ and the latter is a 
weakly singular point, then a linearizing homeomorphism around
$(\bar{x},\bar{u})$ may be chosen to be a map of class \CCC{k},
although not necessarily a \CCC{k}diffeomorphism (its inverse may
fail to be \CCC{k}).
\end{itemize}
\end{thm}
\begin{proof} The first assertion is a consequence of
Theorem \ref{lem-C1} together with Theorems~\ref{jak-res} and \ref{jak-res-quasi}
to come, observing that condition (\ref{2'}) in the latter will automatically 
hold at a regular point by the constant rank theorem. Next, assume that
$\chi:\Omega\to\Omega'$ is a quasi-\CCC{k} diffeomorphism that conjugates 
the \CCC{k} system \eqref{sysnl} 
to the linear controllable system $\dot{z}=Az+Bv$ at some weakly singular
point $(\bar{x},\bar{u})$. By (3) of Proposition 
\ref{lissechi2}, the rank of 
$B$ is $m$ hence it is left invertible; by (1) of the same 
proposition, $\chi_{\mathrm{I\!I}}$ is indeed \CCC{k}.
\end{proof}

Whether Theorem~\ref{th-2}
remains true if ``weakly singular'' gets replaced by ``strongly singular''
is unknown to the authors. This turns out to be equivalent 
to the following question in differential topology 
which is of interest in its own right and seems to have no answer so far. 

\begin{openq}
\label{open-q}
Let $O$ be a neighborhood of the origin in $\RR^{p+q}$ and
$F:O\rightarrow\RR^{p}$ a smooth (resp. real-analytic) map. 
Suppose $G:O\rightarrow\RR^q$ is a \emph{continuous} map 
such that $F\times G:O\rightarrow\RR^p\times\RR^q$ is a local homeomorphism 
at $0$.
\\[0.5ex]
Does there exist another neighborhood $O'\subset O$ of the origin 
and a \emph{smooth} (resp. \emph{real-analytic}) map $H:O'\rightarrow\RR^q$ 
such that $F\times H:O'\rightarrow\RR^p\times\RR^q$ is still
a local homeomorphism at $0$~?
\end{openq}

If the answer to the open question was yes, then 
Definitions~\ref{def-lindiff} and  \ref{def-quasi}
of quasi-smooth (resp. quasi-analytic) linearizability might 
equivalently require
$\chi$ to be smooth (resp. analytic) because, assuming the linear 
system is in normal form
\eqref{syslcan}-\eqref{syslcanJ}, one could set 
$F=\pi_{n+s_1}\circ\chi$ and smoothly (resp. analytically) redefine
the last $m-s_1$ components of $\chi$.

If the answer to the open question was no, then Definition \ref{def-lindiff}
would really be more general than the one obtained by restricting 
$\chi$ to be smooth (resp. analytic).
Indeed, if $F$ provides a counterexample to the open question, say, 
in the $\CCC \infty$ case, we may
consider on
$\RR^p\times O$ the control system
\begin{equation}
  \label{eq:contrex}
  \dot{x}\ \;=\;\ F(u)\ \ ,\ \ \ \ x\in\RR^{p},\;u\in\RR^{p+q}
\end{equation}
which is locally quasi-smoothly linearizable at the origin
because the local homeomorphism
\begin{displaymath}
\ (x,u)\ \mapsto\ (z,v)\ =\ (x,F(u),G(u))
\end{displaymath}
conjugates \eqref{eq:contrex} to
\begin{equation}
\label{linears}
  \dot{z}\ \;=\;\ B\,v ,~~~~{\rm with}~B=(\,I_p\,|0).
\end{equation}
However, no \emph{smooth} homeomorphism
\begin{displaymath}
  \chi\,:\ (x,u)\ \mapsto\ (z,v)\ =\ 
(\chi_{\mathrm{I}}(x),\chi_{\mathrm{I\!I}}(x,u))
\end{displaymath}
exists that quasi-smoothly linearizes  \eqref{eq:contrex} at $0$:
if this was the case, by Corollary \ref{cor-willems} we may assume up to
a linear change of variables that $\chi$ conjugates \eqref{eq:contrex}
to \eqref{linears}. Then conjugacy would imply
\begin{displaymath}
  \frac{\partial\chi_{\mathrm{I}}}{\partial x}(x)F(u)\ \;=\;\ 
\;B\,\chi_{\mathrm{I\!I}}(x,u)
\end{displaymath}
whence in particular
\begin{displaymath}
  F(u)\ \;=\;\ 
\left(\frac{\partial\chi_{\mathrm{I}}}{\partial x}(0)\right)^{-1}
\;B\,\chi_{\mathrm{I\!I}}(0,u),
\end{displaymath}
and the last $q$ components of $\chi_{\mathrm{I\!I}}(0,u)$ would yield a
smooth $H$ such that $F\times H$ is a local homeomorphism at $0$ in $R^{p+q}$,
contrary to the assumption.

\subsection{Geometric characterization of quasi smooth linearization}
\label{sec-JakRes}
Let $\mathcal{X}$ and $\mathcal{U}$ be
two open subsets of $\RR^{n}$ and $\RR^{m}$ respectively, and assume that $f$
is defined on $\mathcal{X}\times\mathcal{U}$.
For each 
$u\in\RR^m$, let $f_u$ be the vector field on $\mathcal{X}$ defined by~:
\begin{equation}
  \label{eq:deffu}
  f_{u}(x) \ \;=\;\  f(x,u).
\end{equation}
Also, for each $(x,u)\in\mathcal{X}\times\mathcal{U}$, we define below
a subspace $D(x,u)$, that coincides with the range of the linear mapping 
$\partial{f}/\partial{u}\,(x,u)$ when its dimension is locally constant.
First, we consider the subset
$\mathcal{L}_{x,u}\subset\RR^n$ (not a vector subspace) given by~:
\begin{eqnarray}
  \label{eq:defLxu}
y\in\mathcal{L}_{x,u}\ \;\Leftrightarrow&
\exists (w_n)\in\mathcal{U}^{\NN},&
\lim_{n\to\infty}w_n= u
\\\nonumber&&\hspace{-4.8ex}
{\rm and}~
\lim_{n\to\infty}\frac{f(x,w_n)-f(x,u)}{\|f(x,w_n)-f(x,u)\|}=y;
\end{eqnarray}
subsequently we put
\begin{eqnarray}
  \label{eq:defD}
  D(x,u) & = &
\mathrm{Span}_{\RR}\,\mathcal{L}_{x,u}\ .
\end{eqnarray}
In words, $D(x,u)$ is the vector space spanned by all limit directions 
of straight lines
through $f(x,u)$ and $f(x,u')$ as $u'$ approaches $u$ in $\RR^m$;
it is of common use in stratified geometry to generalize
the notion of tangent space.
Note that the set $\mathcal{L}_{x,u}$ depends on the norm used in
\eqref{eq:defLxu}, but the subspace $D(x,u)$ does not.
\begin{prop}
\label{prop-Dpartial}
If $f$ is of class $\CCC\infty$ and if we denote by ${\rm Ran}\,L$ the range
of a linear map $L$, we have that
\begin{eqnarray}
  \label{eq:egDpartial}
  D(x,u) &\supset& {\rm Ran}\;\frac{\partial{f}}{\partial{u}}(x,u)
\end{eqnarray}
and equality holds at every $(x,u)$ where 
the rank of $\displaystyle \partial{f}/\partial{u}\,(x,u)$ is
locally constant with respect to $u$.
\end{prop}
\begin{proof}
The inclusion \eqref{eq:egDpartial} holds because any nonzero
element of ${\rm Ran}\,{\partial{f}}/{\partial{u}}(x,u)$ can be written 
${\partial f}/{\partial u}(x,u).h$ for some $h$ in $\RR^m$, and one has
$$
\frac{\partial f/\partial u\,(x,u).h}{\|\partial{f}/\partial{u}\,(x,u).h\|}=
\lim_{t\to0^+}\frac{f(x,u+th)-f(x,u)}{\|f(x,u+th)-f(x,u)\|}\ .
$$
Now fix $(x,u)$ and assume that the rank of $\partial{f}/\partial{u}$ is
locally constant around $(x,u)$, equal to $r\leq m$ and use the constant rank-theorem.
Up to a permutation of coordinates,
$$(h_1,\ldots,h_m)\stackrel\lambda\mapsto(f_1(x,u+h)-f_1(x,u),\ldots,f_r(x,u+h)-f_r(x,u),h_{r+1},\ldots,h_m)$$
is a local diffeomorphism around zero in $\RR^m$ and, setting
$\rho=\lambda^{-1}\circ z\circ\lambda$ with
$z$ 
given by
$z(w_{1},\ldots,w_{m})=(w_{1},\ldots,w_{r},0,\ldots,0)$, there is a
constant $c$ such that 
\begin{equation}
\label{coercDf}
\|\rho(h)\|\leq c\,\|f(x,u+h)-f(x,u)\|
\ \ \mbox{and}\ \ f(x,u+\rho(h))=f(x,u+h)
\end{equation}
for all $h$.
Take $y\in\mathcal{L}_{x,u}$; by definition, there is a sequence
$(h_n)$ converging to zero and satisfying
\eqref{eq:defLxu} with $w_n=u+h_n$; from (\ref{coercDf}), we may re-write it as
\begin{equation}
\label{eq-ylim}
y=\lim_{n\to\infty}\frac{f(x,u+\rho(h_n))-f(x,u)}{\|\rho(h_n)\|}\,
\frac{\|\rho(h_n)\|}{\|f(x,u+h_n)-f(x,u)\|}
\end{equation}
where both ratios are bounded; extracting a sequence such that both
converge, the limit of the first ratio is, by definition of the derivative, 
$\frac{\partial{f}}{\partial{u}}(x,u).h$ with $h$ a limit point of
$\rho(h_n)/\|\rho(h_n)\|$; hence $y\in{\rm Ran}\,\partial{f}/\partial{u}\,(x,u)$.
We have proved that $\mathcal{L}_{x,u}\subset{\rm
  Ran}\,\partial{f}/\partial{u}\,(x,u)$. From \eqref{eq:defD}, this implies the reverse
inclusion of \eqref{eq:egDpartial} because the left-hand side is a
linear subspace.
\end{proof}

\smallskip

We can now characterize smooth (resp. analytic) and quasi-smooth
(resp. quasi-analytic) linearizability in parallel.
The proofs are given in section~\ref{sec-proof-JR}.
\begin{thm}[smooth or analytic linearizability]
\label{jak-res}
Let $k\in\{\infty,\omega\}$ and 
$f$ be of class \CCC{k} on an open set $\Omega\subset\RR^{n+m}$.
The control system \eqref{sysnl} is locally \CCC{k} linearizable at $(\bar{x},\bar{u})\in\Omega$ if, and only if there are open
neighborhoods $\mathcal{X}$ and $\mathcal{U}$ of $\bar{x}$ and
$\bar{u}$ in $\RR^{n}$ and  $\RR^{m}$, with 
$\mathcal{X}\times\mathcal{U}\subset\Omega$,
such that the following conditions are satisfied.
\begin{enumerate}
\item\  \label{jr-it1}
$D(x,u)$
does not depend on $u$ for $(x,u)\in\mathcal{X}\times\mathcal{U}$.
\item\  \label{jr-it3}
The rank of $\displaystyle\frac{\partial f}{\partial u}(x,u)$ is
constant in $\mathcal{X}\times\mathcal{U}$.
\item\  \label{jr-it2}
Defining on $\mathcal{X}$ the distribution $\Delta_{0}$ by
$\Delta_{0}(x)=D(x,u)$ ---~this is possible if point (\ref{jr-it1}) 
holds true~---
and inductively the flag of distributions $(\Delta_k)$ by~:
  \begin{eqnarray}\label{dist-JR}
    \Delta_{k+1} & = &
    \Delta_{k}\;+\;[\,f_{\bar{u}}\,,\,\Delta_{k}\,]
  \end{eqnarray}
where $[~,~]$ denotes the Lie bracket,
then each $\Delta_{k}$ for $0\leq k\leq n-1$ is integrable ({\it i.e.}
has constant dimension over $\RR$ and is closed
  under Lie bracket) and the rank of $\Delta_{n-1}$ is $n$.
\end{enumerate}
\end{thm}

\begin{thm}[quasi-smooth or quasi-analytic linearizability]
\label{jak-res-quasi}
Let $k\in\{\infty,\omega\}$ and 
$f$ be of class \CCC{k} on an open set $\Omega\subset\RR^{n+m}$.
The control system \eqref{sysnl} is locally quasi-\CCC{k}
linearizable at $(\bar{x},\bar{u})\in\Omega$ if, and only if there are open
neighborhoods $\mathcal{X}$ and $\mathcal{U}$ of $\bar{x}$ and
$\bar{u}$ in $\RR^{n}$ and  $\RR^{m}$, with 
$\mathcal{X}\times\mathcal{U}\subset\Omega$,
such that conditions (\ref{jr-it1}) and (\ref{jr-it2}) of Theorem~\ref{jak-res}
are met and, instead of
condition~(\ref{jr-it3}), it holds that
{\renewcommand{\theenumi}{$2'$}
\begin{enumerate}
\item\label{2'}
Denoting by $r_1\leq m$ the constant rank of $\Delta_0$, 
the mapping
\begin{equation}
\label{openmap}
\begin{array}{cccl}
        F\;: & \mathcal{X}\times\mathcal{U} & \rightarrow &
        \mathcal{X}\times\RR^{n}  \\
        & (x,u) & \mapsto & (\,x\,,\,f(x,u)\,)
\end{array}
\end{equation}
restricts to a $\CCC0$ fibration\footnote{\label{foot-fibre}
A $\CCC0$ fibration with fiber $\mathcal{F}$ over $\mathcal{B}$ is a 
continuous map
$g:\mathcal{E}\rightarrow \mathcal{B}$ for which every $\xi\in \mathcal{B}$ 
has a neighborhood
$\mathcal{O}$ in $\mathcal{B}$ such that 
$g^{-1}(\mathcal{O})\subset\mathcal{E}$ is homeomorphic
to $\mathcal{O}\times \mathcal{F}$, the so-called \emph{trivializing}
homeomorphism
$\psi:g^{-1}(\mathcal{O})\rightarrow\mathcal{O}\times \mathcal{F}$ 
being such that $\pi\circ\psi= g$ 
where $\pi:\mathcal{O}\times \mathcal{F}\rightarrow\mathcal{O}$ is the
natural projection onto the first factor.
}
$\mathcal{W}\rightarrow F(\mathcal{W})$
with fiber $\RR^{m-r_1}$
on some neighborhood $\mathcal{W}$ of $(\bar{x},\bar{u})$
in $\mathcal{X}\times\mathcal{U}$.
\end{enumerate}
}
\end{thm}
Theorem~\ref{jak-res} is of course equivalent to the results in
\cite{Jaku-Res80,Hunt-Su-Mey83,VdSc84scl}, but the conditions are stated
here in a slightly different form to parallel
Theorem~\ref{jak-res-quasi}.

\smallskip

\begin{cor}
\label{cor-ana}
  Assume that $f$ is real analytic on some open set $\Omega\subset\RR^{n+m}$. 
If the control system \eqref{sysnl} is locally \CCC\infty (resp. quasi-\CCC\infty)
linearizable at $(\bar{x},\bar{u})\in\Omega$, then it is also \CCC\omega
(resp. quasi-\CCC\omega) linearizable there.
\end{cor}
\begin{proof}
  analyticity does not appear in the conditions of the theorems, except for the
  regularity of $f$ itself.
\end{proof}

\subsection{Linearization versus equivalence to the linear approximation}
\label{lve}
For a control system, \emph{smooth} linearizability at an equilibrium implies
conjugacy \emph{to its linear approximation}:
\begin{prop}
\label{rmk-approxl-lisse}
Let $(\bar{x},\bar{u})$ be an equilibrium point of \eqref{sysnl},
i.e. $f(\bar{x},\bar{u})=0$, and let $\overline{A}=\partial f/\partial x
\,(\bar{x},\bar{u})$, $\overline{B}=\partial f/\partial u
\,(\bar{x},\bar{u})$ so that~:
\begin{equation}
  \label{eq:approx}
  f(x,u)\ =\ \overline{A}\,(x-\bar{x})
\;+\;\overline{B}\,(u-\bar{u})\;+\;\varepsilon(x-\bar{x},u-\bar{u})\;,
\end{equation}
where $\varepsilon$ is little $\rm{o}(\|x-\bar{x}\|+\|u-\bar{u}\|)$.

If system \eqref{sysnl} is locally smoothly linearizable at
$(\bar{x},\bar{u})$, then:
\\1. its linear approximation $(\overline{A},\overline{B})$ is controllable
(cf. Definition~\ref{def-eqlincont}),
\\2. the system is smoothly conjugate to $(\overline{A},\overline{B})$ at 
$(\bar{x},\bar{u})$.
\end{prop}
\begin{proof}
Let $\chi$ be a local diffeomorphism conjugating system \eqref{sysnl}
to $\dot{z}=Az+Bv$ at $(\bar{x},\bar{u})$, and observe
from \eqref{formuleConj} in Remark~\ref{rmk-conjdiff} that smooth 
linearizability translates into \eqref{eq:trsfpart}.
If we write $f$ as in \eqref{eq:approx}, and if we set 
$\overline{P}=\frac{\partial\chi_{\mathrm{I}}}{\partial x}(\bar{x})$,
$\overline{K}=\frac{\partial\chi_{\mathrm{I\!I}}}{\partial x}(\bar{x},\bar{u})$,
$\overline{Q}=\frac{\partial\chi_{\mathrm{I\!I}}}{\partial
  u}(\bar{x},\bar{u})$, we get by differentiating \eqref{eq:trsfpart} with
respect to $x$ and $u$ at $(\bar{x},\bar{u})$, using the 
relation $f(\bar{x},\bar{u})=0$, that
\begin{displaymath}
\overline{P}\,\overline{A}\ \;=\;\ A\,\overline{P}\;+\;B\,\overline{K}\ \ ,\ \ \ \ 
\overline{P}\,\overline{B}\ \;=\;\ B\,\overline{Q}\ .
\end{displaymath}
Since $\overline{P}$ and $\overline{Q}$ are square invertible matrices
by the triangular structure of $\chi$ displayed in (\ref{homeotri}), this
implies that the linear systems $(A,B)$ and
$(\overline{A},\overline{B})$ are linearly conjugate, 
see \eqref{lfeedABtilde}. Since  $(A,B)$ is controllable by definition
so is $(\overline{A},\overline{B})$, thereby achieving the proof.
\end{proof}

Proposition \ref{rmk-approxl-lisse}
has no analog if the control system is only 
topologically linearizable (hence \emph{quasi-smoothly} linearizable 
according to Theorem~\ref{lem-C1}).
For example, the system \eqref{eq:u3} in the introduction is quasi-\CCC\omega
linearizable at $(0,0)$, but its linear approximation $\dot{x}=0$ is not
controllable and it is not topologically equivalent to $\dot{x}=0$.
Apart from such degenerate cases, there also exist systems that are
quasi-analytically linearizable at some point with
controllable linear approximation there,
and still they are not conjugate to this linear approximation. An example
when $m=n=2$ is given by~:
\begin{displaymath}
\dot{x}_1\ =\ u_1\ \ ,\ \ \ \  
\dot{x}_2\ =\ x_1\;+\;u_2^{\,3}\ ,
\end{displaymath}
This system is quasi-analytically conjugate at $(0,0)$ to
\begin{equation}
\label{tata}
\dot{z}_1\ =\ v_1\ \ ,\ \ \ \  
\dot{z}_2\ =\ v_2\ ,
\end{equation}
via $z=x$, $v_1=u_1$, $v_2=u_2^{\,3}\;+\;x_1$. However, its linear
approximation at the origin is $\dot{x}_1=u_1$, $\dot{x}_2\ =\ x_1$, 
which is controllable yet not conjugate to \eqref{tata} 
({\it cf} Theorem~\ref{th-willems}).

\subsection{Non-genericity of linearizability}
\label{sec-generic}

Except when $m\geq n$ or $(n,m)=(2,1)$,
the conditions of Theorem \ref{jak-res} require a certain number of
\emph{equalities} (involving $f$ and its partial derivatives) to hold
\emph{everywhere}.
For example, the integrability of a distribution entails that all Lie brackets
be linearly dependent on the original vector fields, \textit{i.e.}
certain determinants must be identically zero.
This makes smooth (resp. analytic) linearizability of a smooth (resp.
analytic) control system highly 
non-generic in any reasonable sense, because when written in
proper jet spaces it is contained in
a set of infinite co-dimension.
Moreover, small perturbations of a system that does not satisfy these 
condition will not
satisfy them either, while most perturbations of a system which satisfies them 
will fail to do so. Compare for instance \cite{Tcho87} where it is shown
that the equivalence class of any system affine in the control
has infinite co-dimension in some Whitney topology.

From Theorem~\ref{th-2}, quasi-smooth or quasi-analytic
linearizability, 
hence also topological linearizability by Theorem~\ref{lem-C1},
require the same
equalities to hold on an open dense set, although this time some
singularities are allowed. This is no more ``generic'' than smooth 
linearizability,
as opposed to ODE's for which 
the Grobman-Hartman theorem allows one to linearize
around an equilibrium as soon as it is hyperbolic.

\section{Proofs}
\label{sec-proofs}

\subsection{Proof of Theorems~\ref{jak-res} and \ref{jak-res-quasi}}
\label{sec-proof-JR}
We begin with a lemma whose cumbersome index arrangement will be rewarded
later when constructing the Kronecker indices of the linearized system.
\begin{lem}
\label{lem-JR} Let $k\in\{\infty,\omega\}$.
  Let $\Delta_0$ and $f_{\bar{u}}$ be respectively a 
  distribution and a vector field, both of class \CCC{k} on a connected open
  neighborhood of $x\in\RR^n$. Let further $\Delta_i$, $i\geq0$, be the 
  distributions defined according to \eqref{dist-JR}, and set for convenience
$\Delta_{-1}=\{0\}$. Assume
  they satisfy point (\ref{jr-it2}) of Theorem~\ref{jak-res} or
  \ref{jak-res-quasi}.
Put
\begin{equation}
\label{JR-ri}
  r_i\ =\ \mathrm{Rank}\Delta_{i-1}\,,\ i>0,\ \ \ \ 
r_0=0\,,\ r_{-1}=-m\,,
\end{equation}
so that $r_i=n$ for some $i\leq n-1$; let $\rho\in\{3,\ldots,n+1\}$ be the
  smallest integer such that  $r_{\rho-1}=n$.
  Define also
\begin{equation}
\label{JR-nombres}
  s_i=r_i-r_{i-1}\,,\ \ 
\sigma_i=\sum_{j=i}^\rho s_j=n-r_{i-1}\,,\ \ 0\leq i\leq\rho\ 
\end{equation}
(note that $s_{\rho}=\sigma_{\rho}=0$).

 Then, there exists coordinates
$\chi_1,\ldots,\chi_n$ of class \CCC{k} on a neighborhood $\mathcal{X}$ of $x$  such that 
\begin{itemize}
\item $\chi_1,\ldots,\chi_{\sigma_i}$ are independent first
integrals of $\Delta_{i-2}$ for $i\in\{1,\ldots,\rho-1\}$,
\item $\chi_{\sigma_{i}+j}=f_{\bar{u}}\chi_{\sigma_{i+1}+j}$ for all integers
  $i,j$, $2\leq i\leq\rho-1$, $1\leq j\leq s_{i}$
\end{itemize}
($f_{\bar{u}}\chi_{\sigma_{i+1}+j}$ is the Lie derivative of  the
function $\chi_{\sigma_{i+1}+j}$ along the vector field $f_{\bar{u}}$).
\end{lem}
\begin{proof}
Note that when $i=1$, the first point above means that $\chi_1,\ldots,\chi_n$
are indeed local coordinates. 
Now, the Frobenius theorem provides us with  $n-r$ independent \CCC{k} 
first integrals for a \CCC{k} integrable
distribution of rank $r$. This accounts for the regularity of the coordinates
if we construct them as follows.

First pick $n-r_{\rho-2}=\sigma_{\rho-1}$ independent first integrals of
$\Delta_{\rho-3}$ and call them $\chi_1,\ldots,\chi_{\sigma_{\rho-1}}$; 
define further $\chi_{1+\sigma_{\rho-1}},\ldots,\chi_{2\sigma_{\rho-1}}$
by  $\chi_{\sigma_{\rho-1}+j}=f_{\bar{u}}\chi_{j}$ for
$1\leq j\leq \sigma_{\rho-1}=s_{\rho-1}$. 
Clearly, $\chi_1,\ldots,\chi_{\sigma_{\rho-1}+s_{\rho-1}}$ satisfy the 
conditions for $i=\rho-1$.
Then proceed inductively~: assume that, for some $i_0\in\{2,\ldots,\rho-1\}$,
the functions $\chi_1,\ldots,\chi_{\sigma_{i_0}+s_{i_0}}$ have been 
constructed and satisfy the conditions for $i\geq i_0$. 
\emph{We claim} that
the differentials $\D\chi_\ell$ are linearly independent at each point of 
$\mathcal{X}$. Indeed, 
assume that there is $\bar{x}\in\mathcal{X}$ and real 
coefficients $\mu_j$ and  $\lambda_k$  such that
\begin{equation}
\label{eqimp}
\sum_{j=1}^{\sigma_{i_0}} \mu_j\D\chi_j(\bar{x})+
\sum_{k=1+\sigma_{i_0+1}}^{\sigma_{i_0}}
\lambda_k\D\left(f_{\bar{u}}\chi_{k}\right)(\bar{x})=0.
\end{equation}
Put $\omega_1=\Sigma\mu_j\D\chi_j$ and 
$\omega_2=\Sigma\lambda_k\D\chi_k$.
Since $d$ commutes with the Lie derivative, we may rewrite (\ref{eqimp}) as
$\omega_1(\bar{x}) + f_{\bar{u}}\omega_2(\bar{x})=0$. In particular, for any 
$\CCC k$-vector field $X$ in $\Delta_{i_0-2}$, we get as $\omega_1(X)\equiv0$
that $f_{\bar{u}}\omega_2(X)(\bar{x})=0$. Now, by virtue of the formula
\begin{equation}
\label{transcro}
f_{\bar{u}}\left(\omega_2(X)\right)=f_{\bar{u}}\omega_2(X)+
\omega_2([f_{\bar{u}},X]),
\end{equation}
we obtain since $\omega_2(X)\equiv0$ that 
$\omega_2([f_{\bar{u}},X])(\bar{x})=0$, that is, $\omega_2$ annihilates
$\Delta_{i_0-1}$ at $\bar{x}$. But 
$\D\chi_1(\bar{x}),\ldots,\D\chi_{\sigma_{i_0+1}}(\bar{x})$ are a basis of the 
orthogonal space to $\Delta_{i_0-1}(\bar{x})$ by the induction hypothesis,
whereas $\omega_2(\bar{x})$ is a linear combination of the $\D\chi_k(\bar{x})$
for $\sigma_{i_0+1}< k\leq\sigma_{i_0}$. Therefore, since we know 
by the induction hypothesis that
the $\D\chi_\ell$ are point-wise independent for  $1\leq \ell\leq \sigma_{i_0}$,
we get that the $\lambda_k$ are zero and then the $\mu_j$ are also zero
by (\ref{eqimp}). \emph{This proves the claim}. Next, recall that 
$\chi_1,\ldots,\chi_{\sigma_{i_0} }$ are first integrals of $\Delta_{i_0-2}$,
thus \textit{a fortiori} of $\Delta_{i_0-3}$. For $X$ a $\CCC k$-vector field 
in the latter we deduce from (\ref{transcro}), where $\omega_2$ is replaced by 
$\D\chi_\ell$ with $1+\sigma_{i_0+1}\leq \ell\leq\sigma_{i_0}$, that
$\chi_{1+\sigma_{i_0}},\ldots,\chi_{\sigma_{i_0}+s_{i_0}}$ are also 
first integrals of $\Delta_{i_0-3}$. In case 
$\sigma_{i_0}+s_{i_0}<\sigma_{i_0-1}$, pick
$\chi_{\sigma_{i_0}+s_{i_0}+1},\ldots,\chi_{\sigma_{i_0-1}}$ so that
$\chi_\ell$ for $1\leq \ell\leq \sigma_{i_0-1}$ is a complete set of 
independent integrals of $\Delta_{i_0-3}$. If $i_0=2$ we are done, otherwise
define $\chi_{\sigma_{i_0-1}+j}=f_{\bar{u}}\chi_{\sigma_{i_0}+j}$ 
for $1\leq j\leq s_{i_0-1}$ in order to complete the induction step.
\end{proof}

\begin{proof}[Proof of Theorems~\ref{jak-res} and \ref{jak-res-quasi}]
The two proofs run parallel to each other.

We first show necessity, assuming that $k=\infty$ for analyticity does 
not appear in the conclusions. 
Assume local (quasi) smooth linearizability, \textit{cf.}
Definitions~\ref{def-lindiff} and \ref{def-quasi}. 
Without loss of generality, we assume that
$\Omega=\mathcal{X}\times\mathcal{U}$ where $\mathcal{X}$ and $\mathcal{U}$
are open neighborhoods of $\bar{x}$ and $\bar{u}$ in $\RR^n$ and $\RR^m$ 
respectively. Let $\chi:\mathcal{X}\times\mathcal{U}\to\Omega'\subset\RR^{n+m}$
be as in \eqref{homeotri}; recall that $\chi_{\mathrm{I}}$ is a smooth diffeomorphism
$\mathcal{X}\to\chi_{\mathrm{I}}(\mathcal{X})$. 
We may also assume, after composing $\chi$ with a linear invertible map, that
the pair $(A,B)$ is in canonical form
\eqref{syslcan}-\eqref{syslcanJ}, but we still write $A,B$ rather than $A_c,B_c$.
Denote by $B_0,\ldots,B_m$ the columns of $B$ and define the vector fields
$b_0,\ldots,b_m$ on $\RR^n$ by
\begin{equation}
  \label{bobi}
  b_i(z)=B_i\,,\ 1\leq i\leq m\,,\ \ \ \ \ b_0(z)=Az+B\bar{u}\ 
\end{equation}
and the distributions $\Lambda_i$ by
\begin{equation}
  \label{Lambdai}
\Lambda_0(z)=\mathrm{Span}_{\RR}\{b_1(z),\ldots,b_m(z)\}=\mathrm{Ran}B
\,\ \ \ \ 
\Lambda_{i+1}=\Lambda_i+[b_0,\Lambda_i]\,,\ 1\leq i\leq m\,.
\end{equation}
From \eqref{eq:trsfpart}, we have
\begin{equation}
  \label{eq:trsfpart2}
  \frac{\partial\chi_{\mathrm{I}}}{\partial x}(x)\,f(x,u)\ \;=\;\ 
b_0(\chi_{\mathrm{I}}(x))\ +\ B
\left(\chi_{\mathrm{I\!I}}(x,u)-\chi_{\mathrm{I\!I}}(x,\bar{u})\right).
\end{equation}
Since $\chi$ is a triangular homeomorphism,
$\chi_{\mathrm{I\!I}}(x,w)-\chi_{\mathrm{I\!I}}(x,u)$ covers an open 
neighborhood of $0$ in $\RR^m$ when $w$ ranges around $u$
in $\RR^m$. Thus, in view of (\ref{eq:trsfpart2}), 
$\mathcal{L}_{x,u}$ defined by
\eqref{eq:defLxu} contains an open set in 
$\left(\frac{\partial\chi_{\mathrm{I}}}{\partial
    x}(x)\right)^{-1}\mathrm{Ran}B$, and by double inclusion
$$
D(x,u)\ =\ \left(\frac{\partial\chi_{\mathrm{I}}}{\partial
    x}(x)\right)^{-1} \mathrm{Ran}B.
$$
This proves point \ref{jr-it1}, and also proves that the distribution
$\Delta_0$ in point \ref{jr-it2} is the pullback of $\Lambda_0$ by the
diffeomorphism $\chi_{\mathrm{I}}$, \textit{i.e.}
$\left(\chi_{\mathrm{I}}\right)_*\Delta_0=\Lambda_0$.
Since \eqref{eq:trsfpart2} also implies
$\left(\chi_{\mathrm{I}}\right)_*f_{\bar{u}}=b_0$, we have
$\left(\chi_{\mathrm{I}}\right)_*\Delta_i=\Lambda_i$ for all $i$. This gives
point \ref{jr-it2} because it is obviously true with $\Lambda_i$ instead of
$\Delta_i$, and integrability and ranks are preserved by
conjugation with the smooth diffeomorphism $\chi_{\mathrm{I}}$.
In the case of smooth linearizability, point \ref{jr-it3} is easily obtained
by differentiating \eqref{eq:trsfpart2} with respect to $u$ and 
using invertibility of 
$\partial\chi_{\mathrm{I\!I}}/\partial u(x,u)$.

To conclude the proof of necessity, let us prove point \ref{2'} in the case of
quasi-smooth linearizability. Let
$$\mathcal{M}=\{\ \ (x,y)\in\mathcal{X}\times\RR^{n};
~~\frac{\partial\chi_{\mathrm{I}}}{\partial x}(x)\,y
\;-\;
A\chi_{\mathrm{I}}(x)\ \;\in\;\ {\rm Ran}\,B\ \ \}\ .$$
This is a smooth embedded sub-manifold of $\mathcal{X}\times\RR^{n}$ of 
dimension $n+r_1$, where $r_1=\mathrm{Rank}\,B\leq m$.
If we define $F$ as in \eqref{openmap}, 
it is clear from \eqref{eq:trsfpart} that
\begin{displaymath}
F(\mathcal{X}\times\mathcal{U})\ \;\subset\;\ \mathcal{M}\ .
\end{displaymath}
Now, take some $(m-r_1)\times m$ matrix $C$ whose rows 
complement $r_1$ independent rows of $B$ into a basis of $\RR^m$.
Pick matrices $E_{1}$ and $E_{2}$ of appropriate sizes such that
\begin{displaymath}
  E_{1}\,B\;+\; E_{2}\,C\ \;=\;\ I_m\ .
\end{displaymath}
By \eqref{eq:trsfpart} we get
\begin{equation}
\label{Bezoutchi}
        E_{1}
\left[\;
\vphantom{  \frac{\partial\chi_{\mathrm{I}}}{\partial x}  }
\;\frac{\partial\chi_{\mathrm{I}}}{\partial x}(x)\,f(x,u)\;-\;A\,\chi_{\mathrm{I}}(x)
\;\right]
\;+\;E_{2}C\chi_{\mathrm{I\!I}}(x,u)\ \;=\;\ 
\chi_{\mathrm{I\!I}}(x,u)\ .
\end{equation}
Define
$$\psi:\mathcal{X}\times\mathcal{U}
\rightarrow\mathcal{M}\times\RR^{m-r_1}$$
by the formula:
$$\psi(x,u)=(x,f(x,u),C\chi_{\mathrm{I\!I}}(x,u)).$$
From \eqref{Bezoutchi}, this mapping has an inverse given by
\begin{eqnarray*}
  \psi^{-1}:\ \ \psi(\mathcal{X}\times\mathcal{U})&\to&
\mathcal{X}\times\mathcal{U}\\
(x,y,z)&\mapsto&
\chi^{-1}\left(
\ \chi_{\mathrm{I}}(x)\ ,\ E_{1}\left[
\frac{\partial\chi_{\mathrm{I}}}{\partial x}(x)\,y\,-\,A\,\chi_{\mathrm{I}}(x)
\right]\;+\;E_{2}\,z\ \right)
\end{eqnarray*}
so that $\psi$ defines a homeomorphism from $\mathcal{X}\times\mathcal{U}$
onto its image which is open in $\mathcal{M}\times\RR^{m-r_1}$ by
invariance of the domain.
Let $\mathcal{O}$ be a neighborhood of
$(\bar{x},f(\bar{x},\bar{u}))$ in $\mathcal{M}$ and $\mathcal{S}$
an open ball centered at $C\chi_{\mathrm{I\!I}}(\bar{x},\bar{u})$ in
$\RR^{m-r_1}$ such that
$\mathcal{O}\times\mathcal{S}\subset\psi(\mathcal{X}\times\mathcal{U})$,
and take
$\mathcal{W}=\psi^{-1}(\mathcal{O}\times\mathcal{S})$.
Then $F:\mathcal{W}\rightarrow F(\mathcal{W})=\mathcal{O}$ is a $\CCC0$ 
fibration with fiber $\mathcal{S}$ and trivializing homeomorphism 
$\psi:\mathcal{W}\to\mathcal{O}\times\mathcal{S}$.
Since $\mathcal{S}$ is homeomorphic to $\RR^{m-r_1}$, 
condition \ref{2'} follows.

\bigskip

We turn to sufficiency. Points \ref{jr-it1},
\ref{jr-it2}, and either \ref{jr-it3} or \ref{2'} imply, 
for all $x\in\mathcal{X}$,
\begin{equation}
  \label{pr1}
  \Delta_{0}(x)\ =\ 
\mathrm{Span}_{\RR}\{f(x,w)-f(x,u),
(u,w)\in\mathcal{U}\times\mathcal{U}\}\ .
\end{equation}
Indeed the right-hand side always contains $D(x,u)$ because it contains all
the differences $f(x,w_n)-f(x,u)$ in \eqref{eq:defLxu}, and
point~\ref{jr-it1} implies the reverse inclusion because
$f(x,w)-f(x,u)$ can be computed as the integral on the segment $[u,w]\subset\mathcal{U}$
of a function that, thanks to Proposition~\ref{prop-Dpartial}, belongs
constantly to $\Delta_{0}(x)$.

From \eqref{pr1}, the distribution $\Delta_0$ is of class \CCC{k}. 
Considering point~\ref{jr-it2}, we may apply Lemma~\ref{lem-JR}.
We thus obtain some, with $r_i$, $s_i$ and $\sigma_i$ the integers defined by
\eqref{JR-ri} and \eqref{JR-nombres}, some \CCC{k} 
coordinates $\chi_1,\ldots,\chi_n$ on a 
neighborhood of $\bar{x}$ possibly smaller than $\mathcal{X}$ (but that we
continue to denote by $\mathcal{X}$), i.e. a diffeomorphism
$\chi_{\mathrm{I}}:\mathcal{X}\to\chi_{\mathrm{I}}(\mathcal{X})$, with
$\chi_{\mathrm{I}}=(\chi_1,\ldots,\chi_n)$, meeting the conclusions of
Lemma~\ref{lem-JR}.
In particular, $\chi_1,\ldots,\chi_{n-r_1}$ are
first integrals of the distribution $\Delta_0$, and from \eqref{pr1}, this
implies that $\partial \chi_i/\partial x(x)\,f(x,u)$ does
not depend on $u$, and is there fore equal to its value for $u=\bar{u}$~:
\begin{equation}
  \label{pr2}
\frac{\partial\chi_i}{\partial x}(x)f(x,u)
=
f_{\bar{u}}\chi_i\,(x)
\,,\ \ 1\leq i\leq n-r_1\ .
\end{equation}
For larger $i$, the left-hand side depends on $x$ and $u$~: define
$\lambda:\mathcal{X}\times\mathcal{U}\to\RR^{m_1}$ by
\begin{equation}
  \label{pr3}
  \lambda(x,u)\ =\ 
(\;\frac{\partial\chi_{n-r_1+1}}{\partial x}(x)f(x,u)
\;,\,\ldots\,,\;\frac{\partial\chi_{n}}{\partial x}(x)f(x,u)\;)\ .
\end{equation}
Then, defining coordinates $z_1,\ldots,z_n$ by $z=\chi_{\mathrm{I}}(x)$. 
The equations of system (\ref{sysnl})
are as follows (the first line gives the derivatives of the $n-r_1$ first
coordinates and the second line the last $r_1$ ones)~:
\begin{equation}
  \label{pr12}
  \begin{array}{rcll}
\dot{z}_{\sigma_{i+1}+j} &=& z_{\sigma_{i}+j}\,,\ \ \ 
&2\leq i\leq\rho-1\,,\ 1\leq j\leq s_i\,,
\\
\dot{z}_{\ell} &=& \lambda_{n-\ell}(\,\chi_{\mathrm{I}}^{-1}(z)\,,\,u\,)\,,\ \
\ 
&n-m_1+1\leq\ell\leq n\ .
  \end{array}
\end{equation}

If point \ref{jr-it3} is satisfied, the rank of the map
$(x,u)\mapsto(\chi_{\mathrm{I}}(x),\frac{\partial\chi_{\mathrm{I}}}{\partial(x)}f(x,u))$
is constant and thus, according to (\ref{pr1}), it is equal to $n+r_1$, $r_1$
being the rank of $\Delta_0$. From \eqref{pr2}, the map
$(x,u)\mapsto(\chi_{\mathrm{I}}(x),\lambda(x,u))$ has the same constant rank $n+r_1$. Hence there exists
$\phi~:\mathcal{X}\times\mathcal{U}\to\RR^{m-r_1}$ such that
\begin{equation}
  \label{didif}
  (x,u)\mapsto(\,\chi_{\mathrm{I}}(x)\,,\,\lambda(x,u)\,,\,
\phi(x,u)\,)
\end{equation}
is a diffeomorphism of class $\CCC{k}$. 
Obviously, defining $\chi_{\mathrm{I\!I}}$ by 
$
\chi_{\mathrm{I\!I}}(x,u)
=
(\lambda(x,u),
\phi(x,u))$
yields a $\CCC{k}$ diffeomorphism $\chi$ that conjugates
\eqref{sysnl} to a linear controllable system $\dot{z}=Az+Bu$. This proves
sufficiency in Theorem~\ref{jak-res}.

If point \ref{2'} is satisfied instead, let
$\psi:\mathcal{W}\to F(\mathcal{W})\times\RR^{n-r_1}$ be the ``trivializing''
homeomorphism. Recall that, with 
$\pi:F(\mathcal{W})\times\RR^{n-r_1}\to F(\mathcal{W})$ the natural
projection, one has $\pi\circ\psi=F$; call $\phi:\mathcal{W}\to\RR^{n-r_1}$
the map such that $\psi=F\times\phi$
Composing $F$ with
$(x,\xi)\mapsto(\chi_{\mathrm{I}}(x),\frac{\partial\chi_{\mathrm{I}}}{\partial(x)}\xi)$,
one gets that
$(x,u)\mapsto\chi(x,u)=(\chi_{\mathrm{I}}(x),\lambda(x,u),\phi(x,u))$ is a
homeomorphism. It clearly conjugates 
\eqref{sysnl} to a linear controllable system $\dot{z}=Az+Bu$. This proves
sufficiency in Theorem~\ref{jak-res-quasi}.
\end{proof}

\subsection{Proof of Theorem \ref{lem-C1}}
\label{sec-proofth}
This theorem for $k=\omega$ is consequence of this theorem for $k=\infty$ and
of Corollary~\ref{cor-ana}. Hence we only have to prove it for $k=\infty$,
\textit{i.e.} we assume that $f$ is infinitely differentiable and we prove that
topological linearizability implies quasi-\CCC\infty linearizability.

Without loss of generality, we suppose that
$(\bar{x},\bar{u})=(0,0)$. 
Assume there exists a homeomorphism 
$\chi$ from a
neighborhood of the origin in $\RR^{n+m}$ to an open subset of $\RR^{n+m}$
that conjugates system \eqref{sysnl} to the linear controllable system
\begin{equation}
\label{sysl}
        \dot{z}\ \;=\;\ Az\;+\;Bv
\end{equation}
with $z\in\RR^{n}$ and $v\in\RR^{m}$.
Composing $\chi$ with a linear invertible map allows us to suppose that
the pair $(A,B)$ is in canonical form
\eqref{syslcan}-\eqref{syslcanJ}, i.e. that \eqref{sysl} can be read
\begin{equation}
\label{lincan}
  \dot{z}_{\sigma_{i}+k}\ =\ z_{\sigma_{i-1}+k}\ ,\ \ \ 
2\leq i\leq\rho,\ 1\leq k\leq s_{i-1},
\end{equation}
where the integers $s_i$ and $\sigma_i$ were defined in
\eqref{indices} and  \eqref{sigmai} and where, for notational compactness,
we have set~:
\begin{equation}
\label{vconv}
        z_{n+k}\ \stackrel\Delta=\ v_{k}\ ;
\end{equation}
recall here that $s_0=m$, and notice that $s_1<m$ may well occur
as it simply means that 
$\mathrm{Rank}\,B<m$, in which case some
of the controls do not appear in the canonical form.
With the aggregate notation~:
\begin{equation}
\label{aggnot1}
Z_{j}\ \stackrel\Delta= \ \left(
\begin{array}{c}
        z_{\sigma_{j+1}+1}\\\vdots\\z_{\sigma_{j}}
        \end{array}
\right)
,\ 1\leq j\leq \rho-1\ ,\ \ \ \ 
Z_{0}\ \stackrel\Delta=\ 
\left(
\begin{array}{c}
        v_{1}\\\vdots\\v_{m}
        \end{array}
\right)
\ ,
\end{equation}
and the matrices $J_{r}^{s}$ defined in \eqref{syslcanJ}, system 
\eqref{lincan} can be rewritten as
\begin{equation}
\label{agg2}
       \begin{array}{rcl}
\dot{Z}_{\rho-1} & = & J^{s_{\rho-1}}_{s_{\rho-2}}Z_{\rho-2} \\
\dot{Z}_{\rho-2} & = & J^{s_{\rho-2}}_{s_{\rho-3}}Z_{\rho-3} \\
                 & \vdots &   \\
\dot{Z}_{2} & = & J^{s_{2}}_{s_{1}}Z_{1}\\
\dot{Z}_{1} & = & J^{s_{1}}_{s_{0}}Z_{0}
        \end{array}
\end{equation}
and is viewed as a control system with state $(Z_{\rho-1},\ldots,Z_1)$ and
control $Z_0$.
We also make the convention, similar to \eqref{vconv}, that
\begin{equation}
\label{uconv}
        x_{n+k}\ \stackrel\Delta=\ u_{k}\ ,
\end{equation}
and we use for the controls the aggregate notation~:
\begin{equation}
\label{uconvagg}
        X_{0}\ \stackrel\Delta=\ 
\left(
\begin{array}{c}
        x_{n+1}\\\vdots\\x_{n+m}
        \end{array}
\right)
        \ =\ 
\left(
\begin{array}{c}
        u_{1}\\\vdots\\u_{m}
        \end{array}
\right)\ .
\end{equation}

Let us now prove that property
$\mathcal{P}_{\ell}$ below is true for $0\leq \ell\leq \rho-1$.

{\it
\textbf{Property $\mathbf{\mathcal{P}_{\ell}}$ :\ }
there exists a smooth local change of coordinates around $0$ in $\RR^n$, say
\begin{displaymath}
(x_{1},\ldots,x_{n})\ \;\mapsto\;\ 
(\widehat{X},X_{\ell},\dots,X_{2},X_{1}),
\end{displaymath}
with $\widehat{X}\in\RR^{\sigma_{\ell+1}}$ and
$X_{i}\in\RR^{s_{i}}$ for $0\leq i\leq\ell$
(if $\ell=0$ there are no $X_{i}$'s beyond $X_0$ whereas
if $\ell=\rho-1$ there is no $\widehat{X}$),
after which system \eqref{sysnl} reads:
\begin{equation}
\label{agg1rec}
        \begin{array}{rcl}
\dot{\widehat{X}} & = & \widehat{F}(\widehat{X},X_{\ell})  \\
\dot{X}_{\ell} & = & F_{\ell}(\widehat{X},X_{\ell},X_{\ell-1})  \\
                 & \vdots &   \\
\dot{X}_{2} & = & F_{2}(\widehat{X},X_{\ell},\ldots,X_{1})\\
\dot{X}_{1} & = & F_{1}(\widehat{X},X_{\ell},\ldots,X_{1},X_{0})
        \end{array}\,,
\end{equation}
and such that \eqref{agg1rec}, viewed as a control system with state
$(\widehat{X},X_{\ell},\ldots,X_{1})$ and control $X_0$, is locally 
topologically conjugate at $(0,0)$ to system \eqref{agg2}
{\it via} a local homeomorphism
$$(\widehat{X},X_{\ell},\dots,X_{1},X_{0})\mapsto(Z_{\rho-1},\ldots,Z_{0})$$
which is, together with its inverse, of the block triangular form~:
\begin{displaymath}
        \begin{array}{lcrcl}
(Z_{\rho-1},\ldots,Z_{\ell+1}) \ =\ 
\widehat{\Phi}(\widehat{X})                   &&\widehat{X} & = &
\widehat{\Psi}(Z_{\rho-1},\ldots,Z_{\ell+1})  \\
Z_{\ell} \ =\  \Phi_{\ell}(\widehat{X},X_{\ell})       &&X_{\ell} &
= & \Psi_{\ell}(Z_{\rho-1},\ldots,Z_{\ell})  \\
\ \ \ \ \ \;\vdots &&& \vdots &   \\
Z_{1} \ =\  \Phi_{1}(\widehat{X},X_{\ell},\ldots,X_{1})
&&X_{1} & = & \Psi_{1}(Z_{\rho-1},\ldots,Z_{1})  \\
Z_{0} \ =\  \Phi_{0}(\widehat{X},X_{\ell},\ldots,X_{1},X_{0})&&X_{0} & = &
\Psi_{0}(Z_{\rho-1},\ldots,Z_{1},Z_{0})
        \end{array}
\end{displaymath}
where $\Phi_{i}$ and $\Psi_{i}$ are, for $1\leq i\leq \ell$, continuously
differentiable with respect to $X_{i}$ and $Z_{i}$ respectively,
have an invertible derivative, 
and satisfy for $1\leq i\leq \ell$ the relation~:
\begin{equation}
\label{inv1}
        \begin{array}{l}
                F_{i}(\widehat{X},X_{\ell},\ldots,X_{i},X_{i-1})
                 \ \;=\;\  F_{i}(\widehat{X},X_{\ell},\ldots,X_{i},0)
\\
 \ \ \ \ \ 
                 +\;\ 
\left(
          \frac{\partial\Phi_{i}}{\partial X_{i}}
          (\widehat{X},X_{\ell},\ldots,X_{i})\right)^{-1}
          J^{s_{i}}_{s_{i-1}}\,\left(
      \vphantom{\left(\frac{\partial\Phi_{1}}{\partial X_{1}}
         (X_{\rho-1},\ldots,X_{1})\right)^{-1}}
\right.
\\
 \ \ \ \ \ 
\ \ \ \ \ \ \ \ \ \ \ \ 
                 \left.
      \vphantom{\left(\frac{\partial\Phi_{1}}{\partial X_{1}}
         (X_{\rho-1},\ldots,X_{1})\right)^{-1}}
                 \Phi_{i-1}(\widehat{X},X_{\ell},\ldots,X_{i},X_{i-1})
                 -\Phi_{i-1}(\widehat{X},X_{\ell},\ldots,X_{i},0)\right)\ ;
        \end{array}
\end{equation}
furthermore, the partial homeomorphism
\begin{equation}
\label{trsfred}
        (\widehat{X},X_{\ell})\ \;\mapsto\;\ (Z_{\rho-1},\ldots,Z_{\ell})
\end{equation}
locally topologically conjugates, at $(0,0)\in\RR^{\sigma_{\ell+1}+s_\ell}$,
the reduced control system
\begin{eqnarray}
\label{aggrecred1}
\dot{\widehat{X}} & = & \widehat{F}(\widehat{X},X_{\ell}),
\end{eqnarray}
with state $\widehat{X}$ and input $X_{\ell}$, to the reduced linear control
system 
\begin{equation}
\label{aggrecred2}
        \begin{array}{rcl}
\dot{Z}_{\rho-1} & = & J^{s_{\rho-1}}_{s_{\rho-2}}Z_{\rho-2}\ , \\
                 & \vdots &   \\
\dot{Z}_{\ell+1} & = & J^{s_{\ell+1}}_{s_{\ell}}Z_{\ell}\ 
        \end{array}
\end{equation}
wit state
$(Z_{\rho-1},\ldots,Z_{\ell+1})$ and input $Z_{\ell}$.
}

Indeed, $\mathcal{P}_{0}$ is merely the
original assumption on local topological conjugacy 
of systems \eqref{sysnl} and \eqref{agg2}, where the triangular structure
\eqref{homeotri} of the conjugating homeomorphism was taken into account; 
note that, in $\mathcal{P}_{0}$, \eqref{inv1} is empty and that the reduced
system \eqref{aggrecred1} is the original system.
Next, supposing
that $\mathcal{P}_{\ell}$ holds for some $\ell\geq 0$, we apply
Lemmas \ref{lem-recur1} and \ref{lem-recur2} (see below)
to the reduced systems \eqref{aggrecred1}, \eqref{aggrecred2}, and to the 
partial homeomorphism \eqref{trsfred}, with 
$$d=\sigma_{\ell+1},~ r=s_\ell,~s=s_{\ell+1},~ U=X_{\ell},~ 
(x_{1},\ldots,x_{d})=\widehat{X},$$
$$Z^{1}=(Z_{\rho-1},\ldots,Z_{\ell+2}),~ Z^{2}=Z_{\ell+1},~{\rm  and}~
V=Z_{\ell},$$
and then, upon renaming $\widetilde{X}^2$ as $X_{\ell+1}$, 
$\widetilde{f}^{2}$ as 
$F_{\ell+1}$, and choosing $\widetilde{X}^{1}$ to be the new $\widehat{X}$, 
we get $\mathcal{P}_{\ell+1}$.

Now, $\mathcal{P}_{\rho-1}$, where we specialize \eqref{inv1} to $i=1$, 
provides us with 
a \emph{smooth} change of variables around $0$ in $\RR^n$:
\begin{displaymath}
        (x_{1},\ldots,x_{n})\ \;\mapsto\;\   (X_{\rho-1},\dots,X_{2},X_{1})
\end{displaymath}
with $X_{i}\in\RR^{s_{i}}$ such that, in the new coordinates,
system \eqref{sysnl} reads
\begin{equation}
\label{agg1}
        \begin{array}{rcl}
\dot{X}_{\rho-1} & = & F_{\rho-1}(X_{\rho-1},X_{\rho-2})  \\
\dot{X}_{\rho-2} & = & F_{\rho-2}(X_{\rho-1},X_{\rho-2},X_{\rho-3})  \\
                 & \vdots &   \\
\dot{X}_{2} & = & F_{2}(X_{\rho-1},\ldots,X_{1})\\
\dot{X}_{1} & = & F_{1}(X_{\rho-1},\ldots,X_{1},X_{0}),
        \end{array}
\end{equation}
and also such that the local homeomorphism $\Phi$ that topologically conjugates
system \eqref{agg1} to system \eqref{agg2} at $(0,0)$ is, together with its
inverse $\Psi$, of the triangular form~:
\begin{equation}
\label{trsf-agg}
        \begin{array}{rclcrcl}
Z_{\rho-1} & = & \Phi_{\rho-1}(X_{\rho-1})
   &&X_{\rho-1} & = & \Psi_{\rho-1}(Z_{\rho-1})
\\
Z_{\rho-2} & = & \Phi_{\rho-2}(X_{\rho-1},X_{\rho-2})
   &&X_{\rho-2} & = & \Psi_{\rho-2}(Z_{\rho-1},Z_{\rho-2})
\\
  & \vdots &
   && & \vdots &
\\
Z_{1} & = & \Phi_{1}(X_{\rho-1},\ldots,X_{1})
   &&X_{1} & = & \Psi_{1}(Z_{\rho-1},\ldots,Z_{1})
\\
Z_{0} & = & \Phi_{0}(X_{\rho-1},\ldots,X_{1},X_{0})
   &&X_{0} & = & \Psi_{0}(Z_{\rho-1},\ldots,Z_{1},Z_{0}),
        \end{array}
\end{equation}
where the following three properties hold~:
\begin{enumerate}
        \item\ 
Each $\Phi_{k}$ and $\Psi_{k}$ for $k\geq 1$ is continuously
differentiable with respect to $X_{k}$ and $Z_{k}$ respectively;
in particular, $\partial\Phi_k/\partial X_k$ is invertible
throughout the considered neighborhood.
        \item\ 
For $k\geq 2$, $\displaystyle\mathrm{Rank}
\frac{\partial F_{k}}{\partial X_{k-1}}(0,\ldots,0)\ =\ s_{k}\ $, 
i.e. this rank is maximum, equal to the number of rows.

        \item\ 
$F_{1}$ satisfies
\vspace{-1ex}
\begin{eqnarray}
\label{propro}
\!
                F_{1}(X_{\rho-1},\ldots,X_{1},X_{0})
                 & = & F_{1}(X_{\rho-1},\ldots,X_{1},0)  \\
\nonumber
                 &  & \!+\;\ 
                 \left(\frac{\partial\Phi_{1}}{\partial X_{1}}
                     (X_{\rho-1},\ldots,X_{1})\right)^{-1}
                 J^{s_{1}}_{m}\,\left(
      \vphantom{\left(\frac{\partial\Phi_{1}}{\partial X_{1}}
         (X_{\rho-1},\ldots,X_{1})\right)^{-1}}
                 \right.
                 \\
\nonumber
&  & \! \! \! \! \! \!
                 \left.
      \vphantom{\left(\frac{\partial\Phi_{1}}{\partial X_{1}}
(X_{\rho-1},\ldots,X_{1})\right)^{-1}}\Phi_{0}(X_{\rho-1},\ldots,X_{1},X_{0})
                 -\Phi_{0}(X_{\rho-1},\ldots,X_{1},0)\right)\ .
\end{eqnarray}
\end{enumerate}
From the maximum
rank assumption on $\partial F_{\rho-1}/\partial X_{\rho-2}$, 
it is possible to define $Y_{\rho-2}$ whose first
$s_{\rho-1}$ entries are those of $F_{\rho-1}(X_{\rho-1},X_{\rho-2})$ and
whose remaining $s_{\rho-2}-s_{\rho-1}$ entries are suitable components
of $X_{\rho-2}$, in such a way that
\begin{displaymath}
        (X_{\rho-1},\ldots,X_{1})\ \;\mapsto\;\ 
        (X_{\rho-1},Y_{\rho-2},X_{\rho-3}\ldots,X_{1})
\end{displaymath}
is a local \emph{smooth} change of coordinates around $0$ in $\RR^n$.
After performing this change of coordinates and setting 
$Y_{\rho-1}=X_{\rho-1}$ for 
notational homogeneity, system \eqref{agg1} reads
\begin{displaymath}
        \begin{array}{rcl}
\dot{Y}_{\rho-1} & = & J^{s_{\rho-1}}_{s_{\rho-2}}\,Y_{\rho-2} \\
\dot{Y}_{\rho-2} & = &
\widetilde{F}_{\rho-2}(Y_{\rho-1},Y_{\rho-2},X_{\rho-3}) \\
                 & \vdots &   \\
\dot{X}_{2} & = & \widetilde{F}_{2}(Y_{\rho-1},Y_{\rho-2},X_{\rho-3},\ldots,X_{1})\\
\dot{X}_{1} & = & \widetilde{F}_{1}(Y_{\rho-1},Y_{\rho-2},X_{\rho-3},\ldots,X_{1},X_{0})
        \end{array}
\end{displaymath}
where the $\widetilde{F}$'s enjoy the same properties than the $F$'s,
in particular the maximality of 
$\mathrm{Rank}\,\partial \widetilde{F}_{k}/\partial X_{k-1}(0,\ldots,0)$
for $\rho-2\geq k\geq2$.
One may iterate this procedure, limited only by the fact that the 
maximum rank property mentioned above only holds for $k\geq 2$
but not necessarily for $k=1$.
Altogether, this yields a \emph{smooth} local change of coordinates 
around $0$ in $\RR^n$~:
\begin{displaymath}
        (X_{\rho-1},\ldots,X_{1})\ \;\mapsto\;\ (Y_{\rho-1},\ldots,Y_{1}),
\end{displaymath}
after which system \eqref{agg1} is of the form
\begin{equation}
\label{linY}
        \begin{array}{rcl}
\dot{Y}_{\rho-1} & = & J^{s_{\rho-1}}_{s_{\rho-2}}\,Y_{\rho-2} \\
                 & \vdots &   \\
\dot{Y}_{2} & = & J^{s_{2}}_{s_{1}}\,Y_{1}\\
\dot{Y}_{1} & = & F_{1}(Y_{\rho-1},,\ldots,Y_{1},X_{0})\ ,
\end{array}
\end{equation}
where we abuse the notation $F_{1}$ for simplicity because,
although it needs not be the same as in \eqref{agg1},
this new $F_{1}$ enjoys the same property
\eqref{propro} for some suitably redefined $\Phi_{1}$ and $\Phi_{0}$.
Now, we may rewrite \eqref{propro} as
\begin{equation}
\label{contlin}
F_{1}(Y_{\rho-1},,\ldots,Y_{1},X_{0})\ \;=\;\ 
J^{s_{1}}_{m}\,H(Y_{\rho-1},,\ldots,Y_{1},X_{0})
\end{equation} 
where $H$, in the aggregate notation
$Y=(Y_{\rho-1},,\ldots,Y_{1})$, is defined by
\begin{displaymath}
\!\!H(Y,X_{0})\ =\ 
\left(\!\!\!\begin{array}{c}F_{1}(Y,0)\\0\end{array}\!\!\right)
\,+\,
\left(\!\begin{array}{cc}
\frac{\partial\Phi_{1}}{\partial Y_{1}}(Y)^{-1} & \!\!\!\!\!0
\\
0                                               & \!\!\!\!\!I_{m-s_{1}}
\end{array}\!\right)
     \big(\Phi_{0}(Y,X_{0})-\Phi_{0}(Y,0)\big).
\end{displaymath}
Since $\Phi$ has the triangular structure displayed in
\eqref{trsf-agg}, the map
$X_{0}\mapsto \Phi_{0}(Y,X_{0})$ is injective 
for fixed $Y=(Y_{\rho-1},\ldots,Y_{1})$ in the neighborhood of
$0$ where it is defined in $\RR^m$. Consequently,
$(Y,X_{0})\mapsto (Y,H(Y,X_{0}))$ is also injective  in the neighborhood of
$0$ where it is defined in $\RR^{n+m}$; since it is continuous, it is a local
homeomorphism of $\RR^{n+m}$ at $(0,0)$ by invariance of the domain, 
and then \eqref{linY}, \eqref{contlin} make it clear that 
system \eqref{agg1} is locally quasi-smoothly linearizable at this
point.

Since \eqref{agg1} is smoothly conjugate to the original system \eqref{sysnl},
this proves local quasi-smooth linearizability of the latter hence the theorem.

\bigskip

\noindent\textbf{Two lemmas.}\hspace{1ex}
The following two lemmas are applied recursively in the above proof of
Theorem~\ref{lem-C1} to obtain the forms \eqref{agg1}, \eqref{agg2},
and \eqref{trsf-agg}.
Although these lemmas  team up into a single result in the above-mentioned
proof, they have been stated here separately for the sake of clarity.

We will consider two control systems with state in $\RR^{d}$ and control in 
$\RR^{r}$. Expanded in coordinates, the first system reads
        \begin{equation}
          \label{eq:1}
          \begin{array}{rcl}
        \dot{x}_1 &=& f_{1}(x_{1},\ldots,x_{d},x_{d+1},\ldots,x_{d+r})
          \\
        &\vdots&
          \\
        \dot{x}_{d} &=&
        f_{d}(x_{1},\ldots,x_{d},x_{d+1},\ldots,x_{d+r})\ ,
          \end{array}
        \end{equation}
with state variable $(x_1,\ldots,x_d)$ and control variable 
$(x_{d+1},\ldots,x_{d+r})\in\RR^{r}$, the functions $f_1,\cdots,f_d$ being
smooth $\RR^{d+r}\to\RR$. The second system has 
state variable $(z_1,\ldots,z_d)$ and control variable 
$(z_{d+1},\ldots,z_{d+r})\in\RR^{r}$, and it assumes the special form~:
        \begin{equation}
        \label{lin}
          \begin{array}{rcl}
        \dot{z}_1 & = & g_{1}(z_{1},\ldots,z_{d})
        \\  & \vdots & \\
        \dot{z}_{d-s} &=&  g_{d-s}(z_{1},\ldots,z_{d})
          \\
        \dot{z}_{d-s+1} &=& z_{d+1}
         \\  & \vdots & \\
        \dot{z}_{d} &=& z_{d+s}\ ,
          \end{array}
        \end{equation}
where $0<s\leq d$ and $s\leq r$ while $g_1,\cdots,g_{d-s}$ are again smooth
$\RR^d\to\RR$. Nothing prevents us here from having $s<r$,
in which case some of the controls do not enter the equation.
It will be convenient to use the aggregate notations
\[
\begin{array}{ll}
  X\stackrel\Delta=(x_{1},\ldots,x_{d})\,,&~~
  U\stackrel\Delta=(x_{d+1},\ldots,x_{d+r})\,,\\
  Z\stackrel\Delta=(z_{1},\ldots,z_{d})\,,&~~
  V\stackrel\Delta=(z_{d+1},\ldots,z_{d+r})\,,
\end{array}
\]
and to further split $Z$ into $(Z^1,Z^2)$ with
\begin{equation}
\label{decZ}
  Z^{1}\stackrel\Delta=(z_{1},\ldots,z_{d-s})\,,\ \ 
  Z^{2}\stackrel\Delta=(z_{d-s+1},\ldots,z_{d})\,,
\end{equation}
so as to write \eqref{eq:1} in the form
\begin{equation}
  \label{eq:nl-agg}
  \dot{X}\ \;=\;\ f(X,U)
\end{equation}
and \eqref{lin} as
\begin{equation}
\label{syslinagg}
  \begin{array}{rcl}
\dot{Z}^{1}&=& g^{1}(Z^{1},Z^{2})\\
\dot{Z}^{2}&=& J_{r}^{s}\,V\,,
  \end{array}
\end{equation}
with $J^{s}_{r}$ the $s\times r$ matrix, defined in \eqref{syslcanJ}, 
that selects the first $s$ entries of a vector.

\begin{lem}
\label{lem-recur1}
Let $d$, $r$ and $s$ be strictly positive integers with $s\leq d$ and
$s\leq r$.
        Suppose, for some $\varepsilon>0$, that 
\[\varphi:(-\varepsilon,\varepsilon)^{d+r}\to
\RR^{d+r}
\]
is a homeomorphism onto its image, with inverse $\psi$,
that conjugates system \eqref{eq:nl-agg} to system \eqref{syslinagg}.
Then, there exists $0<\varepsilon'<\varepsilon$
and a \emph{smooth} local change of coordinates around
$0\in\RR^{d}$~:
\[\theta:\ (-\varepsilon',\varepsilon')^{d}\ \;\to\;\ 
\theta\bigl((-\varepsilon',\varepsilon')^{d}\bigr)\subset
(-\varepsilon,\varepsilon)^{d}
\]
that fixes the origin and is such that, in the new coordinates
$\widetilde{X}=\theta^{-1}(X)$, 
both the system \eqref{eq:nl-agg} and the conjugating homeomorphism
$\widetilde{\varphi}=\varphi\circ(\theta\times\mathrm{id})$ 
assume a block triangular structure with respect to the partition
$\widetilde{X}=(\widetilde{X}^{1},\widetilde{X}^{2})$, where
$\widetilde{X}^{1}\stackrel\Delta=(\widetilde{x}_1,\ldots,\widetilde{x}_{d-s})$ and $\widetilde{X}^{2}\stackrel\Delta=(\widetilde{x}_{d-s+1},\ldots,\widetilde{x}_{d})$; that is to say,
on $(-\varepsilon',\varepsilon')^{d+r}$, we have that
  \begin{itemize}
\item\  system \eqref{eq:1} reads~:
        \begin{equation}
        \label{eq:1bis}
        \begin{array}{rcl}
        \dot{\widetilde{X}}^{1} &=&
        \widetilde{f}^{1}(\widetilde{X}^{1},\widetilde{X}^{2}) \\
        \dot{\widetilde{X}}^{2} &=&
        \widetilde{f}^{2}(\widetilde{X}^{1},\widetilde{X}^{2},U),
        \end{array}
        \end{equation}
\item\ On their respective domains of definition,  the homeomorphism 
$\widetilde{\varphi}$ and its
inverse $\widetilde{\psi}=(\theta^{-1}\times {\rm id})\circ\psi$
read ~:
        \begin{equation}
        \label{phipsitilde}
\begin{array}{rclcrcl}
Z^{1} & = &
\widetilde{\varphi}^{1}(\widetilde{X}^{1})
&&
\widetilde{X}^{1} & = & \widetilde{\psi}_{1}(Z^{1})
\\
Z^{2} & = &
\widetilde{\varphi}^{2}(\widetilde{X}^{1},\widetilde{X}^{2})
&&
\widetilde{X}^{2} & = & \widetilde{\psi}_{2}(Z^{1},Z^{2})
\\
V & = &
\widetilde{\varphi}^{3}(\widetilde{X}^{1},\widetilde{X}^{2},U)
&&
U & = & \widetilde{\psi}_{3}(Z^{1},Z^{2},V)\,.
\end{array}
\end{equation}
\end{itemize}
\end{lem}

\begin{lem}
\label{lem-recur2}
Let 
\[\widetilde{\varphi}:(-\varepsilon',\varepsilon')^{d+r}\to
\RR^{d+r}
\]
be a homeomorphism onto its image, having the block
triangular structure displayed in
\eqref{phipsitilde}, and assume that it conjugates the smooth
system \eqref{eq:1bis} to the smooth system \eqref{syslinagg}. 
Necessarily then, $\widetilde{\varphi}$
has the following properties~:
\begin{enumerate}
  \item\  \label{it-phi2diff}
The map $\widetilde{\varphi}^{2}$ is continuously differentiable with
respect to its second argument $\widetilde{X}^{2}$, and
$\displaystyle\frac{\partial\widetilde{\varphi}^{2}}{\partial\widetilde{X}^{2}}(0,0)$ is
invertible.
  \item\  \label{it-f2}
On some neighborhood of $0\in\RR^{d+r}$ included in 
$(-\varepsilon',\varepsilon')^{d+r}$, one has~:
\begin{eqnarray}
\label{relf2}
                \widetilde{f}^{2}(\widetilde{X}^{1},\widetilde{X}^{2},U)
                 & = &
\\
\nonumber
\widetilde{f}^{2}(\widetilde{X}^{1},\widetilde{X}^{2},0)
                 &\!+&\!\!
\left(\frac{\partial\widetilde{\varphi}^{2}}{\partial\widetilde{X}^{2}}
                     (\widetilde{X}^{1},\widetilde{X}^{2})\right)^{-1}\!\!\!\!
J^{s}_{r}\,\left(
\widetilde{\varphi}^{3}(\widetilde{X}^{1},\widetilde{X}^{2},U)
-\widetilde{\varphi}^{3}(\widetilde{X}^{1},\widetilde{X}^{2},0)
\right)
\end{eqnarray}
  \item\  \label{it-conjred}
On  some neighborhood of $0\in\RR^d$ included in 
$(-\varepsilon',\varepsilon')^d$, the partial homeomorphism 
\begin{equation}
  \label{homredlem}
  (\,\widetilde{X}^{1}\,,\,\widetilde{X}^{2}\,)\ \;\mapsto\;\ 
(\,\widetilde{\varphi}^{1}(\widetilde{X}^{1})\,,\,\widetilde{\varphi}^{2}(\widetilde{X}^{1},\widetilde{X}^{2})\,)
\end{equation}
conjugates the control system
\begin{equation}
  \label{sysredX}
\dot{\widetilde{X}^{1}} =
        \widetilde{f}^{1}(\widetilde{X}^{1},\widetilde{X}^{2}),
\end{equation}
with state $\widetilde{X}^{1}$ and control $\widetilde{X}^{2}$, to the
control system 
\begin{equation}
  \label{sysredZ}
\dot{Z}^{1} = g^{1}(Z^{1},Z^{2})
\end{equation}
with state $Z^{1}$ and input $Z^{2}$.
\end{enumerate}
\end{lem}
Note that \eqref{sysredX} and \eqref{sysredZ}
are reduced systems from \eqref{eq:1bis} and \eqref{syslinagg}.

\begin{proof}[Proof of Lemma \ref{lem-recur1}]
Since the homeomorphism $\varphi$ conjugates \eqref{eq:nl-agg} to
\eqref{syslinagg}, we know, by Proposition~\ref{prop-0}, that $\varphi$ and
$\psi$ split component-wise into~:
        \begin{equation}
\label{PHIPSI}
\begin{array}{rclcrcl}
Z & = & \varphi_{\mathrm{I}}(X)
&\ \ \ &
X & = & \psi_{\mathrm{I}}(Z)
\\
V & = & \varphi_{\mathrm{I\!I}}(X,U)
&&
U & = & \psi_{\mathrm{I\!I}}(Z,V)\ \ .
\end{array}
\end{equation}

Consider the map $f:\,(-\varepsilon,\varepsilon)^{d+r}\to\RR^d$ given in 
\eqref{eq:nl-agg}, and let us define
$g:\varphi((-\varepsilon,\varepsilon)^{d+r})\rightarrow\RR^d$ analogously
from  \eqref{syslinagg}, namely
$g$ is the concatenated map whose first $d-s$ components are
given by $g^{1}(Z)$ and whose last $s$ components are
given by $J^{s}_{r}V$.
Define two families of continuous vector fields  $\mathcal{F}'$ 
and $\mathcal{G}'$, on 
$(-\varepsilon,\varepsilon)^{d}$ and 
$\varphi_{\mathrm{I}}((-\varepsilon,\varepsilon)^{d})$ respectively,
by the following formulas (compare \eqref{F'gen})~:
\begin{eqnarray}
        \label{F'}
\mathcal{F}' & = & \{\,\delta f_{\alpha_1,\alpha_2}\,;
\alpha_{1},\alpha_{2}\mbox{ feedbacks on $(-\varepsilon,\varepsilon)^{d+r}$}
\,\}\ ,
\\
        \label{G'} 
\mathcal{G}' & = & \{\,\delta g_{\beta_1,\beta_2}\,;
\beta_{1},\beta_{2}\mbox{ feedbacks on
  $\varphi\bigl((-\varepsilon,\varepsilon)^{d+r}\bigr)$} \,\}\ .
\end{eqnarray}
Applying Proposition \ref{prop-01} twice, first to
$\chi=\varphi$ and then to $\chi=\psi$, we see
that each integral curve of a vector field in $\mathcal{F}'$
is mapped by $\varphi_{\mathrm{I}}$ to some
integral curve of a vector field in $\mathcal{G}'$
and {\it vice-versa} upon replacing  $\varphi_{\mathrm{I}}$ by  
$\psi_{\mathrm{I}}$. This shows in particular that uniqueness of solutions to
the Cauchy problem associated to vector fields is preserved, i.e. if we define 
the families of vector fields (compare \eqref{F''gen})~:
\begin{eqnarray}
  \label{F''}
  \mathcal{F}'' & = & 
\{\,Y\in\mathcal{F}'\,,\ Y\mbox{ has a flow}\,\}\ ,
\\
  \label{G''}
  \mathcal{G}'' & = & 
\{\,Y\in\mathcal{G}'\,,\ Y\mbox{ has a flow}\,\}\ ,
\end{eqnarray}
we also have that each integral curve of a vector field in $\mathcal{F}''$
is mapped by $\varphi_{\mathrm{I}}$ to an integral curve of a vector field in
$\mathcal{G}''$ and {\it vice-versa} upon replacing  $\varphi_{\mathrm{I}}$ by  
$\psi_{\mathrm{I}}$. 
By concatenation, using Proposition~\ref{prop-FF'}, it follows that
\begin{equation}
  \label{conjorb}
\left.
  \begin{array}{l}
\textit{for any } X\in\,(-\varepsilon,\varepsilon)^{d},
\;\varphi_{\mathrm{I}}
\textit{ defines a homeomorphism, }\\
\textit{for the orbit topologies, from the orbit of } 
\mathcal{F}''
\textit{ through } $X$\\
\textit{onto the orbit of } \mathcal{G}''
\textit{ through } \varphi_{\mathrm{I}}(X),
  \end{array}
\right\}
\end{equation}
where the orbit topology as described in Proposition~\ref{prop-FF'}
(by definition the restriction of $\varphi_{\mathrm{I}}$ is bi-continuous 
for the topologies induced by the ambient space; 
bi-continuity for the orbit topologies requires the description of these
topologies as given in Proposition~\ref{prop-FF'}). 

Now, the vector fields $\delta g_{\beta_{1},\beta_{2}}$ appearing
in \eqref{G'} inherit from the structure of $g$, displayed 
in \eqref{syslinagg}, the following particular form~:
\begin{equation}
\label{deltag}
\delta g_{\beta_{1},\beta_{2}}(Z)
\ \;=\;\ 
\left(
\begin{array}{c}
0\\\vdots\\0\\
\beta_{1,1}(Z)\,-\,\beta_{2,1}(Z)
\\\vdots\\
\beta_{1,s}(Z)\,-\,\beta_{2,s}(Z)
\end{array}
\right)\ ,
\end{equation}
where $\beta_{i,1},\ldots,\beta_{i,s}$ designate, for $i=1,2$, the first $s$
component of the feedback $\beta_{i}$. 
This will allow for us to describe explicitly the 
orbits of $\mathcal{G}''$, namely~:
\begin{equation}
  \label{eq:orbG}
\left.
  \begin{array}{l}
\textit{the orbit of }\mathcal{G}''\textit{ through }
Z_0=(c_1,\ldots,c_d)
\\
\textit{is the connected component containing $Z_0$ of the set}
\\
  \{Z\in \varphi_{\mathrm{I}}\left((-\varepsilon,\varepsilon)^d\right)\,,\ 
z_1=c_1,\ldots,z_{d-s}=c_{d-s}\}.
  \end{array}
\right\}
\end{equation}
Indeed, the orbit in question is contained in this set, because it is 
connected, and because all the vector fields in  $\mathcal{G}''$
have their first $d-s$ components equal to zero by \eqref{deltag}.

To prove the reverse inclusion, it is enough to show that
the orbit of $\mathcal{G}''$ through $Z_0$, denoted hereafter by
$\mathcal{O}_{\mathcal{G}'',Z_0}$,
contains all the points sufficiently close to $Z_0$ having the same first 
$d-s$ coordinates as $Z_0$. Indeed, since $Z_0$ was arbitrary,
this will imply that the connected
component defined by \eqref{eq:orbG} splits into a disjoint union of
open orbits hence consists of a single one by connectedness.
That is to say, putting $Z_0=(Z_0^1,Z_0^2)$ according to \eqref{decZ},
\ref{eq:orbG} will follow from the existence of a $\rho>0$
such that
\begin{equation}
\label{etaporb}
\{Z_0^1\}\times B(Z_0^2,\rho)=B(Z_0,\rho)\cap\mathcal{O}_{\mathcal{G}'',Z_0}.
\end{equation}

Now, it follows from Remark \ref{rmk-topoprod} that, for sufficiently small
$\rho$, each connected component of
$B(Z_0,\rho)\cap\mathcal{O}_{\mathcal{G}'',Z_0}$
is an embedded sub-manifold of $B(Z_0,\rho)$.
Then, the connected component of 
$B(Z_0,\rho)\cap\mathcal{O}_{\mathcal{G}'',Z_0}$
containing $Z_0$ is, by inclusion, an embedded sub-manifold
of the linear manifold $\{Z_0^1\}\times B(Z_0^2,\rho)$.
In particular, since no strict sub-manifold can be densely embedded in a 
given manifold, we see that \eqref{etaporb} will hold is only we can prove
that
\begin{equation}
\label{etaporbf}
  \begin{array}{l}
\!\!\!\!\!\!\!\!\textit{The connected component containing $Z_0$ of } 
B(Z_0,\rho)\cap\mathcal{O}_{\mathcal{G}'',Z_0}\\
\!\!\!\!\!\!\!\!\textit{is dense in } \{Z_0^1\}\times B(Z_0^2,\rho)\ 
\textit{for the Euclidean topology.}
\\
  \end{array}
\end{equation}
To prove \eqref{etaporbf}, pick $V_0$ such that
$(Z_0,V_0)\in
\varphi\bigl((-\varepsilon,\varepsilon)^{d+r}\bigr)$
and observe, since the latter is an open set, that 
shrinking $\rho$ further, if necessary, allows us to assume
$\overline{B}(Z_0,\rho)\times \overline{B}(V_0,\rho)
\subset\varphi\bigl((-\varepsilon,\varepsilon)^{d+r}\bigr)$.
\emph{We claim that any continuous map 
$\overline{B}(Z_0,\rho)\to \overline{B}(V_0,\rho)$ extends
to a feedback on $\varphi\bigl(
(-\varepsilon,\varepsilon)^{d+r}\bigr)$}. 
Indeed, in view of the one-to-one correspondence 
$\beta\to \psi\carre\beta$
between feedbacks on $\varphi\bigl(
(-\varepsilon,\varepsilon)^{d+r}\bigr)$ 
and feedbacks on $(-\varepsilon,\varepsilon)^{d+r}$ 
({\it cf} the discussion leading to \eqref{defcarre}-\eqref{act-inv}),
it is enough to prove that every continuous map 
$\psi_{\mathrm{I}}\bigl(\overline{B}(Z_0,\rho)\bigr)\to
(-\varepsilon,\varepsilon)^r$ extends to a continuous
map 
$(-\varepsilon,\varepsilon)^d\to(-\varepsilon,\varepsilon)^r$, and this in 
turn follows from the Tietze extension theorem since 
$\psi_{\mathrm{I}}\bigl(\overline{B}(Z_0,\rho)\bigr)$ is closed
in $(-\varepsilon,\varepsilon)^d$ and since
$(-\varepsilon,\varepsilon)^r$ is a poly-interval.
\emph{This proves the claim}.

From the claim, it follows that the restriction to 
$\overline{B}(Z_0,\rho)$ of the $\RR^s$-valued vector field 
$J_r^s(\beta_1(Z)-\beta_2(Z))$,
accounting for the lower half of the right-hand side in \eqref{deltag},
can be assigned \emph{arbitrarily}, by choosing adequately
the feedbacks $\beta_1$ and $\beta_2$, among 
\emph{continuous} vector fields
$\overline{B}(Z_0,\rho)\to \overline{B}(0,\rho)$ 
(take $\beta_2$ to
extend the constant map $V_0$ on $\overline{B}(Z_0,\rho)$).
Of course, the corresponding vector field 
$\delta g_{\beta_{1},\beta_{2}}$ in \eqref{deltag} belongs to
$\mathcal{G}'$ but not necessarily to $\mathcal{G}''$ since
continuous vector fields need not have a flow. However, since 
$\delta g_{\beta_{1},\beta_{2}}$ has a flow at least when
$\beta_1$ and $\beta_2$ are smooth, we deduce from Proposition
\ref{prop-feed-lisse} that the restriction to 
$\overline{B}(Z_0,\rho)$ of the vector fields in  $\mathcal{G}''$
are of the form $\{0\}\times Y$, where $Y$ ranges over a uniformly dense
subset $\Upsilon$ of all $\RR^s$-valued continuous maps
$\overline{B}(Z_0,\rho)\to \overline{B}(0,\rho)$.
Now, every point in $B(Z_0^2,\rho)$ can be attained
from $Z_0^2$ upon integrating, \emph{within} $B(Z_0^2,\rho)$,
a constant vector field of arbitrary small norm.
By Lemma \ref{lem-ode3} applied with $\mathcal{U}=B(Z_0^2,\rho)$ 
and $K=\{Z_0^2\}$, the corresponding trajectory can be approximated
uniformly by integral curves  \emph{that remain in} $B(Z_0^2,\rho)$
of vector fields in $\Upsilon$. Therefore,
every point 
in $\{z_0^1\}\times B(Z_0^2,\rho)$ is the limit of endpoints
of integral curves of $\mathcal{G}''$ \emph{that remain in} 
$\{z_0^1\}\times B(Z_0^2,\rho)$, which proves \eqref{etaporbf}
and thus \eqref{eq:orbG}.
In particular, the orbits of $\mathcal{G}''$ are
\emph{embedded} sub-manifolds in 
$\varphi_{\mathrm{I}}\bigl((-\varepsilon,\varepsilon)^d\bigr)$.

Next, we turn to the orbits of $\mathcal{F}''$, and we designate
by $\mathcal{O}_{\mathcal{F}'',p}$ the orbit of $\mathcal{F}''$ in
$]-\varepsilon,\varepsilon[^{d}$ through the point $p$. 
On the one hand, Proposition \ref{prop-FF'} and Theorem~\ref{th-orb} show that
$\mathcal{O}_{\mathcal{F}'',p}$ is a smooth immersed sub-manifold of
$]-\varepsilon,\varepsilon[^{d}$.  
On the other hand, by \eqref{conjorb}, this immersed sub-manifold is sent 
homeomorphically  by $\varphi_{\mathrm{I}}$, 
both for the orbit topology and the ambient topology, onto 
$\mathcal{O}_{\mathcal{G}'',\varphi_{\mathrm{I}}(p)}$
which is a smooth \emph{embedded} $s$-dimensional sub-manifold of
$\varphi_{\mathrm{I}}\bigl((-\varepsilon,\varepsilon)^{d}\bigr)$, as we saw from \eqref{eq:orbG}.
This entails that 
all orbits of $\mathcal{F}''$ in
$]-\varepsilon,\varepsilon[^{d}$ are \emph{embedded} sub-manifolds of 
dimension $s$.
Consequently, still from Proposition \ref{prop-FF'} and Theorem~\ref{th-orb},
there are coordinates $(\xi_1,\ldots,\xi_d)$ defined
on an open neighborhood $W_{0}$ of the origin in
$]-\varepsilon,\varepsilon[^{d}$ ---this neighborhood may be assumed 
to be of the form 
$\{(\xi_1,\ldots,\xi_d),\,|\xi_i|<\varepsilon'\}$ ---
such that, in these coordinates, 
\begin{displaymath}
W_{0}\cap \mathcal{O}_{\mathcal{F}'',0}\ \;=\;\ 
\{\,(\xi_1,\ldots,\xi_d),\,\mbox{with}
\;(\xi_{s+1},\ldots,\xi_{d})\in T\,\}\ ,
\end{displaymath}
with $T$ a subset of $]-\varepsilon',\varepsilon'[^{d-s}$ containing 
$(0,\ldots,0)$,
the tangent space to $W_{0}\cap \mathcal{O}_{\mathcal{F}'',0}$ at 
each of its points being \emph{spanned} by
$\partial/\partial\xi_{1},\ldots,\partial/\partial\xi_{s}$,
while at any point $p\in W_{0}$ the vector fields 
$\partial/\partial\xi_{1},\ldots,\partial/\partial\xi_{s}$
\emph{belong} to the tangent space of $\mathcal{O}_{\mathcal{F}'',p}$.
But since we saw that \emph{all} orbits are smooth
sub-manifolds of dimension $s$, these vector fields actually \emph{span}
the tangent space to the orbit at every point. 
Hence all the vector fields $\delta f_{\alpha_1,\alpha_2}$ 
in $\mathcal{F}''$
have their last $d-s$ components equal to zero on $W_{0}$
in the $\xi$ coordinates, and this holds in particular when $\alpha_1$,
$\alpha_2$ range over all constant feedbacks 
$(-\varepsilon,\varepsilon)^d\to(-\varepsilon,\varepsilon)^r$.
This implies, by the very definition 
of $\delta f_{\alpha_1,\alpha_2}$, 
that $(\dot{\xi}_{s+1},\ldots,\dot{\xi}_{d})$ ---
as computed from  \eqref{eq:nl-agg} upon performing the change of variable
$X\mapsto(\xi_1,\ldots,\xi_d)$ ---
does not depend on the control variable $U$.
Choose for $\widetilde{X}$ the $\xi$ 
coordinates arranged in
reverse order, and let $\widetilde{f}$ be the analog
of $f$ in the new coordinates $(\widetilde{X},U)$. 
Then the first $d-s$ components of $\widetilde{f}$ 
do not depend on $U$ so that \eqref{eq:1bis} holds.
Moreover, if $\widetilde{\varphi}$ denotes the new homeomorphism that
conjugates \eqref{eq:1bis} to \eqref{syslinagg} over 
$(-\varepsilon,\varepsilon)^{d+r}$, 
$\widetilde{\varphi}((-\varepsilon,\varepsilon)^{d+r})$, and if
$\widetilde{\psi}$ denotes its inverse, it follows from 
\eqref{conjorb} and the above characterization of the orbits that
$\widetilde{\varphi}_{\mathrm{I}}$  maps the sets where
$\tilde{x}_1,\ldots,\tilde{x}_{d-s}$
are constant to those where $z_1,\ldots,z_{d-s}$ are constant, thus
the functions
$\widetilde{\varphi}_{1},\ldots,\widetilde{\varphi}_{d-s}$ and
$\widetilde{\psi}_{1},\ldots,\widetilde{\psi}_{d-s}$ 
depend only on their $d-s$ first arguments whence \eqref{phipsitilde}
follows.
\end{proof}

\begin{proof}[Proof of Lemma \ref{lem-recur2}]
We use again the concatenated notation
$\widetilde{\varphi}_{\mathrm{I}}=(\widetilde{\varphi}^1,\widetilde{\varphi}^2)$,
$\widetilde{\psi}_{\mathrm{I}}=(\widetilde{\psi}^1,\widetilde{\psi}^2)$,
these partial homeomorphisms being inverse of each other.
Let $(Z_0,V_0)\in\widetilde{\varphi}((-\varepsilon',\varepsilon')^{d+r})$
and $\varepsilon''$ be so small that
the product neighborhood $(Z_0,V_0)+(-\varepsilon'',\varepsilon'')^{d+r}$ lies
entirely within $\widetilde{\varphi}((-\varepsilon',\varepsilon')^{d+r})$.
The restriction to $(Z_0,V_0)+(-\varepsilon'',\varepsilon'')^{d+r}$
of $\widetilde{\psi}$ conjugates \eqref{syslinagg} to \eqref{eq:1bis}. 
Consequently,
for any $\overline{V}\in(-\varepsilon'',\varepsilon'')^{r}$,
we may apply Proposition \ref{prop-01} to this restriction
and to the constant feedbacks
$\alpha_{1}(Z)=V_0+\overline{V}$ and $\alpha_{2}(Z)=V_0$;
this yields that $\widetilde{\psi}_{\mathrm{I}}$, given by
\begin{displaymath}
(Z^{1},Z^{2})\ \;\mapsto\;\ (\widetilde{X}^{1},\widetilde{X}^{2})
\ \;=\;\ (\widetilde{\psi}^{1}(Z^{1}),\widetilde{\psi}^{2}(Z^{1},Z^{2})),
\end{displaymath}
maps every solution of 
\begin{equation}
\label{ss12left}
\dot{Z}^{1} \ =\  0 \ \ ,\ \ \ \ 
\dot{Z}^{2} \ =\  J_{r}^{s}\overline{V}
\end{equation}
that remains in $Z_0+(-\varepsilon'',\varepsilon'')^{d}$ to a solution of
\begin{equation}
\label{ss12right}
        \begin{array}{ll}
\hspace{-2em}\dot{\widetilde{X}}^{1}= 0  \ ,\ \ \ 
\dot{\widetilde{X}}^{2} & \!\!= 
\widetilde{f}^{2}(\widetilde{X}^{1},\widetilde{X}^{2},
\widetilde{\psi}^{3}(\widetilde{\varphi}^{1}(\widetilde{X}^{1}),
\widetilde{\varphi}^{2}(\widetilde{X}^{1},\widetilde{X}^{2}),V_0+\overline{V}))
\\&\ -\,\widetilde{f}^{2}(\widetilde{X}^{1},\widetilde{X}^{2},
\widetilde{\psi}^{3}(\widetilde{\varphi}^{1}(\widetilde{X}^{1}),
\widetilde{\varphi}^{2}(\widetilde{X}^{1},\widetilde{X}^{2}),V_0))
        \end{array}
\end{equation}
that remains in $\widetilde{\psi}_{\mathrm{I}}
(Z_0+(-\varepsilon'',\varepsilon'')^{d})$, and {\it vice versa} upon 
applying Proposition \ref{prop-01} in the other direction.

Integrating \eqref{ss12left} explicitly with initial condition
$Z(0)=Z_0$, we get that
\begin{displaymath}
        t\ \;\mapsto\ \;\left(
\begin{array}{l}
        \widetilde{\psi}^{1}(Z_0^{1})  \\
        \widetilde{\psi}^{2}(Z_0^{1},Z_0^{2}\,+\,tJ_{r}^{s}\overline{V})
\end{array}
\right)
\end{displaymath}
solves \eqref{ss12right} for sufficiently small $t$, hence
$\widetilde{\psi}^{2}(Z^1,Z^2)$ is differentiable at  $Z_0$
with respect to its second argument in the direction 
$J_r^s\overline{V}$, with directional derivative 
\begin{eqnarray}
\nonumber
\frac{\partial\widetilde{\psi}^{2}}{\partial Z^{2}}(Z_0^{1},Z_0^{2})\,
        J_{r}^{s}\overline{V} & = &
        \widetilde{f}^{2}(\widetilde{\psi}^{1}(Z_0^{1}),
        \widetilde{\psi}^{2}(Z_0^{1},Z_0^{2}),
                 \widetilde{\psi}^{3}(Z_0^{1},Z_0^{2},V_0+\overline{V}))
\\\label{ddpsi2} &&
-\,        \widetilde{f}^{2}(\widetilde{\psi}^{1}(Z_0^{1}),
                 \widetilde{\psi}^{2}(Z_0^{1},Z_0^{2}),
                 \widetilde{\psi}^{3}(Z_0^{1},Z_0^{2},V_0))\ .
\end{eqnarray}
In particular, since $Z_0$ can be any member 
of $\widetilde{\varphi}_{\mathrm{I}}((-\varepsilon',\varepsilon')^{d})$
while $J_r^s\overline{V}$ can be assigned arbitrarily in 
$(-\varepsilon'',\varepsilon'')^{s}$, we conclude that 
$\partial\widetilde{\psi}^{2}/\partial Z^{2}(Z^{1},Z^{2})$ 
exists and is continuous since this holds for the partial derivatives.
Next we prove that $\partial\widetilde{\psi}^{2}/\partial Z^{2}$
is invertible at every point
by showing that its kernel reduces to zero. 
In fact, if the left-hand side of \eqref{ddpsi2} vanishes, so does the 
right-hand side 
which is also the value of the right-hand side of \eqref{ss12right} for 
$\widetilde{X}=\widetilde{\psi}_{\mathrm{I}}(Z_0)$. Therefore the constant map
$t\mapsto \widetilde{\psi}_{\mathrm{I}}(Z_0)$ is a solution to 
\eqref{ss12right} over a suitable time interval, and by conjugation
the constant map $t\mapsto Z_0$ is a solution to \eqref{ss12left}
over that time interval which clearly entails $J_r^s\overline{V}=0$, as
desired. Now, since $\partial\widetilde{\psi}^{2}/\partial Z^{2}$
is invertible at every $(Z^1,Z^2)\in
\widetilde{\varphi}_{\mathrm{I}}((-\varepsilon',\varepsilon')^{d})$, the
triangular structure of \eqref{phipsitilde} and the inverse function theorem 
together imply that 
\begin{equation}
\label{matderinv}
\frac{\partial\widetilde{\varphi}^{2}}{\partial \widetilde{X}^{2}}
(\widetilde{X}^{1},\widetilde{X}^{2})=
\left(\frac{\partial\widetilde{\psi}^{2}}{\partial
  Z^{2}}(\widetilde{\varphi}^{1}(\widetilde{X}^1),\widetilde{\varphi}^{2}(\widetilde{X}^1,\widetilde{X}^{2}))\right)^{-1}
\end{equation}
continuously exists and is invertible for
$(\widetilde{X}^1,\widetilde{X}^{2})\in (-\varepsilon',\varepsilon')^{d}$.
This proves point \ref{it-phi2diff}.

Let us turn to point \ref{it-f2}. Select an open neighborhood 
$\mathcal{W}$ of $0$ having compact closure in 
$(-\varepsilon',\varepsilon')^d$, so there is $\eta>0$ such that
$\widetilde{\varphi}(\widetilde{X},0)+(-\eta,\eta)^{d+r}$ is included
in $\widetilde{\varphi}((-\varepsilon',\varepsilon')^{d+r})$
whenever $\widetilde{X}\in\mathcal{W}$. If $\overline{V}\in (-\eta,\eta)^r$,
we can apply \eqref{ddpsi2} to 
$(Z_0,V_0)=\widetilde{\varphi}(\widetilde{X},0)$ with
$\widetilde{X}\in\mathcal{W}$, and we obtain in view of \eqref{matderinv}~:
\begin{eqnarray}
\label{ddpsi2linU}
\left(\frac{\partial\widetilde{\varphi}^{2}}{\partial \widetilde{X}^{2}}
(\widetilde{X}^{1},\widetilde{X}^{2})\right)^{-1}\,
        J_{r}^{s}\overline{V} & = &
        -\,        \widetilde{f}^{2}(\widetilde{X}^{1},
                 \widetilde{X}^{2}, 0)\\
\nonumber &&
\hspace{-10em}+\,\widetilde{f}^{2}\left(\widetilde{X}^{1},
        \widetilde{X}^{2},
                 \widetilde{\psi}^{3}\bigl(\widetilde{\varphi}^1(\widetilde{X}^1),
\widetilde{\varphi}^2(\widetilde{X}^1,\widetilde{X}^2),\widetilde{\varphi}^3
(\widetilde{X}^1,\widetilde{X}^2,0)+\overline{V}\bigr)\right)\,.
\end{eqnarray}
Set 
\begin{equation}
\label{defUb}
U=\widetilde{\psi}^{3}(\widetilde{\varphi}^1(\widetilde{X}^1),
\widetilde{\varphi}^2(\widetilde{X}^1,\widetilde{X}^2),\widetilde{\varphi}^3
(\widetilde{X}^1,\widetilde{X}^2,0)+\overline{V})
\end{equation}
and observe that $(\widetilde{X},\overline{V})\mapsto(\widetilde{X},U)=\widetilde{\psi}(\widetilde{\varphi}(\widetilde{X},0)+(0,\overline{V}))$
defines a continuous map 
$h:\mathcal{W}\times(-\eta,\eta)^{r}\to(-\varepsilon',\varepsilon')^{d+r}$,
such that $h(0)=0$, which is injective. By invariance of the domain, $h$
is a homeomorphism onto some open neighborhood of $0$, say
$\mathcal{N}\subset(-\varepsilon',\varepsilon')^{d+r}$.
For $(\widetilde{X},U)\in\mathcal{N}$, \eqref{defUb} can be inverted as
\begin{equation}
\label{invB}
\overline{V}=\widetilde{\varphi}^3(\widetilde{X},U)-\widetilde{\varphi}^3
(\widetilde{X},0),
\end{equation}
and substituting  \eqref{defUb} and \eqref{invB}
in \eqref{ddpsi2linU} yields \eqref{relf2}.

Finally we prove point \ref{it-conjred}, keeping in mind 
the previous definitions and properties of $h$, $\mathcal{W}$, $\eta$ and 
$\mathcal{N}$.
For $\widetilde{X}=(\widetilde{X}^{1},\widetilde{X}^{2})\in(-\varepsilon',\varepsilon')^d$,
define $\overline{V}(\widetilde{X})\in\RR^s\times\{0\}\subset\RR^r$ by the formula~: 
\begin{equation}
\label{relf2pa}
J_r^s\overline{V}(\widetilde{X})=
\frac{\partial\widetilde{\varphi}^{2}}{\partial\widetilde{X}^{2}}
                     (\widetilde{X}^{1},\widetilde{X}^{2})
                \bigl(\widetilde{f}^{2}(0,0,0)-\widetilde{f}^{2}(\widetilde{X}^{1},\widetilde{X}^{2},0)\bigr).
\end{equation}
Clearly $\overline{V}:(-\varepsilon',\varepsilon')^d\to\RR^r$ is continuous
and $\overline{V}(0)=0$, so there exists an open neighborhood
$\mathcal{V}\subset\mathcal{W}$ of $0$ in $\RR^{d}$ such that 
$\overline{V}(\widetilde{X})\in(-\eta,\eta)^r$ as soon as
$\widetilde{X}\in\mathcal{V}$; then, if we set 
$h(\widetilde{X},\overline{V}(\widetilde{X}))=(\widetilde{X},U(\widetilde{X}))\in\mathcal{N}$,
it follows from \eqref{relf2pa}, \eqref{invB}, and  \eqref{relf2} that
\begin{equation}
\label{constU}
\widetilde{f}^{2}(\widetilde{X}^{1},\widetilde{X}^{2},U(\widetilde{X}))
=\widetilde{f}^{2}(0,0,0),\ \ \ \ \ \ \widetilde{X}\in\mathcal{V}.
\end{equation}
\emph{We will show, using Proposition
{\rm \ref{prop-uspecial}}, that the restriction of $\widetilde{\varphi}_{\mathrm{I}}$ to any relatively compact open subset $\mathcal{X}$ of $\mathcal{V}$
conjugates \eqref{sysredX} and
\eqref{sysredZ} over $\mathcal{X}$, $ \widetilde{\varphi}(\mathcal{X})$,
and this will achieve the proof}.
To this effect, let $\mathcal{C}$
to be the collection of all piecewise affine maps $\RR\to \RR^{s}$ 
with constant slope $\widetilde{f}^{2}(0,0,0)$ ({\it cf} the discussion before
Proposition \ref{prop-uspecial}) and note that, for any open 
set $\mathcal{O}\subset\RR^s$ and any compact interval $J\subset\RR$,
the restriction of $\mathcal{C}$ to $J$ contains, 
in its uniform closure, the set all piecewise continuous maps
$J\to\mathcal{O}$. Now, consider a solution 
$\gamma~:I\to \mathcal{V}$ of the control system~:
\begin{eqnarray}
  \label{sysredXp}
\dot{\widetilde{X}}^{1} &=&
        \widetilde{f}^{1}(\widetilde{X}^{1},\Upsilon)
\end{eqnarray}
with state $\widetilde{X}^{1}$ and control $\Upsilon$;
hereafter,
$\mathcal{V}_{\mathrm{I}}\subset\RR^{d-s}$ and
$\mathcal{V}_{\mathrm{I\!I}}\subset\RR^s$ will indicate the projections of 
$\mathcal{V}$ onto the first $d-s$ and the last $s$ components respectively,
and similarly for any other open set in $\RR^d$.
Assume that the control function 
$\gamma_{\mathrm{I\!I}}~:I\to\mathcal{V}_{\mathrm{I\!I}}$ 
is the restriction to 
$I$ of some member of $\mathcal{C}$.
By definition, if $a,b$ are the endpoints of
$I$ (that may belong to $I$ or not), there are time instants 
$a=t_0<t_1<\cdots<t_N=b$,
and vectors $\bar{\xi}_1,\ldots,\bar{\xi}_N\in\RR^{s}$ such that,
for $1\leq j<N$, one has
\begin{equation}
\label{pwaffine}
t_{j-1}<t<t_{j}\ \;\Rightarrow\;\ 
\gamma_{\mathrm{I\!I}}(t)\;=\;\bar{\xi}_{j}+t\widetilde{f}^{2}(0,0,0),
\end{equation}
while at the points $t_j$ themselves 
$\gamma_{\mathrm{I\!I}}$ is either right or left continuous when $1<j<N$.
\emph{We claim} that $\widetilde{\varphi}_{\mathrm{I}}(\gamma(t))$ 
is a solution that
remains in $\widetilde{\varphi}_{\mathrm{I}}(\mathcal{V})$ 
of the control system~:
\begin{eqnarray}
  \label{sysredZp}
\dot{Z}^{1} &=&
        g^{1}(Z^{1},\Gamma)
\end{eqnarray}
with state $Z^{1}$ and control $\Gamma$. In fact, since
$\gamma_{\mathrm{I}}$ is continuous by definition of a solution,
so is $\widetilde{\varphi}^{1}(\gamma_{\mathrm{I}})$
and therefore, as
$\widetilde{\varphi}_{\mathrm{I}}(\gamma(t))$ lies in 
$\widetilde{\varphi}_{\mathrm{I}}(\mathcal{V})$ for all $t\in I$ 
by construction,
it is enough to check that
\begin{equation}
\label{eqdir}
\widetilde{\varphi}^{1}(\gamma_{\mathrm{I}}(T_{2}))-
\widetilde{\varphi}^{1}(\gamma_{\mathrm{I}}(T_{1}))\!=
        \!\!\int_{T_{1}}^{T_{2}}
\!\!g^1\!\!\left(\widetilde{\varphi}^{1}(\gamma_{\mathrm{I}}(t)),
\widetilde{\varphi}^{2}(\gamma_{\mathrm{I}}(t),\gamma_{\mathrm{I\!I}}(t))\right)\Dt\ 
\ 
\end{equation}
whenever $t_{j-1}<T_1<T_2<t_j$ for some $j>1$. However, the restriction
of $\gamma(t)$ to $(t_{j-1},t_j)$ is a solution that remains in 
$\mathcal{V}$ of the differential equation~:
$$
  \begin{array}{rcl}
\dot{\gamma}_{\mathrm{I}} &=&
\widetilde{f}^{1}(\gamma_{\mathrm{I}},\gamma_{\mathrm{I\!I}}) \\
\dot{\gamma}_{\mathrm{I\!I}} &=&
\widetilde{f}^{2}(0,0,0),
  \end{array}
$$
hence $(\gamma(t), U(\gamma(t))$ is, by \eqref{constU}, a solution 
of \eqref{eq:1bis} that remains in $\mathcal{N}$, and therefore
\eqref{eqdir} follows from the triangular structure
\eqref{phipsitilde} of $\widetilde{\varphi}$ and the fact that
it conjugates system \eqref{eq:1bis} to system \eqref{syslinagg}.
\emph{This proves the claim}.

In the other direction, we observe since it is included in $\mathcal{W}$
that $\mathcal{V}$ has compact closure in $(-\varepsilon',\varepsilon')^d$,
and therefore that $\widetilde{\varphi}_{\mathrm{I}}(\mathcal{V})$ 
in turn has compact closure in $\widetilde{\varphi}_{\mathrm{I}}
\bigl((-\varepsilon',\varepsilon')^d\bigr)$. Pick
$\eta'>0$ such that
$\widetilde{\varphi}_{\mathrm{I}}(\mathcal{V})\times(-\eta',\eta')^{r}\subset\widetilde{\varphi}((-\varepsilon',\varepsilon')^{d+r})$, 
and let $\mathcal{C}'$ denote the
collection of all piecewise smooth maps $\RR\to \RR^{s}$ whose derivative
is strictly bounded by $\eta'$ component-wise. The 
restriction of $\mathcal{C}'$ to any compact real interval $J$ is
uniformly dense in the set all piecewise continuous maps
$J\to\mathcal{O}$, for any open set $\mathcal{O}\subset\RR^s$.
Clearly, any solution
$\gamma':I\to\widetilde{\varphi}_{\mathrm{I}}(\mathcal{V})$ of 
system \eqref{sysredZp}, whose control function 
$\gamma'_{\mathrm{I\!I}}:I\to\bigl(\widetilde{\varphi}_{\mathrm{I}}(\mathcal{V})\bigr)_{\mathrm{I\!I}}$ 
is the restriction to $I$ of
some member of $\mathcal{C}'$, satisfies the differential equation
$$
  \begin{array}{rcl}
\dot{\gamma'}_{\mathrm{I}} &=&
g^{1}(\gamma'_{\mathrm{I}},\gamma'_{\mathrm{I\!I}}) \\
\dot{\gamma}'_{\mathrm{I\!I}} &=&
J_r^s\left(d{\gamma}'_{\mathrm{I\!I}}/dt\,,\,0\right)
  \end{array}
$$
on every interval where it is smooth. By the very definition of $\eta'$ and
$\mathcal{C}'$, it follows that
$\bigl(\gamma'(t),(d\gamma'_{\mathrm{I\!I}}(t)/dt,0)\bigr)$ is,
on such intervals, a solution to 
\eqref{syslinagg} that remains in
$\widetilde{\varphi}((-\varepsilon',\varepsilon'))^{d+r}$ and, since
$\widetilde{\psi}$ conjugates system \eqref{syslinagg} 
to system \eqref{eq:1bis}, we argue as before to the effect that
$\widetilde{\psi}_{\mathrm{I}}(\gamma')$ is a solution to 
system \eqref{sysredXp} that remains in $\mathcal{V}$.
Appealing to Proposition \ref{prop-uspecial}, we conclude that
$\widetilde{\varphi}_{\mathrm{I}}$ conjugates system \eqref{sysredXp}
to system \eqref{sysredZp} on relatively compact open subsets of $\mathcal{V}$,
as desired.
\end{proof}

\appendix

\section{Four lemmas on ODEs}
\label{blode1}

Throughout this section, we let $\mathcal{U}$ be an open subset of $\RR^d$.
We say that a continuous vector field $X:\mathcal{U}\to\RR^d$
\emph{has a flow} if the
Cauchy problem $\dot{x}(t)=X(x(t))$  
with initial condition $x(0)=x_0$ has a unique solution, defined for
$t\in(-\varepsilon,\varepsilon)$ with $\varepsilon=\varepsilon(x_0)>0$.
The flow of $X$ at time $t$ is denoted by $X_t$, in other words
we have with the preceding notations that $X_t(x_0)=x(t)$.
It is easy to see that the domain of definition of 
$(t,x)\mapsto X(t,x)$ is open in $\RR\times\mathcal{U}$.
\begin{lem}
\label{lem-ode2}
If $X:\mathcal{U}\to\RR^d$ is a continuous vector field that has a flow, 
the map 
$(t,x)\mapsto X_t(x)$ is continuous on the open subset of 
$\RR\times\mathcal{U}$ where it is defined.
\end{lem}
\begin{proof}
This is an easy consequence of the Ascoli-Arzela theorem, and actually a
special case of \cite[chap. V, Theorem 2.1]{Hart82}. 
\end{proof}
\begin{lem}
\label{lem-ode3}
Assume that the sequence of continuous vector fields 
$X^{k}:\mathcal{U}\to\RR^d$
converges to $X$, uniformly on compact subsets of $\mathcal{U}$,
and that all the $X^{k}$ as well as $X$ itself have a flow. 
Suppose that $X_t(x)$ is defined for all $(t,x)\in[0,T]\times K$ with
$T>0$ and $K\subset\mathcal{U}$ compact. 
Then $X^{k}_t(x)$ is also defined on $[0,T]\times K$ for $k$ large
enough, and the sequence of mappings $(t,x)\mapsto X^{k}_t(x)$ 
converges to $(t,x)\mapsto X_t(x)$, uniformly on $[0,T]\times K$.
\end{lem}
\begin{proof}
By assumption,
\[K_1=\{X_t(x);~(t,x)\in[0,T]\times K\}\]
is a well-defined subset of  $\mathcal{U}$ that contains $K$,
and it is compact by Lemma \ref{lem-ode2}.
Let $K_0$ be another compact subset of $\mathcal{U}$ whose interior
contains $K_1$, and put $d(K_1,\mathcal{U}\setminus K_0)=\eta>0$
where $d(E_1,E_2)$ indicates the distance between two 
sets $E_1$, $E_2$. From the hypothesis there is $M>0$ such that
$\|X^k\|\leq M$ on $K_0$ for all $k$, hence the maximal solution to
$\dot{x(t)}=X^k(x(t))$ with initial condition $x(0)=x_0\in K$
remains in $K_0$ as long as $t\leq \eta/2M$. Consequently the flow 
$(t,x)\mapsto X^{k}_t(x)$ is defined  on 
$[0,\eta/2M]\times K$ for all $k$, with values in $K_0$.
\emph{We claim} that it is a bounded 
equicontinuous sequence of functions there. Boundedness is clear since these
functions are $K_0$-valued, so we must show that, to every
$(t,x)\in[0,\eta/2M]\times K$ and every $\varepsilon>0$, there is
$\alpha>0$ such that $\|X^k(t',x')-X^k(t,x)\|<\varepsilon$ for all $k$
as soon as $|t-t'|+\|x-x'\|<\alpha$. By the mean-value theorem and the 
uniform majorization $\|X^k(X^k_t(x))\|\leq M$, it is sufficient to prove this when
$t=t'$. Arguing by contradiction, assume for some subsequence $k_l$
and some sequence $x_l$ converging to $x$ in $K$ that
\begin{equation}
\label{ecartX}
\|X_t^{k_l}(x)-X_t^{k_l}(x_l)\|\geq\varepsilon~~~~\mathrm{for~all~} l\in\NN.
\end{equation}
Then, by Lemma \ref{lem-ode2}, the index $k_l$ tends to infinity with $l$.
Next consider the sequence of maps $F_l:[0,\eta/2M]\to K_0$ defined by
$F_l(t)=X^{k_l}_t(x_l)$. Again, by the mean value theorem, it is
a bounded equicontinuous family of functions and,
by the Ascoli-Arzela theorem, it is relatively
compact in the topology of uniform convergence
(compare \cite[chap. II, Theorem 3.2]{Hart82}).
But if $\Phi:[0,\eta/2M]\to K_0$ is the uniform limit of some 
subsequence $F_{l_j}$, and
since $X^{k_{l_j}}$ converges uniformly to $X$ 
on $K_0$ as $j\to\infty$,
taking limits in the relation
\[X^{k_{l_j}}_t(x_{l_j})=x_{l_j}+\int_0^t X^{k_{l_j}}(X^{k_{l_j}}_s(x_{l_j}))
\,ds\]
 gives us
\[\Phi(t)=x+\int_0^t X(\Phi(s))\,ds\]
so that $\Phi(t)=X_t(x)$ since $X$ has a flow. Altogether
$F_l(t)$ converges uniformly to $X_t(x)$ on $[0,\eta/2M]$ because
this is the only accumulation point, and then \eqref{ecartX} becomes absurd.
\emph{This proves the claim}. From the claim it follows, using
the Ascoli-Arzela theorem again, that the family of functions
$(t,x)\mapsto X^{k}_t(x)$ is relatively compact for the topology of
uniform convergence $[0,\eta/2M]\times K\to K_0$,
and in fact it converges to $(t,x)\mapsto X_t(x)$
because, by the same limiting argument as was used to prove the claim,
every accumulation point $\Phi(t,x)$ 
must be a solution to
\[\Phi(t,x)=x+\int_0^t X(s,\Phi(s,x))\,ds\]
hence for fixed $x$ is an integral curve of $X$ with initial condition
$x$. In particular, by definition of $K_1$, we shall have that
$d(X^k_t(x),K_1)<\eta/2$ for all $(t,x)\in[0,\eta/2M]\times K$
as soon as $k$ is large enough. For such $k$ the flow 
$(t,x)\mapsto X^k_t(x)$ will be defined on $[0,\eta/M]\times K$
with values in $K_0$, and we can repeat the whole argument again
to the effect that $X^k_t(x)$ converges uniformly to $X_t(x)$ there.
Proceeding inductively, we obtain after $[2TM/\eta]+1$ steps at most
that $(t,x)\mapsto X^{k}_t(x)$ is defined  on 
$[0,T]\times K$ with values in $K_0$ for $k$ large enough, and converges
uniformly to $(t,x)\mapsto X_t(x)$ there, as was to be shown.
\end{proof}
The next lemma stands analogous to Lemma \ref{lem-ode3} for time-dependent
vector fields, assuming that the convergence holds boundedly
almost everywhere in time. The assumption that the vector fields
have a flow is replaced here by a local Lipschitz condition that we now
comment upon. 

By definition, a time-dependent vector field
$X:[t_1,t_2]\times\mathcal{U}\to\RR^d$ is locally Lipschitz
with respect to the second variable if every 
$(t_0,x_0)\in[t_1,t_2]\times\mathcal{U}$ 
has a neighborhood there such that
$\|X(t,x')-X(t,x)\|< c \|x'-x\|$, for some constant $c$, whenever
$(t,x)$ and $(t,x')$ belong to that neighborhood. This of course
entails that $X$ is bounded on compact subsets of
$[t_1,t_2]\times\mathcal{U}$. Next,
by the compactness of $[t_1,t_2]$, the local Lipschitz character of $X$ 
strengthens to the effect that each $x_0\in\mathcal{U}$ has a
neighborhood $\mathcal{N}_{x_0}$ such that
$\|X(t,x')-X(t,x)\|< c_{x_0} \|x'-x\|$, for some constant $c_{x_0}$, whenever
$x,x'\in\mathcal{N}_{x_0}$ and $t\in[t_1,t_2]$. 
If now $\mathcal{K}\subset\mathcal{U}$ is compact, we can cover it by finitely
many $\mathcal{N}_{x_{0,k}}$ as above and find $\varepsilon>0$ such that
$x,x'\in\mathcal{K}$ and $\|x-x'\|<\varepsilon$ is impossible 
unless $x,x'$ lie in some common $\mathcal{N}_{x_0}$. 
Consequently there is 
$c_{\mathcal{K}}>0$ such that $\|X(t,x')-X(t,x)\|< c_{\mathcal{K}} \|x'-x\|$
whenever $x,x'\in\mathcal{K}$ and $t\in[t_1,t_2]$, because if
$\|x-x'\|<\varepsilon$ we can take  $c_{\mathcal{K}}\geq \max_k c_{x_{0,k}}$,
whereas if $\|x-x'\|\geq\varepsilon$ it is enough to take
$c_{\mathcal{K}}> 2M/\varepsilon$ where $M$ is a bound for
$\|X\|$ on $[t_1,t_2]\times\mathcal{K}$. Finally, if $X(t,x)$ 
happens to vanish identically for $x$ outside some compact
$\mathcal{K}'\subset\mathcal{U}$, we can choose $\mathcal{K}$ such that
$$\mathcal{K}'\ \;\subset\;\ \stackrel{\circ}{\mathcal{K}}
\ \;\subset\;\ \mathcal{K}\ \;\subset\;\ \mathcal{U}
$$
and construct $c_{\mathcal{K}}$ as before except that we also
pick $\varepsilon>0$ so small that $\|x-x'\|<\varepsilon$ is 
impossible for $x\in\mathcal{K}'$ and $x'\notin\mathcal{K}$. 
Then it holds that  $\|X(t,x')-X(t,x)\|< c_{\mathcal{K}} \|x'-x\|$
for all $x,x'\in\mathcal{U}$ and all $t\in[t_1,t_2]$, that is to say
$X(t,x)$ becomes globally Lipschitz with respect to $x$.
These remarks will be used in the proof to come.
\begin{lem}
\label{lem-ode1}
Let $t_1<t_2$ be two real numbers and 
$X^{k}:[t_1,t_2]\times\mathcal{U}\to\RR^d$
a sequence of time-dependent vector fields, measurable with respect to
$t$, locally Lipschitz continuous with respect to 
$x\in\mathcal{U}$, and bounded on compact subsets of
$[t_1,t_2]\times\mathcal{U}$ independently of $k$.
Let $X:[t_1,t_2]\times\mathcal{U}\to\RR^d$
be another time-dependent vector field, measurable with respect to
$t$, locally Lipschitz continuous with respect to 
$x\in\mathcal{U}$, and assume that, to each compact
$\mathcal{K}\subset\mathcal{U}$, there is 
$E_{\mathcal{K}}\subset[t_1,t_2]$ of zero measure
such that, whenever $t\notin E_{\mathcal{K}}$, 
the sequence $X^{k}(t,x)$ converges to
$X(t,x)$ as $k\to\infty$, uniformly with respect to $x\in\mathcal{K}$.
Suppose finally that $\gamma:[t_1,t_2]\to\mathcal{U}$ is, for
some $(t_0,x_0)\in [t_1,t_2]\times\mathcal{U}$, a solution 
to the Cauchy problem 
\begin{equation}
\label{repXl}
\dot{\gamma}(t)\ =\ X(t,\gamma(t))\ ,\ \ \ \ \gamma(t_0)\ =\ x_0.
\end{equation}
Then, for $k$ large enough, there is a unique solution 
$\gamma_k:[t_1,t_2]\to\mathcal{U}$ to the Cauchy problem
\begin{equation}
\label{equamk}
\dot{\gamma}_k(t)\ =\ X^k(t,\gamma_k(t))\ ,\ \ \ \ \gamma_k(t_0)\ =\ x_0,
\end{equation}
and the sequence $(\gamma_k)$ converges to $\gamma$, 
uniformly on $[t_1,t_2]$.
\end{lem}
\begin{proof}
Upon multiplying $X^k(t,x)$ and $X(t,x)$ by a smooth function $\varphi(x)$
which is compactly supported $\mathcal{U}\to\RR$ 
and identically $1$ on a
neighborhood of $\gamma([t_1,t_2])$, we may assume
in view of the discussion preceding the lemma that
$X(t,x)$ and $X^k(t,x)$ are defined and bounded 
$[t_1,t_2]\times\RR^d\to\RR^d$ independently 
of $k$, measurable with respect to $t$, and (globally)
Lipschitz continuous with respect to $x$.

Then, by classical results \cite[Proposition C 3.8., Theorem 54]{Sont98},  
the solution to \eqref{equamk}, say $\gamma_k$ 
uniquely exists $[t_1,t_2]\to\RR^d$ for each $k$~:
\begin{equation}
\label{limrepg}
\gamma_k(t)=x_0+\int_{t_0}^t X^k(s,\gamma_k(s))\,ds,~~~~~~t\in[t_1,t_2].
\end{equation}
From the boundedness of $X^k$, it is clear that $\gamma_k$ is an 
equicontinuous and bounded family of functions,
hence it is relatively
compact in the topology of uniform convergence on 
$[t_1,t_2]$. All we have to prove then is that
every accumulation point of $\gamma_k$ coincides with $\gamma$.
Extracting a subsequence if necessary, let us assume that
$\gamma_k$ converges to some $\bar{\gamma}$, uniformly on $[t_1,t_2]$.
Let $\mathcal{K}\subset\RR^d$ be a compact set containing 
$\gamma_k([t_1,t_2])$ for all $k$; such a set exists by
the boundedness of $\gamma_k$. If we let $E_{\mathcal{K}}\subset[t_1,t_2]$ 
be the set of zero measure granted by the hypothesis, there exists to
each $s\in[t_1,t_2]\setminus E_{\mathcal{K}}$ and each $\varepsilon>0$
an integer $k_{s,\varepsilon}$ such that
$\|X^k(s,x)-X(s,x)\|<\varepsilon$ as soon as $x\in\mathcal{K}$ and 
$k>k_{s,\varepsilon}$. In another connection, 
the Lipschitz character of $X$ with respect to the second argument and the
uniform convergence of $\gamma_k$ to $\bar{\gamma}$ shows that
that $\|X(s,\gamma_k(s))-X(s,\bar{\gamma}(s))\|<\varepsilon$ for $k$ large
enough. Altogether, by a 2-$\varepsilon$ majorization , we find that 
$$\lim_{k\to\infty}\|X^k(s,\gamma_k(s))-X(s,\bar{\gamma}(s))\|=0,$$
that is to say the integrand in the right-hand side of \eqref{limrepg}
converges point-wise almost everywhere to $X(s,\bar{\gamma}(s))$. Since $X^k$
is bounded we can apply the dominated convergence theorem and, taking limits
on both sides of \eqref{limrepg} as $k\to\infty$, 
we find that $\bar{\gamma}$ is a solution
to \eqref{repXl} whereas the latter is unique. Hence $\bar{\gamma}-\gamma$ as
desired. 
\end{proof}

The following averaging lemma for continuous vector fields is less classical
than in the locally Lipschitz case, where the Cauchy problem has a unique 
solution. 
\begin{lem}
\label{lem-aver}
Let $t_1<t_2$ be real numbers and $(X^{1,\ell})_{\ell\in\NN}$,
$(X^{2,\ell})_{\ell\in\NN}$, be two sequences of 
continuous time-dependent vector fields
$[t_1,t_2]\times\RR^d\rightarrow\RR^d$, uniformly bounded 
with respect to $\ell$,
that converge uniformly on compact subsets of $[t_1,t_2]\times\RR^d$ 
to some vector fields $X^1$ and $X^2$ respectively. Denoting by 
$L=t_2-t_1$ the length of the time interval, define,
for each $\ell\in\NN$, the ``average'' vector field
$G_\ell:[t_1,t_2]\times\RR^d\rightarrow\RR^d$ by~:
\begin{equation}
    \label{eq:a2}
    \begin{array}{l}
t\in [t_1+\frac j \ell L\,,\,t_1+\frac{2j+1}{2\ell} L )
    \ \Rightarrow\ 
    G_\ell(t,x)\ =\ X^{1,\ell}(t,x)\;,
\\[0.2em]
t\in [t_1+\frac{2j+1}{2\ell} L\,,\,t_1+\frac{j+1}{\ell} L )
    \ \Rightarrow\ 
    G_\ell(t,x)\ =\ X^{2,\ell}(t,x)\ \ ,
    \end{array}
\end{equation}
for $j\in\{0,\ldots,\ell-1\}$ and, say, $G_\ell(t_2,x)=X^{2,\ell}(t_2,x)$ for definiteness.

Let $\gamma_\ell:[t_1,t_2]\rightarrow\RR^d$ be a
solution to 
\begin{equation}
  \label{eq:eqn}
  \gamma_\ell(t)\,-\,\bar{x} \ \;=\;\ \int_{t_1}^t
  G_\ell(\tau,\gamma_\ell(\tau)) \D\tau\ .
\end{equation}
Then the sequence $(\gamma_\ell)$ is compact in $\CCC0([t_1,t_2],\RR^d)$, and every
accumulation point $\gamma_\infty$ is a solution to
\begin{equation}
  \label{eq:eqlim}
  \gamma_\infty(t)\,-\,\bar x \ \;=\;\ \frac 1 2 \int_{t_1}^t
\left(X^1(\tau,\gamma_\infty(\tau))
+X^2(\tau,\gamma_\infty(\tau))\right)\D\tau\ .
\end{equation}
\end{lem}
\begin{proof}
Let
\begin{equation}
\label{defM}
  M\ \;=\;\ \sup_{t,x,i,\ell} \|X^{i,\ell}(t,x)\|\ .
\end{equation}
From \eqref{eq:a2}-\eqref{eq:eqn}, it is clear that $M$ is a Lipschitz
constant for $\gamma_\ell$, regardless of $\ell$.  In particular $\gamma_\ell(t)$ stays
in a fixed compact ball $B$ of radius $ML$, and the family
$(\gamma_\ell)$ is equicontinuous.
From Ascoli-Arzela's theorem this implies compactness of the sequence 
$(\gamma_\ell)$ in the uniform topology on $[t_1,t_2]$.

Rewrite \eqref{eq:eqn} as
\begin{eqnarray}
\nonumber
  \gamma_\ell(t)\,-\,\bar{x} &=&\int_{t_1}^t
\left( G_\ell(\tau,\gamma_\ell(\tau))\,-\,
\frac { X^{1,\ell}(\tau,\gamma_\ell(\tau))+X^{2,\ell}(\tau,\gamma_\ell(\tau))}{2}
\right)\D\tau
\\
\nonumber
  &&\hspace{-5em}+\;\int_{t_1}^t
\left( \frac { X^{1,\ell}(\tau,\gamma_\ell(\tau))+X^{2,\ell}(\tau,\gamma_\ell(\tau))}{2}
\,-\,\frac { X^1(\tau,\gamma_\ell(\tau))+X^2(\tau,\gamma_\ell(\tau))}{2}
\right)\D\tau
\\
\label{bof}
&&\hspace{-5em}+\;\int_{t_1}^t \frac{
X^1(\tau,\gamma_\ell(\tau))+X^2(\tau,\gamma_\ell(\tau))}{2}\,\D\tau\ .
\end{eqnarray}
By the uniform convergence of $X^{i,\ell}$ to $X^i$, it will clearly follow that
any accumulation point $\gamma_\infty$ of $(\gamma_\ell)$ satisfies
\eqref{eq:eqlim} if only we can show that the first integral in the right-hand
side of \eqref{bof} converges to zero as $\ell\rightarrow\infty$.

To prove this, we compute, from the definition of $G_\ell$~:
\begin{eqnarray}
\nonumber
&&
\int_{t_1+\frac j \ell L}^{t_1+\frac {j+1} \ell L}
\left( G_\ell(\tau,\gamma_\ell(\tau))\,-\,
\frac { X^{1,\ell}(\tau,\gamma_\ell(\tau))+X^{2,\ell}(\tau,\gamma_\ell(\tau))}{2}
\right)\D\tau
\\
\label{nunu}
&&=\ \;
\int_{t_1+\frac j \ell L}^{t_1+\frac{2j+1}{2\ell}L}
\frac { X^{1,\ell}(\tau,\gamma_\ell(\tau))-X^{2,\ell}(\tau,\gamma_\ell(\tau))}{2}\,\D\tau
\\
\nonumber
&&\ \ \ \ \ \ \ \ \ \ 
\;-\;
\int_{t_1+\frac{2j+1}{2\ell}L}^{t_1+\frac {j+1} \ell L}
\frac{ X^{1,\ell}(\tau,\gamma_\ell(\tau))-X^{2,\ell}(\tau,\gamma_\ell(\tau))}{2}\,\D\tau
\\
\nonumber
&&=
\int_{t_1+\frac j \ell L}^{t_1+\frac{2j+1}{2\ell}L}
\left(
\Delta_\ell(\tau,\gamma_\ell(\tau))-\Delta_\ell
(\tau+{\textstyle\frac{L}{2\ell}},\gamma_\ell(\tau+{\textstyle\frac{L}{2\ell}}))
\right)\D\tau
\end{eqnarray}
with $\Delta_\ell=\frac 1 2 (X^{1,\ell}-X^{2,\ell})$.
On the compact set $[t_1,t_2]\times B$, the vector field $\Delta_\ell$ is
uniformly continuous with a modulus of continuity that does not depend on
$\ell$~; consequently, by the uniform Lipschitz property of $\gamma_\ell$, 
we see for
arbitrary $\varepsilon>0$ that the norm of the last integral is less that
$\varepsilon/2\ell$ as soon as $\ell$ is large enough, independently of $j$.

Now, the first integral in \eqref{bof} can be decomposed into a sum of at most
$\ell$ integrals like these we just studied plus an integral over an interval of
length smaller that $1/\ell$. Since the norm of the integrand is bounded by $2M$,
the norm of the last term is less than $2M/\ell$.  Summing over $j$, the above
estimates tell us that, for $t\in[t_1,t_2]$ and for $\ell$ is large enough,
\begin{displaymath}
  \int_{t_1}^t
\left( G_\ell(\tau,\gamma_\ell(\tau))\,-\,
\frac { X^1(\tau,\gamma_\ell(\tau))+X^2(\tau,\gamma_\ell(\tau))}{2}
\right)\D\tau\ \;\leq\;\ \frac{\varepsilon}{2}\;+\;\frac{2M}{\ell}\,.
\end{displaymath}
This achieves the proof since $\varepsilon>0$ was arbitrary.
\end{proof}

\section{Orbits of families of vector fields}
\label{app-orbits}

In the proof of lemma \ref{lem-recur1} we need results from
\cite{Suss73} on orbits\footnote{
One of the motivations in \cite{Suss73} was to generalize 
the notion of integral manifolds to vector fields that are smooth but
not real analytic.
Note that the orbits of a family of \emph{real analytic} vector fields actually
coincide with
the maximal integral manifolds of the closure of this family under
Lie brackets~\cite{Suss73,Lobr70,Naga66}.
However, even if we assume the control system \eqref{sysnl}
to be real analytic, integral manifolds are of no help to us
because topological conjugacy does
not preserve tangency nor Lie brackets.
Using orbits of families of vector fields instead
is much more efficient,
because topological conjugacy does preserve integral curves.}
of families of smooth vector fields, that were
recently exposed in the textbook \cite[chapter II]{Jurd97}. 
We recall them below, in a slightly expanded form.

Let $\mathcal{F}$ be a family of smooth vector fields defined on 
an open subset $U$ of $\RR^{d}$. 
For any positive integer $N$ and vector fields $X^{1},\ldots,X^{N}$
belonging to $\mathcal{F}$, given $m\in U$, consider the map $F$
given by
\begin{equation}
\label{deforb}
(t_{1},\ldots,t_{N})\ \;\mapsto\;\ 
        X^{1}_{t_{1}}(X^{2}_{t_{2}}(\cdots(X^{N}_{t_{N}}(m))\cdots))
\end{equation}
where the standard notation $X_t(x)$ indicates the flow of $X$ from $x$
at time $t$;
of course, $F$ depends on the choice of the vector fields $X^{j}$ and 
of the point $m$. 
This map is defined on some open connected neighborhood of the origin,
hereafter denoted by 
$\mathrm{dom}(F)$, and takes values in $U$. In fact,
$(t_{1},\ldots,t_{N})\in\mathrm{dom}(F)$ if, and only if, for every
$j\in\{1,\ldots,N\}$, the solution $x(\tau)$ to $\dot{x}=X^j(x)$, 
with initial condition
$x(0)= X^{j-1}_{t_{j-1}}(\cdots(X^{1}_{t_{1}}(m))\cdots)$, exists in $U$
for all $\tau\in[0,t_{j}]$ (or $[t_{j},0]$ if $t_{j}<0$).

The \emph{orbit} of the family $\mathcal{F}$ through a point $m\in U$ is the
set of all points that lie in the image of $F$ for at least one choice of the
vector fields $X^{1},\ldots,X^{N}$.
In words, the orbit of the family $\mathcal{F}$ through $m$ is the set of points
that may be linked to $m$ in $U$ upon concatenating finitely many integral
curves of vector fields in the family. 
We shall denote by $\mathcal{O}_{\mathcal{F},p}$
the orbit of $\mathcal{F}$ through $m$.

Note that the definition depends on $U$ in a slightly subtle manner~: 
if $\mathcal{F}$ defines by
restriction a family of vector fields $\mathcal{F}|_{V}$ on
a smaller open set $V\subset U$ and if $m\in V$, then
\begin{equation}
  \label{calObis}
  V\cap\mathcal{O}_{\mathcal{F},m}\ \;\supset\;\ 
\mathcal{O}_{\mathcal{F}|_{V},m},
\end{equation}
but the inclusion is generally strict because of the requirement that the
integral curves used to construct $\mathcal{O}_{\mathcal{F}|_{V},m}$
should lie entirely in $V$. 

We turn to topological considerations. 
The topology of $U$ is the usual Euclidean topology.
The topology of $\mathcal{O}_{\mathcal{F},m}$ \emph{as an orbit} is the finest
that makes all the maps $F$, arising from \eqref{deforb}, continuous on 
their respective domains of definition, the latter 
being endowed with the Euclidean topology.
The classical smoothness of the flow
implies that each $F$ is continuous $\mathrm{dom}(F)\to\RR^d$,
hence the topology  of $\mathcal{O}_{\mathcal{F},m}$ as an orbit  
is finer than the Euclidean topology induced by the ambient
space $U$. \emph{It can be strictly finer}, and this is why we speak
of the \emph{orbit topology}, as opposed to the \emph{induced topology}.

Starting from $\mathcal{F}$, one defines a larger family of vector fields
$P_{\mathcal{F}}$, consisting of all the push-forwards\footnote{
Recall that the push-forward of a vector field $X:V\to\RR^d$
through a diffeomorphism $\varphi:V\to\varphi(V)$ is the vector field
$\varphi_\star X$ on $\varphi(V)$ whose flow at each time is the conjugate of
the flow of $X$ under the diffeomorphism $\varphi$; it can be defined as
$\varphi_\star X(\varphi(x))=D\varphi(x)X(x)$, where $D\varphi(x)$ is the
derivative of $\varphi$ at $x\in V$.
} 
of vector fields in  $\mathcal{F}$ through all local 
diffeomorphisms of the form 
$X^{1}_{t_{1}}\circ X^{2}_{t_{2}}\circ\cdots\circ X^{N}_{t_{N}}$ where 
$X^1,\ldots,X^N$ belong to $\mathcal{F}$. That is to say, vector fields in
$P_{\mathcal{F}}$ are of the form 
\begin{equation}
\label{vftrans}
  \left(
   X^{1}_{t_{1}}\circ\cdots\circ X^{N}_{t_{N}}
   \right)_{\star} X^{0}
\end{equation}
where $X^0,X^1,\ldots,X^N$ belong to $\mathcal{F}$. 
\begin{rem}
\label{remdomP}
Note that a member of $P_{\mathcal{F}}$ is defined on an open set
which is generally a strict subset of $U$, 
whereas members of $\mathcal{F}$ are defined over the whole of $U$, and it
is understood that 
a curve $\gamma:I\to U$, where $I$ is a real interval, will 
be called an integral curve of $Y\in P_{\mathcal{F}}$ only when 
$\gamma(I)$ is included in the domain of definition of $Y$.
\end{rem}
For $x\in U$, we denote by $P_{\mathcal{F}}(x)$ the subspace of $\RR^d$
spanned by all the vectors $Y(x)$, where $Y\in P_{\mathcal{F}}(x)$ 
is defined in a neighborhood of $x$.

Theorem \ref{th-orb} below, which is the central result in this appendix, 
describes the topological nature of the orbits. To interpret
the statement correctly, it is necessary to recall (see for instance
\cite{Spiv79a}) that an \emph{immersed} sub-manifold of a manifold 
is a subset of the latter which is a manifold
in its own right, and is such that the inclusion map is
an immersion. This allows one to naturally identify the tangent space to 
an immersed sub-manifold at a given point with a linear subspace of the tangent
space to the ambient manifold at the same point.
The topology of an immersed sub-manifold is in general finer than the one 
induced by the ambient manifold; when 
these two topologies coincide, the sub-manifold is
called \emph{embedded}.

\begin{thm}[Orbit Theorem, Sussmann \cite{Suss73}]
\label{th-orb}
Let $\mathcal{F}$ be a family of smooth vector fields defined on an open set
$U\subset\RR^{d}$, and $m$ be a point in $U$. If $\mathcal{O}_{\mathcal{F},m}$
 denotes the orbit of $\mathcal{F}$ through $m$, then:
\begin{itemize}
\item [(i)]
Endowed with the orbit topology, $\mathcal{O}_{\mathcal{F},m}$
has a unique differential structure that makes it a smooth
connected immersed sub-manifold of $U$, for which the maps 
\eqref{deforb} are smooth.
\item [(ii)] 
The tangent space to $\mathcal{O}_{\mathcal{F},m}$ at $x\in
\mathcal{O}_{\mathcal{F},m}$ is $P_{\mathcal{F}}(x)$. 
\item [(iii)] 
There exists an open neighborhood $W$ of $m$ \underline{in $U$}, and smooth
local
coordinates $\xi:W\rightarrow\,(-\eta,\eta)^d\subset\RR^d$, with
$\xi(m)=0$, such that
\begin{itemize}
\item [(a)] in 
these coordinates, $W\cap \mathcal{O}_{\mathcal{F},m}$ is a product~:
\begin{equation}
\label{WcapS}
  W\cap \mathcal{O}_{\mathcal{F},m}\ \;=\;\ (-\eta,\eta)^q\times T
\end{equation}
where $\eta>0$, $q$ is the dimension of $\mathcal{O}_{\mathcal{F},m}$, and
$T$ is some subset of $(-\eta,\eta)^{d-q}$ containing the origin.
The orbit topology of $ \mathcal{O}_{\mathcal{F},m}$ induces
on $W\cap \mathcal{O}_{\mathcal{F},m}$ the product 
topology where $(-\eta,\eta)^{q}$ is endowed with
the usual Euclidean topology and $T$ with the discrete topology.
\item[(b)] 
if $\gamma:[t_1,t_2]\to W\cap \mathcal{O}_{\mathcal{F},m}$ is an integral 
curve of a vector field
$Y\in P_{\mathcal{F}}$ (see remark \ref{remdomP}), then
$t\mapsto\xi_{i}(\gamma(t))$, $q+1\leq i\leq d$, are constant mappings,
\item[(c)] 
the tangent space to $\mathcal{O}_{\mathcal{F},m}$ at each point
$p\in W\cap \mathcal{O}_{\mathcal{F},m}$ is spanned by the vector fields
$\partial/\partial\xi_{1},\ldots,\partial/\partial\xi_{q}$,
\item[(d)] 
at any point $p\in W$, the vector fields
$\partial/\partial\xi_{1},\ldots,\partial/\partial\xi_{q}$ belong to the
tangent space to the orbit of $\mathcal{F}$ through $p$.
\end{itemize}
\end{itemize}
\end{thm}
\begin{rem}
\label{rmk-topoprod}
Another description of the product topology in point $(iii)-(a)$ is as follows.
The connected components of $W\cap \mathcal{O}_{\mathcal{F},m}$ are the sets
\begin{equation}
  \label{eq:SWa}
  S_{W,a}\ \;=\;\ (-\eta,\eta)^q\times\{a\}
\end{equation}
for $a\in T$, and the topology on each of
these connected components is the topology induced by the ambient Euclidean
topology. In particular each $S_{W,a}$ is an \emph{embedded} sub-manifold of $U$.
\end{rem}
\begin{proof}[Proof of Theorem \ref{th-orb}]
Assertion $(i)$ is the standard form of the orbit theorem (cf {\it e.g.}
\cite[Chapter 2, Theorem 1]{Jurd97}), while assertion $(ii)$ is a rephrasing
of \cite[Theorem 4.1, point (b)]{Suss73}. Assertion $(iii)$ apparently
cannot be referenced exactly in this form, but we shall deduce it from
the previous ones as follows.

By point $(ii)$, the tangent space to $\mathcal{O}_{\mathcal{F},m}$ 
at $m\in S$ is 
the linear span over $\RR$ of $Y^{1}(m),\ldots,Y^{q}(m)$, where 
$Y^{1},\ldots,Y^{q}$ are $q$ vector fields belonging to $P_{\mathcal{F}}$,
defined on some neighborhood of $m$, and such that $Y^{1}(m),\ldots,Y^{q}(m)$
are linearly independent (recall that $q$ is the dimension of 
$\mathcal{O}_{\mathcal{F},m}$).
Let us write
\[
   Y^{j}\ \;=\;\  \left(
   X^{j,1}_{t_{j,1}}\circ\cdots\circ X^{j,N_{j}}_{t_{j,N_{j}}}
   \right)_{\star} X^{j,0},~~~~1\leq j\leq q,
\]
where $X^{j,k}\in\mathcal{F}$ for $0\leq k\leq N_j$, and where the
$t_{j,k}$'s are real numbers for which the concatenated flow
exists, locally around $m$ (compare \eqref{vftrans}).

Since $Y^{1}(m),\ldots,Y^{q}(m)$ are linearly independent, 
one may
complement them into a basis of $\RR^d$ by adjunction of $d-q$ 
independent vectors that may, without loss of generality, be 
regarded as values at $m$ of $d-q$ smooth vector fields in $U$, say
$Y^{q+1},\ldots,Y^{d}$. Then, the smooth map
\begin{equation}
\label{coordorb}
  L(\xi_1,\ldots,\xi_d)\ \;=\;\ \left(
Y^{1}_{\xi_{1}}\circ\cdots\circ Y^{q}_{\xi_q}
\circ Y^{q+1}_{\xi_{q+1}}\circ\cdots\circ Y^{d}_{\xi_{d}}
   \right)(m)
\end{equation}
defines a diffeomorphism from some
poly-interval $\mathcal{I}_\eta=\{(\xi_1,\ldots,\xi_d)\,,\;|\xi_i|<\eta\}$
onto an open neighborhood $W$ of $m$ in $U$, simply
because the derivative of $L$ is invertible at the origin as
$Y^{1}(m),\ldots,Y^{d}(m)$ are linearly independent by construction.
Let $\xi:W\to\mathcal{I}_\eta$ denote its inverse. 

By the characteristic property of push-forwards, 
we locally have, for $1\leq j\leq q$, that
\begin{equation}
\label{pfY}  Y^{j}_{{\xi_{j}}}=
X^{j,1}_{t_{j,1}}\circ\cdots\circ X^{j,N_{j}}_{t_{j,N_{j}}}
   \circ X^{j,0}_{\xi_{j}} \circ
X^{j,N_{j}}_{-t_{j,N_{j}}}\circ\cdots\circ X^{j,1}_{-t_{j,1}}\;.
\end{equation}
This implies that, in \eqref{coordorb}, 
the images under $L$ of those $d$-tuples sharing a common
value of $\xi_{q+1},\ldots,\xi_d$ all lie in the same
orbit $\mathcal{O}_{\mathcal{F},L(0,\ldots,0,\xi_{q+1},\ldots,\xi_d)}$. In particular, the map 
\begin{displaymath}
\tau_1,\ldots,\tau_q\ \;\mapsto\;\ \left(
Y^{1}_{\tau_1+\xi_{1}}\circ\cdots\circ Y^{q}_{\tau_q+\xi_q}
\circ Y^{q+1}_{\xi_{q+1}}\circ\cdots\circ Y^{d}_{\xi_{d}}
   \right)(m)
\end{displaymath}
is defined $\Pi_{j=1}^q(-\eta-\xi_j,\eta-\xi_j)\to
W\cap\mathcal{O}_{\mathcal{F},L(\xi_{1},\ldots,\xi_d)}$,
and this map is smooth from the Euclidean to the orbit topology
by \eqref{pfY} and point $(i)$.
If we compose it with the immersive injection
$J_W:W\cap\mathcal{O}_{\mathcal{F},L(\xi_{1},\ldots,\xi_d)}\to W$ 
(keeping in mind that $W\cap\mathcal{O}_{\mathcal{F},L(\xi_{1},\ldots,\xi_d)}$
is open in $\mathcal{O}_{\mathcal{F},L(\xi_{1},\ldots,\xi_d)}$ since the orbit 
topology is finer than the Euclidean one),
and if we subsequently apply $\xi$, we get
the affine map
\begin{equation}
\label{afficomp}
\tau_1,\ldots,\tau_q\ \;\mapsto\;\ (
\tau_1+\xi_{1},\cdots,\tau_q+\xi_q,\xi_{q+1},\cdots,\xi_{d}
   ).
\end{equation}
Thus the derivative of \eqref{afficomp} factors through the derivative of
$\xi\circ J_W$ at $ L(\xi_{1},\ldots,\xi_d)$,
which implies $(d)$; from this $(c)$ follows, because $q$ is the
dimension of the orbit through $m$.
If $Y\in P_{\mathcal{F}}$ is defined over an open subset of $W$, 
and if we write in the $\xi$ coordinates 
$Y(\xi)=\sum_{i}a_i(\xi)\partial/\partial\xi_{i}$, 
then, since $Y(\xi)$ is tangent to $\mathcal{O}_{\mathcal{F},\xi}$ by
$(ii)$, we deduce from $(c)$, that the functions $a_{q+1},\ldots,a_{d}$ 
vanish on $\mathcal{O}_{\mathcal{F},m}$, whence $(b)$ holds.

We finally prove $(a)$. Considering \eqref{coordorb} and \eqref{pfY},
a moment's thinking will convince the reader
that $W\cap \mathcal{O}_{\mathcal{F},m}$ consists exactly, 
in the $\xi$ coordinates, of those
$(\xi_1,\ldots,\xi_d)$ such that 
\begin{equation}
\label{descxi}
\left(
Y^{q+1}_{\xi_{q+1}}\circ\cdots\circ Y^{d}_{\xi_{d}}
   \right)(m)\in \mathcal{O}_{\mathcal{F},m},
\end{equation}
which accounts for \eqref{WcapS} where $T$
is the set of $(d-q)$-tuples $(\xi_{q+1},\ldots,\xi_d)$ such that
\eqref{descxi} holds. 
To prove that the orbit topology is the product topology on
$(-\eta,\eta)^q\times T$ where $T$ is discrete, consider a map 
$F$ as in \eqref{deforb}, and pick
$\bar{t}=(\bar{t}_1,\ldots,\bar{t}_N)\in\mathrm{dom}(F)$ such that
$F(\bar{t})\in W$ (hence $F(\bar{t})\in W\cap\mathcal{O}_{\mathcal{F},m}$)~; 
then
$F$ is continuous at $\bar{t}$ for the product topology
because, for $t$ close enough to $\bar{t}$, the values
$\xi_{q+1}(F(t)),\ldots,\xi_d(F(t))$ do not depend on $t$ 
by $(b)$ (moving $t_i$ means following the flow of
a vector field in $P_{\mathcal{F}}$, namely the push-forward of $X^i$
through $X_{t_1}^1\circ\cdots\circ X_{t_{i-1}}^{i-1}$) while
$\xi_{1}(F(t)),\ldots,\xi_q(F(t))$ vary continuously with $t$ according
to  the continuous dependence on time and initial conditions of solutions
to differential equations.
Since this is true for all maps $F$, the orbit topology on 
$W\cap\mathcal{O}_{\mathcal{F},m}$ is finer than the
product topology. To show that it cannot be strictly finer,
it is enough to prove that the orbit topology coincides with the Euclidean
topology on each set  $S_{W,a}$ defined in \eqref{eq:SWa}, a basis
of which consists of the sets $O\times\{a\}$ where $O$ is open in
$(-\eta,\eta)^q$. Being open for the product topology, these sets
are open the orbit topology as well by what precedes and, since
$\mathcal{O}_{\mathcal{F},m}$ is a manifold by $(i)$, each point 
$(y,a)\in O\times\{a\}$ has, in the orbit topology,
a neighborhood $\mathcal{N}_y\subset O\times\{a\}$
which is homeomorphic to an open ball of $\RR^q$ {\it via} some 
coordinate map. When viewed in these coordinates,
the injection $\mathcal{N}_y\to O\times\{a\}$ from the orbit topology to 
the Euclidean topology is a continuous injective map from an open ball in 
$\RR^q$ into $\RR^q$, and therefore it is a homeomorphism onto its image by
invariance of the domain. As $(y,a)$ was arbitrary in 
$O\times\{a\}$, this shows the latter is a union of open sets for the orbit
topology, as desired.
\end{proof}

Consider now the control system~:
\begin{equation}
  \label{sysnl-ter}
  \dot{x}\ \;=\;\ f(x,u),
\end{equation}
with state $x\in\RR^{d}$ and control 
$u\in\RR^{r}$, the function $f$ being smooth on 
$\RR^{d}\times\RR^{r}$. Let $\Omega$ be an open subset of
$\RR^{d}\times\RR^{r}$ and,
following the notation introduced in section \ref{sec-topequiv}, put
$\Omega_{\RR^{d}}$ to denote its projection onto the first factor.
In the proof of Theorem \ref{lem-C1}, we shall be concerned
with the following family of vector fields on $\Omega_{\RR^{d}}$~:
\begin{equation}
  \label{F'gen}
  \mathcal{F}'\ =\ \{\,
\delta f_{\alpha_1,\alpha_2}\,,\ 
\alpha_{1},\alpha_{2}\mbox{ feedbacks on $\Omega$}\,\}\ ,
\end{equation}
where feedbacks on $\Omega$ were introduced in Definition \ref{def-feed}
and the notation $\delta f_{\alpha_1,\alpha_2}$ was fixed in 
\eqref{falpha}, \eqref{deltafalpha}.

Since feedbacks are only required to be \emph{continuous},
$\mathcal{F}'$  is a family of continuous \emph{but
not necessarily differentiable} vector fields on  $\Omega_{\RR^{d}}$ and,
though  the existence of solutions to differential equations with continuous 
right-hand side makes it still possible to define the orbit as the collection
of endpoints of all concatenated integrations like \eqref{deforb}, 
Theorem \ref{th-orb} does not apply in this case.

To overcome this difficulty, we will consider instead of 
$\mathcal{F}'$ the smaller family~:
\begin{equation}
  \label{F''gen}
  \mathcal{F}''\ =\ \{\,
X\in\mathcal{F}'\,,\ X\mbox{ has a flow}\,\}\ ,
\end{equation}
where the sentence ``$X$ has a flow'' means, as in appendix \ref{blode1},
that the Cauchy problem $\dot{x}(t)=X(x(t))$, $x(0)=x_0$, has a unique 
solution,  defined for $|t|<\varepsilon_0$ where $\varepsilon_0$
may depend on $x_0$, whenever $x_0$ lies in the domain of definition of $X$.
Let us consider
the orbit $\mathcal{O}_{\mathcal{F}'',m}$ of
$\mathcal{F}''$ through $m\in \Omega_{\RR^{d}}$, which is still defined as the 
union of images of all maps \eqref{deforb} where $X^j\in\mathcal{F}''$,
the domain of each such map $F$ being again a connected open neighborhood 
$\mathrm{dom}(F)$ of the
origin in $\RR^N$ by repeated application of  Lemma~\ref{lem-ode2}.
As before, we define the orbit
topology on $\mathcal{O}_{\mathcal{F}'',m}$ to be the finest that makes all
the maps \eqref{deforb} continuous,
and since uniqueness of solutions implies continuous dependence on initial
conditions (see Lemma~\ref{lem-ode2}),
the orbit topology is again finer than the Euclidean topology.
{\it A priori}, we know very little about $\mathcal{O}_{\mathcal{F}'',m}$
and its orbit topology as Theorem \ref{th-orb} does not apply.
However, Proposition \ref{prop-FF'} below will establish 
that these notions coincide with those arising from the 
family $\mathcal{F}$ of \emph{smooth} vector fields obtained by setting~:
\begin{equation}
  \label{Fgen}
  \mathcal{F}\ =\ \{\,\delta f_{\alpha_{1},\alpha_{2}}\,,\ 
\alpha_{1},\alpha_{2}\mbox{ \underline{smooth} feedbacks on $\Omega$}\,\}.
\end{equation}

Note that, from the definitions 
\eqref{F'gen}, \eqref{F''gen} and \eqref{Fgen}, we obviously have
\begin{equation}
  \label{Fordre}
  \mathcal{F}\ \subset\ \mathcal{F}''\ \subset\ \mathcal{F}'\ ,
\end{equation}
hence the orbits of these families through a given point obey the same
inclusions.
\begin{rem}
It may of course happen that the family $\mathcal{F}'$ is empty because
$\Omega$ admits no feedback at all.
However, if  $\mathcal{F}'$ is not empty, then $\mathcal{F}$ is not empty 
either by Proposition \ref{prop-feed-lisse}.
\end{rem}

\begin{prop}
\label{prop-FF'} 
Suppose that $f:\RR^{d}\times\RR^{r}\to\RR^d$ is smooth,
and let $\Omega$ be an open subset of $\RR^{d}\times\RR^{r}$.
Let $\mathcal{F}''$ be defined by \eqref{F'gen}-\eqref{F''gen}.

For any $m\in\Omega_{\RR^{d}}$, the orbit $\mathcal{O}_{\mathcal{F}'',m}$ 
of $\mathcal{F}''$ through $m$ coincides with
the orbit through $m$ of the family $\mathcal{F}$ of \emph{smooth} vector
fields defined by \eqref{Fgen}, and the topology of
$\mathcal{O}_{\mathcal{F}'',m}$, as an orbit of
$\mathcal{F}$, coincides with its topology as an orbit of 
$\mathcal{F}''$.
In particular, the conclusions of Theorem \ref{th-orb} hold if we replace
$\mathcal{F}$ by $\mathcal{F}''$ and $U$ by $\Omega_{\RR^d}$.
\end{prop}

\begin{rem}
With a limited amount of extra-work, it is possible to show that the orbits
of $\mathcal{F}'$ also coincide with those of $\mathcal{F}$.
Hence they turn out to be manifolds despite the possible non-uniqueness of 
solutions to the Cauchy problem. However, \eqref{deforb} is no longer
convenient to define the orbit topology in this case
because the maps $F$ may be multiply-valued when $X^j\in\mathcal{F}'$,
and it is simpler to work with the family $\mathcal{F}''$ anyway.
\end{rem}

The proof of the proposition is based on the following lemma.
\begin{lem}
\label{prop-iteres} 
For $m\in\Omega_{\RR^d}$ and $X^1,\ldots,X^N\in\mathcal{F}''$, 
let $F:\mathrm{dom}(F)\to\Omega_{\RR^d}$ be defined by 
\eqref{deforb}. Fix 
$\bar{t}=(\bar{t}_1,\ldots,\bar{t}_N)\in\mathrm{dom}(F)$ 
and set $\overline{m}=F(\bar{t})$.

Then, there is a neighborhood $\mathcal{T}$ of $\bar{t}$ in 
$\mathrm{dom}(F)$, with 
$F(\mathcal{T})\subset \mathcal{O}_{\mathcal{F},\overline{m}}$,
such that $F:\mathcal{T}\to \mathcal{O}_{\mathcal{F},\overline{m}}$ is
continuous from the Euclidean topology to the orbit topology.
\end{lem}
Assuming the lemma for a while, we first prove the proposition.
\begin{proof}[Proof of Proposition \ref{prop-FF'}]
We noticed already from \eqref{Fordre} that the orbit of
$\mathcal{F}''$ through $m$ contains the orbit of $\mathcal{F}$
through $m$. To get the reverse inclusion, consider the  map $F$ 
defined by \eqref{deforb} for some vector fields $X^1,\ldots,X^N$ 
belonging to $\mathcal{F}''$. Then, observe from Lemma~\ref{prop-iteres}
that $F$ takes values in a disjoint union of orbits of $\mathcal{F}$,
and that it is continuous if each orbit in this union is endowed with
the orbit topology. Since $\mathrm{dom}(F)$ is connected, $F$ takes 
values in a single orbit, which can be none but
$\mathcal{O}_{\mathcal{F},m}$. As $F$ was arbitrary, we conclude
that $\mathcal{O}_{\mathcal{F}'',m}\subset\mathcal{O}_{\mathcal{F},m}$ and
therefore the two orbits agree as sets.
Moreover, since each map $F$ was continuous $\mathrm{dom}(F)\to
\mathcal{O}_{\mathcal{F},m}$, the orbit topology of 
$\mathcal{O}_{\mathcal{F}'',m}$ is by definition finer than the
orbit topology of $\mathcal{O}_{\mathcal{F},m}$; but since it is also 
coarser, by definition of the orbit topology on 
$\mathcal{O}_{\mathcal{F},m}$,
because $\mathcal{F}\subset\mathcal{F}''$, the two topologies
in turn agree as desired.
\end{proof}

\begin{proof}[Proof of Lemma~\ref{prop-iteres}]
Theorem \ref{th-orb} applied to the family $\mathcal{F}$, at the point
$\overline{m}=F(\bar{t})$, yields an open
neighborhood $W$ of $\overline{m}$ in $\Omega_{\RR^d}$ and 
smooth local coordinates
$(\xi_1,\ldots,\xi_d):W\to(-\eta,\eta)^d$ satisfying properties 
$(iii)-(a)$ to $(iii)-(d)$ of that theorem.
For $\varepsilon>0$ denote by $\mathcal{T}_{\varepsilon}$ the 
compact poly-interval~:
\begin{displaymath}
  \mathcal{T}_{\varepsilon}\ \;=\;\
  \{t=(t_1,\ldots,t_N)\in\RR^N,\;|t_i-\bar{t}_i|\leq\varepsilon\}\;.
\end{displaymath}
By Lemma~\ref{lem-ode2}, $F$ is
continuous $\mathrm{dom}(F)\to\Omega_{\RR^d}$ and, since
$\mathrm{dom}(F)$ is an open neighborhood of $\overline{t}$
in $\RR^N$, we can pick $\varepsilon>0$ such that  
\begin{displaymath}
  \mathcal{T}_{\varepsilon}\subset\mathrm{dom}(F)
\ \ \ \textrm{and}\ \ \ 
F(\mathcal{T}_{\varepsilon})\subset W\ .
\end{displaymath}
As $X^1,\ldots,X^N$ belong to 
$\mathcal{F}''\subset\mathcal{F}'$, we can write
\begin{displaymath}
  X^{\ell}\ \;=\;\ \delta f_{\alpha^{\ell}_1,\alpha^{\ell}_2}\ ,\ \ 1\leq {\ell}\leq N
\end{displaymath}
for some collection of feedbacks 
$\alpha^{\ell}_1$, $\alpha^{\ell}_2$ on $\Omega$.
From Proposition~\ref{prop-feed-lisse}, there exists for each
$({\ell},l)\in\{1,\ldots,N\}\times\{1,2\}$ a sequence of 
\emph{smooth} feedbacks on $\Omega$, 
say $(\beta^{{\ell},k}_l)_{k\in\NN}$, 
converging to $\alpha^{\ell}_l$ uniformly on $\Omega_{\RR^d}$.
Subsequently, we let $Y^{{\ell},k}$ denote, for $1\leq {\ell}\leq N$ 
and $k\in\NN$, the smooth vector field on $\Omega_{\RR^d}$
\begin{displaymath}
  Y^{{\ell},k}\ \;=\;\ \delta f_{\beta^{{\ell},k}_1,\beta^{{\ell},k}_2}\ .
\end{displaymath}
Clearly $Y^{{\ell},k}\in\mathcal{F}$ and, for each 
${\ell}$, we have that $Y^{{\ell},k}$ converges to
$X^{\ell}$ as $k\to\infty$, uniformly on compact subsets of $\Omega_{\RR^d}$.

 Now, pick $j\in\{1,\ldots,N\}$ and consider a $N$-tuple
$t^{(j)}\in\mathcal{T}_{\varepsilon}$ of the form~:
\begin{displaymath}
  t^{(j)}=(\bar{t}_{1},\ldots,\bar{t}_{j-1},t_{j},\ldots,t_N)\,,~~~~
|t_\ell-\bar{t}_\ell|\leq\varepsilon~\mbox{ for $j\leq\ell\leq N$}.
\end{displaymath}
Let also $\mathbf{1}_j$ designate, for simplicity, the $N$-tuple 
$(0,\ldots,1,\ldots,0)$ with zero entries except for the 
$j$-th one which is $1$. Then, for $|\lambda|\leq\varepsilon$, we have that
\begin{displaymath}
  t^{(j)}+\lambda\mathbf{1}_j=
(\bar{t}_{1},\ldots,\bar{t}_{j-1},\bar{t}_{j}+\lambda,t_{j+1},\ldots,t_N)\ 
\in\mathcal{T}_{\varepsilon},
\end{displaymath}
and a simple computation allows us to rewrite $F(t+\lambda\mathbf{1}_j)$ as~:
\begin{displaymath}
\!\!\!\!\!\!\!\!\!\!\!\!
  F(t^{(j)}+\lambda\mathbf{1}_j) \ =\ 
X^{1}_{\bar{t}_1} \circ\cdots\circ X^{j-1}_{\bar{t}_{j-1}}
\circ X^{j}_{\lambda} \circ
X^{j-1}_{-\bar{t}_{j-1}} \circ\cdots\circ X^{1}_{-\bar{t}_1}
(F(t)).
\end{displaymath}
Let us set
\begin{displaymath}
  \label{Anlamb}\!\!\!\!\!\!
  A_k(\lambda)\;=\;Y^{1,k}_{\bar{t}_1} \circ\cdots\circ Y^{j-1,k}_{\bar{t}_{j-1}}
\circ Y^{j,k}_{\lambda} \circ
Y^{j-1,k}_{-\bar{t}_{j-1}} \circ\cdots\circ Y^{1,k}_{-\bar{t}_1}
(F(t)).
\end{displaymath}
Repeated applications of Lemmas \ref{lem-ode2} and \ref{lem-ode3}
show that, for fixed $j$ and $t^{(j)}$,
the map $\lambda\mapsto A_k(\lambda)$ is well-defined 
$[-\varepsilon,\varepsilon]\to W$
as soon as the integer $k$ is sufficiently large, and moreover that
$A_k(\lambda)$ converges to $F(t^{(j)}+\lambda\mathbf{1}_j)$
as $k\to+\infty$, uniformly with respect to
$\lambda\in[-\varepsilon,\varepsilon]$.
Now, by the characteristic property push forwards,
$\lambda\mapsto A_k(\lambda)$ is an integral curve of the smooth vector field 
\begin{displaymath}
  Z^k\ \;=\;\ \left(
Y^{1,k}_{\bar{t}_1} \circ\cdots\circ Y^{j-1,k}_{\bar{t}_{j-1}}
\right)_\star Y^{j,k}\ ,
\end{displaymath}
which is defined on a neighborhood of 
$\{F(t^{(j)}+\lambda\mathbf{1}_j);~|\lambda|\leq\varepsilon\}$ in $W$.
Since $Z^k\in P_{\mathcal{F}}$
({\it cf} equation \eqref{vftrans}), it follows from
point $(iii)-(b)$ of Theorem \ref{th-orb} that, for $k$ large enough,
\begin{displaymath}
  \xi_i\circ A_{k}(\lambda)\;=\;\xi_i\circ A_{k}(0)\;,\ \ 
\forall\lambda\in[-\varepsilon,\varepsilon],
\;i\in\{q+1,\ldots,d\}.
\end{displaymath}
It is clear from the definition that $A_{k}(0)=F(t^{(j)})$; hence, using the
continuity of $\xi_i$ and taking, in the above equation,
the limit as $k\to+\infty$, we get
\begin{equation}
\label{descont}
\xi_i\circ  F(t^{(j)}+\lambda\mathbf{1}_j)=\xi_i\circ F(t^{(j)}),\ 
\forall\lambda\in[-\varepsilon,\varepsilon],
\;i\in\{q+1,\ldots,d\}.
\end{equation}
Since $\xi_{q+1}\circ F(\bar{t})=\cdots=\xi_{d}\circ F(\bar{t})=0$
by definition of $W$, successive applications of \eqref{descont}
for $j=N,\ldots,1$ lead us to the conclusion that
\begin{equation}
\label{ancont}
\xi_{q+1}\circ F(t)=\cdots=\xi_d\circ F(t)=0,\ \ \forall t\in
\mathcal{T}_{\varepsilon}\ .
\end{equation}
Equation \eqref{ancont} means that, in the $\xi$-coordinates,
$F(\mathcal{T}_{\varepsilon})\subset (-\eta,\eta)^q\times \{0\}$. Hence,
from the local description of the orbits
in \eqref{WcapS} (where $m$ is to be
replaced by $\overline{m}$), we deduce that
$F(\mathcal{T}_{\varepsilon})\subset \mathcal{O}_{\mathcal{F},\overline{m}}$.
Actually, with the notations of \eqref{eq:SWa},
we even get the stronger conclusion that
\begin{displaymath}
  F(\mathcal{T}_{\varepsilon})\ \subset\  S_{W,0}
\end{displaymath}
which achieves the proof of the lemma, with
$\mathcal{T}=\mathcal{T}_{\varepsilon}$, because
the orbit topology on $S_{W,0}$ is the
Euclidean topology by Remark \ref{rmk-topoprod}.
\end{proof}

\subsection*{Acknowledgments}
The authors wish to acknowledge fruitful discussions with Prof. M. Chyba, from
University of Hawaii (USA).
Thanks are also due to Prof. C.T.C. Wall from the University of 
Liverpool for his comments on the open question of section \ref{subsec-diff}.


\end{document}